 \documentstyle[preprint,eqsecnum,aps,prb,amssymb,xy]{revtex}
\xyoption{graph}
\xyoption{poly}
\def\dpst{\displaystyle }
\def\tsum{\mathop{\textstyle \sum }}
\def\tprod{\mathop{\textstyle \prod }}
\def\QATOPD#1#2#3#4{{#3 \atopwithdelims#1#2 #4}}

\tightenlines
\begin{document}
\draft
\title{The multiple sum formulas for 9$j$ and 12$j$ coefficients \\
of SU(2) and $u_q$(2)}
\author{Sigitas Ali\v{s}auskas}
\address{Institute of Theoretical Physics and Astronomy, \\
A.\ Go\v stauto 12, Vilnius 2600, Lithuania}
\date{\today}
\maketitle

\vspace{-1.8cm}
\begin{abstract}
Seven different triple sum formulas for $9j$ coefficients of the quantum
algebra $u_q(2)$ are derived, using for these purposes the usual expansion
of $q$-$9j$ coefficients in terms of $q$-$6j$ coefficients and recent
summation formulas of twisted $q$-factorial series (resembling the very
well-poised basic hypergeometric $_5\phi _4$ series) as 
$q$-generalizations of Dougall's summation formula of the very 
well-poised hypergeometric $_4F_3(-1)$ series. This way for $q=1$ the new 
proof of the known triple sum formula is proposed, as well as six new 
triple sum formulas for $9j$ coefficients of the SU(2) group, in the 
angular momentum theory. The mutual rearrangement possibilities of the 
derived triple sum formulas by means of the Chu--Vandermonde summation 
formulas are considered and applied to derive several versions of double 
sum formulas for the stretched $q$-$9j$ coefficients, which give new 
rearrangement and summation formulas of special Kamp\'e de F\'eriet 
functions and their $q$-generalizations. Several fourfold sum formulas 
[with each sum of the $_5F_4(1)$ or $_5\phi _4$ type] for the $12j$ 
coefficients of the second kind (without braiding) of the SU(2) and 
$u_q(2)$ are proposed, as well as expressions with five sums [of the 
$_4F_3(1)$ and $_3F_2(1)$ or $_4\phi _3$ and $_3\phi _2$ type] for the 
$12j$ coefficients of the first kind (with braiding) instead of the usual 
expansion in terms of $q$-$6j$ coefficients. Stretched and doubly 
stretched $q$-$12j$ coefficients [as triple, double or single sums, 
related to composed or separate hypergeometric $_4F_3(1)$ and $_5F_4(1)$ 
or $_4\phi _3$ and $_5\phi _4$ series, respectively] are considered. 
\end{abstract}
  

\section{INTRODUCTION}

The Wigner $9j$ coefficients arise as the recoupling coefficients of the
four irreducible representations (irreps) of the SU(2) group and play the
important roles in the quantum mechanical angular momentum theory.\cite
{Ed57,JLV60,JB77,BL81a,BL81b} Many applications have also the $12j$
coefficients of both kinds\cite{JLV60,JB77} as the recoupling 
coefficients of the five irreps of the SU(2) group. There are many known 
expressions for $9j$ coefficients as multiple series. Nevertheless, the 
most compact formula (however, which does not represent any symmetry of 
$9j$-symbol) was derived originally by Ali\v sauskas and Jucys\cite{AlJ71} 
as a triple sum series, in frames of resolution of the multiplicity-free 
(semistretched) coupling problem\cite{AlJ71,AlJ69} for the states of 
irreducible representations of the Sp(4) [SO(5)] group restricted to 
SU(2)$\times $SU(2). In Ref.\ \onlinecite{JB77} it was also proved after 
tedious rearrangement of the fourfold sum\cite{JB65} in frames of the 
usual angular momentum [SU(2) representation theory] technique, by means 
of the Chu--Vandermonde summation formulas. Different computational, 
\cite{ZhZ88,S-RRC89,LCh92,Rt93} polynomial,\cite{S-RR88,S-RPV-J96}
rearrangement,\cite{S-RV-J94,V-JPS-R94,PV-J96,S-R98,V-JPS-R98,S-RV-J98} 
specification,\cite{W73,S-RR89} and other \cite{N89c} aspects of these 
multiple sum series were considered. Their analytical continuation was 
also adapted\cite{Al73} for the isoscalar factors of the Clebsch--Gordan 
(CG) coefficients of the Lorentz or $SL(2,C)$ group.

Recently Rosengren\cite{R98,R99t,R99} proposed two new proofs of the 
triple sum formula for $9j$ coefficients of SU(1,1). His first 
proof\cite{R98,R99t} was based on the use of the explicit coupling 
kernels in the $\frak {su}$(1,1) algebra, rather then in the 
$\frak {su}$(2) algebra, when in the second case\cite{R99t,R99} the 
usual expansion\cite{Ed57,JLV60,JB77,BL81a} of $9j$ coefficients of 
$\frak {su}$(1,1) in terms of $6j$ coefficients was rearranged using 
the appropriate expressions for the Racah coefficients in terms of the 
balanced hypergeometric $_4F_3(1)$ series and Dougall's summation 
formula \cite{Sl66} of the very well-poised $_4F_3(-1)$ series [which, 
in other words, corresponds to the factorial sum weighted with 
factor $(2j+1)$].

For the quantum algebra $u_q(2)$, the expansion of the $q$-$9j$
coefficients in terms of $q$-$6j$ coefficients was generalized by Nomura 
\cite{N89q,N89bq,N90,N91} and Smirnov et al\cite{STK90} and extended to 
$q$-$3nj$ coefficients (particularly, of the first and the second kind) 
by Nomura,\cite{N89bq,N90} who discussed their role in frames of the 
Yang--Baxter Relations. The corresponding summation formula of the 
twisted $q$-factorial series [generalizing Dougall's summation formula 
and resembling (but not equivalent with) the very well-poised basic 
hypergeometric $_5\phi _4$ series, depending on 3 parameters] needed for 
our purpose was derived by Ali\v sauskas\cite{Al97} and the twisted very 
well-poised $q$-factorial series, resembling the basic hypergeometric 
$_7\phi _6$ series (depending on 5 parameters) appear in a new approach 
\cite{A-NS96} to the Clebsch--Gordan coefficients of $u_q(2)$. 
In the $u_q(3)$ context, Ali\v sauskas\cite{Al97} also used the summation 
formula of the $q$-factorial series depending on 4 parameters which 
correspond to Dougall's summation formula of the very well-poised 
hypergeometric $_5F_4(1)$ or basic hypergeometric $_6\phi _5$ 
series.\cite{GR90}

The main purpose of the present paper is to derive the all independent 
expressions with the triple sums for the $q$-$9j$ and the usual $9j$ 
coefficients, as well as to rearrange expressions for the $q$-$12j$ and 
the usual $12j$ coefficients of the both kinds into more convenient 
forms, with minimal number of sums, or at least eliminating the 
cumbersome factorial sums weighted with factors $[2j+1]$ or $(2j+1)$ 
from the compositions of the $q$-$6j$ or usual $6j$ coefficients expanded 
in different forms. 

In Section II, the appropriate expressions for the $6j$ coefficients of 
SU(2) and $u_q(2)$ are presented, as well as the rearranged expansions of 
$9j$ coefficients of $u_q(2)$ in terms of $6j$ coefficients, allowing us
to generalize Rosengren's \cite{R99} approach with help of the new 
summation formulas, weighted with factor $[2j+1]$ (see Appendix A). 
In Section III, seven different triple sum formulas for $9j$ coefficients 
of $u_q(2)$ are derived, their summation intervals and other properties 
are compared. For $q=1$ they turn either to known expression of 
Ali\v sauskas and Jucys,\cite{AlJ71} or to six new triple sum formulas 
for $9j$ coefficients of SU(2). In Section IV, the mutual rearrangement 
possibilities of new expressions by means of the Chu--Vandermonde 
summation formulas (see Appendix B) are considered. Particularly, several 
versions (of different classes) of the double sum formulas for the 
stretched $9j$ coefficients of $u_q(2)$ and SU(2) are derived, which 
enable to get new relations and summation formulas (presented in Appendix 
C) for special Kamp\'e de F\'eriet functions \cite{K-F21} and their 
$q$-generalizations (cf.\ Refs.\ \onlinecite{V-JPS-R94,PV-J96}).

Section V is devoted to rearrangement of the usual expansion formula (in 
terms of $6j$ and $q$-$6j$ coefficients) of $12j$ and $q$-$12j$ 
coefficients of the second kind\cite{JLV60,JB77,N90} (i.e., without 
braiding\cite{N89bq}) into the fourfold sums, using Dougall's summation 
formula\cite{Sl66} of the very well-poised $_5F_4(1)$ series, depending 
on 4 parameters. Also, specific stretched and doubly stretched $12j$ 
and $q$-$12j$ coefficients of the second kind are studied. In Sec.~VI, 
expressions of $12j$ and $q$-$12j$ coefficients of the first 
kind\cite{JLV60,JB77} (with braiding\cite{N90}) in terms of $q$-$6j$ 
coefficients are rearranged using the transformation formula\cite{Al97} 
of the very well-poised $_6F_5(-1)$ series or $q$-factorial sums 
(depending on 5 parameters and weighted with factors $(2j+1)$ or 
$[2j+1]$), resembling the very well-poised basic hypergeometric 
$_7\phi _6$ series. Variety of the stretched and doubly stretched 
$12j$ and $q$-$12j$ coefficients of the first kind are also considered.

\section{PRELIMINARIES}

\subsection{Expressions for the $6j$ coefficients of SU(2) and $u_q(2)$}

The appropriate for our purpose expressions for the $6j$ (Racah) 
coefficients of SU(2) (with the different or coinciding signs of 
summation parameters in 3 numerator $q$-factorial arguments) were derived 
originally by Bandzaitis {\it et al}\cite{BZhMJ64} (see also Refs.\ 
\onlinecite{JB77,JB65}), when Smirnov {\it et al}\cite{STK91b,AsST96} 
rederived them for the Racah coefficients of $u_q(2)$. These two 
expressions, each in two different versions, may be written as follows:
\begin{mathletters}
\begin{eqnarray}
\left\{ \begin{array}{lll}
a & b & e \\ 
d & c & f
\end{array} \right\} _{\!\!q} &=&
\frac{\nabla [acf]\nabla [dbf]}{\nabla [abe]\nabla [dce]}\sum_{z}
\frac{(-1)^{a+b+c+d+z}[c+f-a+z]!}{[z]![a+c-f-z]![b+d-f-z]!}  \nonumber \\
&&\times \frac{[b+f-d+z]![a+d+e-f-z]!}{[e+f-a-d+z]![2f+z+1]!}  
\label{f6jb} \\
&=&\frac{\nabla [acf]\nabla [dbf]}{\nabla [abe]\nabla [dce]}\sum_{z}
\frac{(-1)^{b+c+e+f+z}[c+d-e+z]!}{[z]![c-d+e-z]![b+e-a-z]!}  \nonumber \\
&&\times \frac{[a+b-e+z]![2e-z]!}{[a+d-e-f+z]![a+d-e+f+z+1]!}  
\label{f6jbb}
\end{eqnarray}
\end{mathletters}
and
\begin{mathletters}
\begin{eqnarray}
\left\{ \begin{array}{lll}
a & b & e \\ 
d & c & f
\end{array} \right\} _{\!\!q}^{\!(\prime )} &=&
\frac{\nabla [eab]\nabla [fbd]}{\nabla [ecd]\nabla [fac]}\sum_{z}
\frac{(-1)^{b+c+e+f+z}[2b-z]!}{[z]![b+e-a-z]![b+f-d-z]!}  \nonumber \\
&&\times \frac{[b+e+f-c-z]![b+c+e+f-z+1]!}{[a+b+e-z+1]![b+d+f-z+1]!},
\label{f6jc} \\
&=&\frac{\nabla [eab]\nabla [fbd]}{\nabla [ecd]\nabla [fac]}\sum_{z}
\frac{(-1)^{b+e-a+z}[a+b-e+z]!}{[z]![b+e-a-z]![a-e+f-d+z]!}  \nonumber \\
&&\times \frac{[a-c+f+z]![a+c+f+z+1]!}{[2a+z+1]![a+d-e+f+z+1]!},
\label{f6jcb}
\end{eqnarray}
\end{mathletters}
where $\nabla [abc]$ is asymmetric triangle coefficient, 
\begin{equation}
\nabla [abc]=\left( \frac{[a+b-c]![a-b+c]![a+b+c+1]!}{[b+c-a]!}\right)
^{1/2},  \label{nabla}
\end{equation}
and Eqs.\ (\ref{f6jc})--(\ref{f6jcb}) are less effective than 
(\ref{f6jb}) and (\ref{f6jbb}). Here and in what follows $[x],\;[x]!$, 
$(\alpha |q)_n$, and $\QATOPD[ ] {n}{r}_{q}$ are, respectively, the 
$q$-numbers, $q$-factorials, $q$-Pochhammer symbols, and $q$-binomial 
coefficients,
\begin{eqnarray}
&&[x]=(q^x-q^{-x})/(q-q^{-1}),\quad [x]!=[x][x-1]...[2][1],  \label{qnf} \\
&&(\alpha |q)_{n}=\prod_{k=0}^{n-1}[\alpha +k],\quad [1]!=[0]!=(\alpha
|q)_{0}=1,  \label{qf} \\
&&\QATOPD[ ] {n}{r}_{q}=\frac{[n]!q^{n(n-r)}}{[r]![n-r]!}    \label{qbc} 
\end{eqnarray}
where (\ref{qnf}) and (\ref{qf}) are invariant under substitution 
$q\leftrightarrow q^{-1}$ and turn into usual integers $x$, factorials 
$x!$ and binomial coefficients $\QATOPD( ) {n}{r}$ for $q=1$.

We see that each parameter $b,c$, or $e$ appears only twice in the 
factorial arguments under the summation sign in (\ref{f6jb}), as well as 
parameters $b,c$, or $f$ in (\ref{f6jbb}) [which is obtained after some 
shift of summation parameter in (\ref{f6jb})]. Similarly each parameter 
$a,c$, or $d$ appears only twice in the factorial arguments under the 
summation sign in (\ref{f6jc}), as well as parameters $b,c$, or $d$ in 
(\ref{f6jcb}). Otherwise, each parameter $a$ or $d$ appears four times 
in the factorial arguments under the summation sign in (\ref{f6jb}) and 
(\ref{f6jbb}), as well as parameters $e$ or $f$ in (\ref{f6jc}) and 
(\ref{f6jcb}) and all the parameters in the most symmetric 
\cite{Ed57,JLV60} (Racah\cite{R42}) and the remaining expressions for 
$6j$ and $q$-$6j$ coefficients, \cite{JB77,JB65,STK91b,AsST96} which 
include only usual symmetric triangle coefficients $\Delta [abc]$ in the 
numerator and denominator before the summation sign. Note that some 
triangular conditions restrict the summation intervals in (\ref{f6jb}) 
and (\ref{f6jc}), or they are represented by definite differences of 
factorial arguments in numerator and denominator, for example, 
$(c+f-a+z)-(e+f-a-d+z)\geq 0$ in (\ref{f6jb}), or 
$(b+d+f-z+1)-(b+e+f-c-z)-1\geq 0$ in (\ref{f6jc}).

It should be also noted, that only the expressions presented above 
(\ref{f6jb}) and (\ref{f6jc}) are correlated with the Racah polynomials 
as introduced by Askey and Wilson, see Ref.\ \onlinecite{GR90}. In 
contrast, the most symmetric and the remaining expressions for the $6j$ 
and $q$-$6j$ coefficients turn into the Racah polynomials only after some 
Whipple (Bailey) or Sears transform \cite{GR90} of the balanced 
hypergeometric $_4F_3(1)$ or $_4\phi _3$ series are used.

\subsection{Rearrangement of expansions for $q$-$9j$ coefficients}

We use here the definition of the $q$-$9j$ coefficients of $u_q(2)$ 
introduced by Nomura\cite{N89q,N91} in contrast with definition by 
Smirnov {\it et al} \cite{STK90}, when the substitution $q\rightarrow 
q^{-1}$ is necessary. These coefficients ($q$-$9j$ symbols) may be 
extracted from the recoupling-braiding coefficients of the states of four 
irreps and are invariant under even permutations of their rows or columns 
and under transposition of $3\times 3$ array (interchange of their rows 
and columns),
\begin{equation}
\left\{ \begin{array}{lll}
a & b & e \\ 
c & d & f \\ 
h & k & g
\end{array} \right\} _{\!\!q}=\left\{ \begin{array}{lll}
e & a & b \\ 
f & c & d \\ 
g & h & k 
\end{array} \right\} _{\!\!q}=\left\{ \begin{array}{lll}
c & d & f \\ 
h & k & g \\
a & b & e
\end{array} \right\} _{\!\!q}=\left\{ \begin{array}{lll}
a & c & h \\ 
b & d & k \\ 
e & f & g
\end{array} \right\} _{\!\!q}={\rm etc.}  \label{symev}
\end{equation}
Taking into account the braiding, in the case of odd permutations of 
their rows or columns, the $q$-$9j$ coefficients obey \cite{N89q,N91}
\begin{equation}
\left\{ \begin{array}{lll}
a & b & e \\ 
c & d & f \\ 
h & k & g
\end{array} \right\} _{\!\!q}=A\left\{ \begin{array}{lll}
a & e & b \\ 
c & f & d \\ 
h & g & k 
\end{array} \right\} _{\!\!q^{-1}}=A\left\{ \begin{array}{lll}
c & d & f \\ 
a & b & e \\
h & k & g
\end{array} \right\} _{\!\!q^{-1}}={\rm etc.},  \label{symod}
\end{equation}
where
\[
A=(-1)^{a+b+c+d+e+f+g+h+k}q^{Z_{deh}+Z_{bcg}+Z_{afk}}
\]
and
\[
Z_{deh}=-d(d+1)-e(e+1)-h(h+1).
\]

Let us consider some different versions of expansions 
\cite{N89q,N91,STK90} of the $q$-$9j$ coefficients of $u_q(2)$, written 
after applying some symmetries of $q$-$6j$ coefficients,
\begin{mathletters}
\begin{eqnarray}
\left\{ \begin{array}{lll}
a & b & e \\ 
c & d & f \\ 
h & k & g
\end{array} \right\} _{\!\!q}&=&\sum_{j}(-1)^{2j}q^{Z_{deh}-j(j+1)}[2j+1]
\left\{ \begin{array}{lll}
a & c & h \\ k & g & j
\end{array} \right\} _{\!\!q}\left\{ \begin{array}{lll}
k & j & c \\ f & d & b
\end{array} \right\} _{\!\!q}\left\{ \begin{array}{lll}
a & g & j \\ f & b & e
\end{array} \right\} _{\!\!q}  \label{vr6a} \\
&=&\sum_{j}(-1)^{2j}q^{Z_{deh}-j(j+1)}[2j+1]
\left\{ \begin{array}{lll}
k & g & h \\ a & c & j
\end{array} \right\} _{\!\!q}\left\{ \begin{array}{lll}
j & b & f \\ d & c & k
\end{array} \right\} _{\!\!q}^{\!\!(\prime )}\left\{ \begin{array}{lll}
f & b & j \\ a & g & e
\end{array} \right\} _{\!\!q}  \label{vr6b} \\
&=&\sum_{j}(-1)^{2j}q^{Z_{deh}-j(j+1)}[2j+1]
\left\{ \begin{array}{lll}
h & c & a \\ j & g & k
\end{array} \right\} _{\!\!q}^{\!\!(\prime )}\left\{ \begin{array}{lll}
k & j & c \\ f & d & b
\end{array} \right\} _{\!\!q}\left\{ \begin{array}{lll}
j & g & a \\ e & b & f
\end{array} \right\} _{\!\!q}^{\!\!(\prime )}  \label{vr6c} \\
&=&\sum_{j}(-1)^{2j}q^{Z_{deh}-j(j+1))}[2j+1]
\left\{ \begin{array}{lll}
a & j & g \\ k & h & c
\end{array} \right\} _{\!\!q}\left\{ \begin{array}{lll}
j & b & f \\ d & c & k
\end{array} \right\} _{\!\!q}^{\!(\prime )}\left\{ \begin{array}{lll}
f & j & b \\ a & e & g
\end{array} \right\} _{\!\!q}  \label{vr6d} \\
&=&\sum_{j}(-1)^{2j}q^{Z_{deh}-j(j+1)}[2j+1]
\left\{ \begin{array}{lll}
a & c & h \\ k & g & j
\end{array} \right\} _{\!\!q}\left\{ \begin{array}{lll}
k & j & c \\ f & d & b
\end{array} \right\} _{\!\!q}\left\{ \begin{array}{lll}
b & e & a \\ g & j & f
\end{array} \right\} _{\!\!q}^{\!\!(\prime )}  \label{vr6e} \\
&=&\sum_{j}(-1)^{2j}q^{Z_{deh}-j(j+1)}[2j+1]
\left\{ \begin{array}{lll}
g & j & a \\ c & h & k
\end{array} \right\} _{\!\!q}^{\!\!(\prime )}\left\{ \begin{array}{lll}
j & b & f \\ d & c & k
\end{array} \right\} _{\!\!q}^{\!\!(\prime )}\left\{ \begin{array}{lll}
a & g & j \\ f & b & e
\end{array} \right\} _{\!\!q}  \label{vrx6a} \\
&=&\sum_{j}(-1)^{2j}q^{Z_{deh}-j(j+1)}[2j+1]
\left\{ \begin{array}{lll}
g & j & a \\ c & h & k
\end{array} \right\} _{\!\!q}^{\!\!(\prime )}\left\{ \begin{array}{lll}
j & b & f \\ d & c & k
\end{array} \right\} _{\!\!q}^{\!\!(\prime )}\left\{ \begin{array}{lll}
b & e & a \\ g & j & f
\end{array} \right\} _{\!\!q}^{\!\!(\prime )}  \label{vrx6b}
\end{eqnarray}
\end{mathletters}
where the summation parameters $j$ are restricted by the triangular 
conditions, 
\begin{equation}
\max (|a-g|,|f-b|,|k-c|)\leq j\leq \min (a+g,b+f,c+k).  \label{limun}
\end{equation}
When we use expressions (\ref{f6jb}) or (\ref{f6jbb}) for nonprimed 
$q$-$6j$ coefficients and expressions (\ref{f6jc}) or (\ref{f6jcb}) for 
primed $q$-$6j$ coefficients, the asymmetric triangle coefficients 
depending on the summation parameter $j$ are distributed in their 
numerators and denominators in expansions (\ref{vr6a}), (\ref{vr6b}), 
and (\ref{vr6e})--(\ref{vrx6b}) as follows:
\begin{mathletters}
\begin{equation}
\frac{\nabla [agj]\nabla [kcj]}{1}\times \frac{\nabla [fbj]}{\nabla [kcj]}
\times \frac{1}{\nabla [agj]\nabla [fbj]},  \label{dndabe}
\end{equation}
when in expansions (\ref{vr6c}) and (\ref{vr6d}) as follows: 
\begin{equation}
\frac{\nabla [kcj]}{\nabla [agj]}\times \frac{\nabla [fbj]}{\nabla [kcj]}%
\times \frac{\nabla [agj]}{\nabla [fbj]}  \label{dndcd}
\end{equation}
\end{mathletters}
Particularly, they cancel if we express all but the first $q$-$6j$ 
coefficients in (\ref{vr6a}) by means of (\ref{f6jb}), as well as the 
first and the last $q$-$6j$ coefficients in (\ref{vr6b}) and 
(\ref{vr6d}), the second $q$-$6j$ coefficient in (\ref{vr6c}), and the 
first two $q$-$6j$ coefficients in (\ref{vr6e}), when the remaining 
(primed) $q$-$6j$ coefficients $\QATOPD\{ \} {\cdot \ \cdot \ \cdot }{%
\cdot \ \cdot \ \cdot }_q^{\!(\prime )}$ in (\ref{vr6b})--(\ref{vr6e}) 
are expressed by means of (\ref{f6jc}) and the first $q$-$6j$ coefficient 
of (\ref{vr6a}) by means of (\ref{f6jbb}). It is expedient to use the 
inverse order of summation [with substituted by $z\to a+c-f-z$ parameters 
in (\ref{f6jb})] for the second $q$-$6j$ coefficients in (\ref{vr6a}), 
(\ref{vr6c}), (\ref{vr6e}), and the first and last $q$-$6j$ 
coefficients in (\ref{vr6d}), with $j$ appearing in the upper middle 
position of the corresponding $6j$-symbol.

Now the summation formulas (\ref{smDqb}) or (\ref{smDqa}) of the 
twisted very well-poised $q$-factorial series\cite{Al97} (see Appendix A) 
may be used in (\ref{vr6a}) or in (\ref{vr6b}) and (\ref{vr6c}), 
respectively, if the summation parameters $j$ are restricted naturally by 
the non-negative integer values of the denominator factorial arguments,
\begin{mathletters}
\begin{eqnarray}
&&\max (a-g,|f-b|,k-c) \leq j\leq \min (a+g,c+k),  \label{lima} \\
&&\max (f-b,a-g) \leq j\leq \min (a+g,b+f,c+k),  \label{limb} \\
&&\max (b-f,k-c) \leq j\leq \min (a+g,c+k),  \label{limc}
\end{eqnarray}
respectively. In (\ref{vr6d}) and (\ref{vr6e}), parameters $j$ are,
respectively, restricted by the natural limits,
\begin{eqnarray}
&&\max (|a-g|,f-b,c-k)\leq j\leq b+f,  \label{limd} \\
&&\max (k-c,b-f)\leq j\leq \min (a+g,c+k).  \label{lime}
\end{eqnarray}
\end{mathletters}
In these two last cases, the summation formulas (\ref{smDqb}) and 
(\ref{smDqc}), respectively, may be used.

However, the formal summation intervals (\ref{lima}--\ref{lime}) may 
exceed the interval (\ref{limun}), determined by triangular conditions. 
Of course, separate $q$-$6j$ coefficients with spoiled triangular conditions 
in (\ref{vr6a}--\ref{vr6e}) vanish, but such vanishing is not evident 
for the corresponding pure $q$-factorial sums of the type (\ref{f6jb}) or 
(\ref{f6jc}). We need to consider each case separately, for example, 
when $j=b+f+1,b+f+2,...$ [i.e., for $b+f<j\leq \min (a+g,c+k)$], or 
$\max (g-a,c-k)>j\geq \max (a-g,|f-b|,k-c)$, the second or the first sum 
of the type (\ref{f6jb}) or (\ref{f6jbb}) in expansion of (\ref{vr6a}) 
turns into 0, in accordance with Karlsson's summation 
formula,\cite{K71,AlJJ72} or its $q$-version,\cite{GR90,AlS94}
\begin{equation}
\sum_{s}\frac{(-1)^sq^{(n-m-1)s}}{[s]!\ [n-s]!}\prod_{j=1}^m[A_{j}-s]=
\delta _{m,n}q^{-n(n+1)/2+\sum_{j=1}^{m}A_{j}}  \label{smCq}
\end{equation}
where $m\leq n$ are integers [cf.\ applications of (\ref{smCq}) for the
multiplicity-free isoscalar factors\cite{AlJJ72,AlS94} of SU($n$) and 
$u_q(n)$]. [Note that the third factorial sum in expansion of 
(\ref{vr6a}) may be nonvanishing in spite of spoiling of the triangular 
conditions.]

\section{NEW EXPRESSIONS FOR 9$\lowercase{j}$ COEFFICIENTS OF SU(2) 
AND \lowercase{$u_{q}$}(2)}

\subsection{Expressions with the full triangle restrictions of 
summation intervals}

Hence, using the expansions (\ref{vr6a})--(\ref{vr6e}), alternative 
expressions for the $q$-$6j$ coefficients, and summation formulas 
(\ref{smDqa})--(\ref{smDqc}), at first we obtained five different 
expressions for the $q$-$9j$ coefficients,
\begin{mathletters} 
\begin{eqnarray}
\left\{ 
\begin{array}{lll}
a & b & e \\ 
c & d & f \\ 
h & k & g
\end{array}
\right\} _{\!\!q} &=&(-1)^{c+h-a}
\frac{\nabla [abe]\nabla [feg]\nabla [kbd]}{\nabla [ach]\nabla [fcd]
\nabla [kgh]}  \nonumber \\
&&\times q^{(f+h-e-k)(a+d-e+k+1)-(a-e+f)(a-e+f+1)+Z_{deh}} 
\nonumber \\
&&\times \sum_{z_1,z_2,z_3}\frac{(-1)^{z_1+z_2+z_3}[g-h+k+z_1]!
[a+c-h+z_1]!}{[z_1]![g+h-k-z_1]![c-a+h-z_1]!}  \nonumber \\
&&\times \frac{[2h-z_1]![2d-z_2]![c-d+f+z_2]!}{[z_2]![d-b+k-z_2]!
[c+d-f-z_2]![b+d+k-z_2+1]!}  \nonumber \\
&&\times \frac{[b+e-a+z_3]![e-f+g+z_3]!q^{-z_1(a+d-e+k-z_2-z_3+1)}%
}{[z_3]![a+b-e-z_3]![f+g-e-z_3]![2e+z_3+1]!}  \nonumber \\
&&\times \frac{q^{z_2(e-f-h+k+z_3)+z_3(a-d+f-h)}[a+d-e+k-z_2-z_3]!%
}{[a-d+f-h+z_1+z_2]![e-f-h+k+z_1+z_3]!}  \label{trsa} \\
&=&(-1)^{e-f-h+k}\frac{\nabla [abe]\nabla [feg]\nabla [kbd]}{\nabla
[ach]\nabla [fcd]\nabla [kgh]}  \nonumber \\
&&\times q^{(b+e-a)(e-f-h+k)-(a-e+f+1)(a-e+f)+Z_{deh}} 
\nonumber \\
&&\times \sum_{z_1,z_2,z_3}\frac{(-1)^{z_1+z_2}[a+c-h+z_1]![g-h+k+z_1]!}{%
[z_1]![c+h-a-z_1]![g+h-k-z_1]![z_2]!}  \nonumber \\
&&\times \frac{[2h-z_1]![2b-z_2]![b-c+f+k-z_2]!}{%
[b-d+k-z_2]![b+d+k-z_2+1]![z_3]![a+b-e-z_3]!}  \nonumber \\
&&\times \frac{[b+c+f+k-z_2+1]![b+e-a+z_3]![e-f+g+z_3]!}{[f+g-e-z_3]!
[2e+z_3+1]![e+k-f-h+z_1+z_3]!}  \nonumber \\
&&\times \frac{q^{z_1(b+e-a-z_2+z_3)+z_3(a+b+f-h+k-z_2+1)-z_2(e+k-f-h)}}{%
[b+e-a-z_2+z_3]![a+b+f-h+k+z_1-z_2+1]!}  \label{trsb} \\
&=&(-1)^{c-d+f}\frac{\nabla [ach]\nabla [feg]\nabla [kbd]}{\nabla
[abe]\nabla [fcd]\nabla [kgh]}  \nonumber \\
&&\times q^{(c-d+f)(a+g)-(d-f+k)(c+k+1)+Z_{deh}}  \nonumber \\
&&\times \sum_{z_1,z_2,z_3}\frac{(-1)^{z_1+z_3}[a+c-g+k-z_1]![a+c+g+k
-z_1+1]!}{[z_1]![a+c-h-z_1]![a+c+h-z_1+1]![z_2]!}  \nonumber \\
&&\times \frac{[2c-z_1]![2d-z_2]![c-d+f+z_2]!}{[d-b+k-z_2]![c+d-f-z_2]!
[b+d+k-z_2+1]![z_3]!}  \nonumber \\
&&\times \frac{[a-b+f+g-z_3]![a+b+f+g-z_3+1]![2g-z_3]!}{[f+g-e-z_3]!
[e+f+g-z_3+1]![c-d+f-z_1+z_2]!}  \nonumber \\
&&\times \frac{q^{-z_1(a-d+f+g-k+z_2-z_3)+z_2(a+c+g+k-z_3+1)-z_3(c-d+f)}%
}{[a-d+f+g-k+z_2-z_3]![a+c+g+k-z_1-z_3+1]!}  \label{trsc} \\
&=&(-1)^{b-a+f-g}\frac{\nabla [ach]\nabla [feg]\nabla [kbd]}{\nabla
[abe]\nabla [fcd]\nabla [kgh]}  \nonumber \\
&&\times q^{-(a+h-k+1)(a+b-e)-(e+f-a)(b+f+1)+Z_{deh}}  \nonumber \\
&&\times \sum_{z_1,z_2,z_3}\frac{(-1)^{z_1+z_2+z_3}[2h-z_1]!
[g-h+k+z_1]!}{[z_1]![a-c+h-z_1]![g+h-k-z_1]![a+c+h-z_1+1]!}  \nonumber \\
&&\times \frac{[2b-z_2]![b-c+f+k-z_2]![b+c+f+k-z_2+1]!}{[z_2]!
[b-d+k-z_2]![b+d+k-z_2+1]!}  \nonumber \\
&&\times \frac{[2e-z_3]![a+b-e+z_3]!q^{z_1(a+b-e-z_2+z_3)+z_2(e+f+h-k
-z_3+1)}}{[z_3]![e+f-g-z_3]![b-a+e-z_3]![e+f+g-z_3+1]!}  \nonumber \\
&&\times \frac{q^{z_3(b-a+f+h+k)}[e+f-k+h-z_1-z_3]!}{[b+f-h+k-a+z_1-z_2]!
[a+b-e-z_2+z_3]!}  \label{trsd} \\
&=&(-1)^{c-b+e-g+k}\frac{\nabla [abe]\nabla [feg]\nabla [kbd]}{\nabla
[ach]\nabla [fcd]\nabla [khg]}  \nonumber \\
&&\times q^{Z_{deh}-(a-b-h+f+k)(a+b+e+2)+(f-b)(f-b+1)+(d-b+k)(e+f+h-k+1)} 
\nonumber \\
&&\times \sum_{z_1,z_2,z_3}\frac{[k+g-h+z_1]![a+c-h+z_1]!
[2h-z_1]!}{[z_1]![g+h-k-z_1]![c-a+h-z_1]!}  \nonumber \\
&&\times \frac{(-1)^{z_2+z_3}[2d-z_2]![c-d+f+z_2]![2e-z_3]!}{[z_2]!
[d-b+k-z_2]![c+d-f-z_2]![b+d+k-z_2+1]!}  \nonumber \\
&&\times \frac{q^{z_3(a-d+f-h+z_1+z_2)-z_2(e+f+h-k-z_1+1)-z_1(a+d+e+k+2)}%
}{[z_3]![a-b+e-z_3]![e+f-g-z_3]![e+f+g-z_3+1]!}  \nonumber \\
&&\times \frac{[a+d+e+k-z_2-z_3+1]![e+f+h-k-z_1-z_3]!}{[a+b+e-z_3+1]!
[a-d+f-h+z_1+z_2]!}.  \label{trse}
\end{eqnarray}
\end{mathletters}
Expression (\ref{trsa}) for $q=1$ is equivalent to the known triple sum
formula\cite{AlJ71} for the $9j$ coefficients of SU(2). The
numerator-denominator distributions of factorials, depending on the 
summation parameters $z_1,z_2,z_3$, are different in all expressions 
(\ref{trsa})--(\ref{trse}). All the terms in the last sum of 
(\ref{trsb}), in the second sum of (\ref{trsc}), and in the first sum of 
(\ref{trse}) are of the same sign. The separate sums correspond to the 
finite basic hypergeometric series, 
\begin{equation}
_{p+1}F_{p}\!\left[ \begin{array}{c}
\alpha _1,\alpha _2,...,\alpha _{p+1} \\ 
\beta _1,...,\beta _p
\end{array};q,x\right] =\sum_k\frac{(\alpha _1|q)_k(\alpha _2|q)_k
\cdot \cdot \cdot (\alpha _{p+1}|q)_k}{(\beta _1|q)_k\cdot \cdot \cdot 
(\beta _p|q)_k(1|q)_k}x^k,  \label{fbhs}
\end{equation}
with $p=3,\;x=q^{\pm (c+1)},\;c=\sum_{i=1}^{p+1}\alpha _i-\sum_{j=1}^p
\beta _j$, as defined (with a minor correction) by \'{A}lvarez-Nodarse 
and Smirnov \cite{A-NS96}, instead of the standard basic hypergeometric 
functions $_{p+1}\phi _p$ (see Gasper and Rahman \cite{GR90}). Parameters 
$c=-1$ and $x=1$ for the balanced basic hypergeometric series, which 
appear in expressions for $q$-$6j$ coefficients.\cite{STK91b,AsST96}

The intervals for summation parameters $z_i$ ($i=1,2,3$) are mainly 
restricted by six [in (\ref{trsa}) and (\ref{trse})], five [in 
(\ref{trsb}) and (\ref{trsd})], or four [in (\ref{trsc})] triangle linear 
combinations of the type $a+b-c$, respectively. As result of their 
vanishing we may write 23 different (independent) expressions as double 
sums for the stretched $9j$ coefficients as compositions of 
$_4F_3[\cdot \cdot \cdot ;q,x]$ and $_3F_2[\cdot \cdot \cdot ;q,x]$ 
series [with latter in the 10 cases corresponding to the CG coefficients 
of $u_q(2)$]. Although Minton's summation formulas (\ref{sMint}) or 
(\ref{sMinta}) (see Ref.\ \onlinecite{GR90}) may be used (12 times) for 
separate alternating sums in (\ref{trsa})--(\ref{trse}), satisfying 
special conditions, in the case of some stretched triangles [e.g., for 
$d=b+k$ or $e=b-a$ in Eq.\ (\ref{trsa})], but the expressions obtained 
are equivalent to some derived previously (although using the different 
triple sum expressions). In contrast with Eq.\ (32.13) of Ref.\ 
\onlinecite{JB77} and its $q$-generalization \cite{AlS94} [appearing, 
e.g., in context of the stretched isofactors of $u_q(3)$], which are 
expressed as compositions of two generic $_3F_2[\cdot \cdot \cdot ;q,x]$ 
series, they are less symmetric and more complicate. Although these 
$_3F_2[\cdot \cdot \cdot ;q,x]$ series in all 23 new expressions may be 
rearranged into other forms separately [cf.\ Refs.\ \onlinecite{JB77,%
AsST96,GR90}] in such ways that the double sums turn into compositions of 
two generic $_3F_2[\cdot \cdot \cdot ;q,x]$ series, we use more universal 
approach in Sec.~IV.

The doubly stretched $q$-$9j$ coefficients with the adjacent consecutive 
stretched triangles may be expressed without sum,
\begin{eqnarray}
\left\{ \begin{array}{ccc}
c\!+\!h\, & b & b\!+\!c\!+\!h \\ 
c & d & f \\ 
h & k & g
\end{array}\right\} _{\!\!q} &=&\frac{(-1)^{c+d-f}q^{2bc+Z_{efhk}}
\nabla [efg]}{\nabla [hgk]\nabla [cdf]\nabla [bdk]}\left( \frac{[2b]!
[2c]![2h]!}{[2a+1][2e+1]!}\right) ^{1/2},  \label{st2a} 
\end{eqnarray}
where $a=c+h$ and $e=a+b$ [cf.\ Eq.\ (32.21) of Ref.\ \onlinecite{JB77}], 
taking into account that in related 11 cases two summation parameters are 
fixed and the last summation may be performed using either the 
Chu--Vandermonde formulas \cite{JB77,AsST96,GR90} [in 9 cases, e.g., in 
(\ref{trsa}) for $k=g+h=b-d$, or in (\ref{trsb}) for $k=g+h=d-b$, see 
Appendix B] or Karlsson's \cite{K71,AlJJ72} summation formula 
(\ref{smCq}) [in (\ref{trsd}) for $g=e+f=k-h$, or in (\ref{trse}) for 
$g=e+f=k-h$]. 

The different versions of the doubly stretched $q$-$9j$ coefficients 
with single sums in expressions [mainly as generalizations of Eqs.\ 
(32.15), (32.17), (32.17a), (32.18), and (32.20) of Ref.\ 
\onlinecite{JB77}] may be obtained straightforwardly from (\ref{trsa})--%
(\ref{trse}) with fixed couples of summation parameters. 

Additional restrictions for $z_i\pm z_j$ in generic expressions 
(\ref{trsa})--(\ref{trse}) may be represented as some couples of triangle 
linear combinations. No formula does represent any usual symmetry of 
$9j$-symbol, but expressions (\ref{trsa})--(\ref{trse}) are mutually 
related by some ``mirror reflection'' ($j\to -j-1$) 
symmetries.\cite{JB77,JB65}

\subsection{Expressions with the partial triangle restrictions of 
summation intervals}

Summation formulas (\ref{smDqb}) and (\ref{smDqc}) also may be used, 
when the first (primed) $q$-$6j$ coefficients in (\ref{vrx6a}) and 
(\ref{vrx6b}) are expressed by means of (\ref{f6jcb}) and remaining 
(primed or nonprimed) ones by means of (\ref{f6jc}) or (\ref{f6jb}), 
respectively. It is impossible to get the definite summation interval 
for $j$ when expressing all three $q$-$6j$ coefficients by means of 
primed $q$-$6j$ coefficients (\ref{f6jc}) or (\ref{f6jcb}) with the 
numerator -- denominator distribution of the type (\ref{dndcd}).
Hence we derive in addition two more triple sum expressions for $q$-$9j$ 
coefficients,
\begin{mathletters}
\begin{eqnarray}
\left\{ \begin{array}{lll}
a & b & e \\ 
c & d & f \\ 
h & k & g
\end{array} \right\} _{\!\!q} &=&(-1)^{a+c-h}
\frac{\nabla [abe]\nabla [feg]\nabla [kbd]}{\nabla [ach]\nabla [fcd]
\nabla [kgh]}  \nonumber \\
&&\times q^{Z_{deh}-(a+b+f-g)(f+g-e+1)-(g-a)(g-a+1)}  
\nonumber \\
&&\times \sum_{z_1,z_2,z_3}\frac{(-1)^{z_1+z_2+z_3}[g-h+k+z_1]!
[g+h+k+z_1+1]!}{[z_1]![g+k-a-c+z_1]![c+g+k-a+z_1+1]!}  \nonumber \\
&&\times \frac{[2b-z_2]![b+f+k-c-z_2]![b+f+k+c-z_2+1]!}{[2g+z_1+1]![z_2]!
[b+k-d-z_2]![b+d+k-z_2+1]!}  \nonumber \\
&&\times \frac{[b+e-a+z_3]![e-f+g+z_3]!q^{-z_1(b+e-a-z_2+z_3)}}{[z_3]!
[a+b-e-z_3]![f+g-e-z_3]![2e+z_3+1]!}  \nonumber \\
&&\times \frac{q^{z_2(f+g-e-z_3+1)+z_3(a+b+f-g)}[f+g-e+z_1-z_3]!}{%
[a+b+f-g-z_1-z_2]![b+e-a-z_2+z_3]!}  \label{trsxa} \\
&=&(-1)^{a+c-h-e-f+g}\frac{\nabla [abe]\nabla [feg]\nabla [kbd]}{\nabla
[ach]\nabla [fcd]\nabla [kgh]}  \nonumber \\
&&\times q^{Z_{deh}-(a+b+f-g)(e+f+g+2)-(g-a)(g-a+1)}  
\nonumber \\
&&\times \sum_{z_1,z_2,z_3}\frac{(-1)^{z_1+z_3}[g-h+k+z_1]![g+h+k+z_1+1]!%
}{[z_1]![g+k-a-c+z_1]![c+g+k-a+z_1+1]!}  \nonumber \\
&&\times \frac{[2b-z_2]![b+f+k-c-z_2]![b+f+k+c-z_2+1]!}{[2g+z_1+1]![z_2]!
[b+k-d-z_2]![b+d+k-z_2+1]!}  \nonumber \\
&&\times \frac{[2e-z_3]!q^{z_1(a-b+e+z_2-z_3+1)+z_2(e+f+g-z_3+2)+
z_3(a+b+f-g)}}{[z_3]![a-b+e-z_3]![e+f-g-z_3]![[a+b+e-z_3+1]!}  
\nonumber \\
&&\times \frac{[e+f+g+z_1-z_3+1]![a-b+e+z_2-z_3]!}{[e+f+g-z_3+1]!
[a+b+f-g-z_1-z_2]!}.  \label{trsxb}
\end{eqnarray}
\end{mathletters}
Here the summation intervals for $z_2$ and $z_3$ are restricted by one or
two triangle conditions, but $z_1+z_2$ restricted only both together by
some couples of triangular linear combinations. In these cases we may 
write 6 more expressions for the stretched $q$-$9j$ coefficients as 
double sums, but only five of them correspond to compositions of generic 
$_4F_3[\cdot \cdot \cdot ;q,x]$ and $_3F_2[\cdot \cdot \cdot ;q,x]$ 
series, and only four times all the summation intervals in these  
expressions are restricted by some triangle conditions. In the remaining 
cases some couples of triangular linear combinations appear as the 
summation intervals.

\section{REARRANGEMENT OF THE TRIPLE SUM EXPRESSIONS FOR  
$\lowercase{q}$-$9\lowercase{j}$ COEFFICIENTS AND STRETCHED 
$\lowercase{q}$-$9\lowercase{j}$ COEFFICIENTS}

\subsection{Search for other rearrangement of the triple sum expressions}

We may identify such three blocks (quintuplets) of factorials under the
summation sign in numerators and denominators of each expression 
(\ref{trsa})--(\ref{trse}), which may be expanded using the 
Chu--Vandermonde summation formulas,\cite{JB77,AsST96,GR90} given in 
Appendix B. For example, expressions (\ref{trsa}), (\ref{trsb}), and 
(\ref{trsd}) may be expanded as follows: 
\begin{mathletters}
\begin{eqnarray}
\left\{ \begin{array}{lll}
a & b & e \\ 
c & d & f \\ 
h & k & g
\end{array} \right\} _{\!\!q} &=&(-1)^{c+h-a}
\frac{\nabla [abe]\nabla [feg]\nabla [kbd]}{\nabla [ach]\nabla [fcd]
\nabla [kgh]}q^{(f+h-e-k)(a+d-e+k+1)-(a-e+f)(a-e+f+1)}  \nonumber \\
&&\times q^{-(g+h-k)(f+g-e)-(c-a+h)(c+d-f)+(b+e-a+1)(d-b+k)+Z_{deh}} 
\nonumber \\
&&\times \sum_{z_1,z_2,z_3,s_1,s_2,s_3}\frac{q^{-z_1(a-c-g+k+1)-z_2(b-c+
f-k+1)+z_3(a-b+f+g)}[2h-z_1]![2d-z_2]!}{[z_1]![z_2]![z_3]![2e+z_3+1]!\
[s_1]![s_2]![s_3]!}  \nonumber \\
&&\times \frac{(-1)^{z_1+z_2+z_3+s_1+s_2+s_3}q^{s_1(a+d-e+k-z_2-z_3)}[b+e
-a+z_3+s_1]!}{[d-b+k-z_2-s_1]![2b+s_1+1]![g+h-k-z_1-s_2]!}  \nonumber \\
&&\times \frac{q^{-s_2(e-f-h+k+z_1+z_3+1)-s_3(a-d+f-h+z_1+z_2+1)}[2g-s_2]!
[2c-s_3]!}{[f+g-e-z_3-s_2]![c-a+h-z_1-s_3]![c+d-f-z_2-s_3]!}  
\label{trxa} \\
&=&(-1)^{e-f-h+k}\frac{\nabla [abe]\nabla [feg]\nabla [kbd]}{\nabla [ach]
\nabla [fcd]\nabla [kgh]}q^{(c+h-a)(b-c+f+k+1)}  \nonumber \\
&&\times q^{(b+e-a)(e-f-h+k)-(f+g-e)(g+h-k)-(a-e+f+1)(a-e+f)+Z_{deh}}  
\nonumber \\
&&\times \sum_{z_1,z_2,z_3}\frac{(-1)^{z_1+z_2}q^{-z_1(a-c-g+k+1)-z_2
(b+c-f+k)+z_3(a+b+f+g)}[2h-z_1]!}{[z_1]![b-d+k-z_2]![b+d+k-z_2+1]!
[z_3]![2e+z_3+1]!}  \nonumber \\ 
&&\times \sum_{s_1,s_2,s_3}\frac{(-1)^{s_1+s_2}q^{-s_1(b+e-a-z_2+z_3+1)
-s_2(e-f-h+k+z_1+z_3+1)}[2b-s_1]![2g-s_2]!}{[s_1]![z_2-s_1]![a+b-e-z_3
-s_1]![s_2]![g+h-k-z_1-s_2]!}  \nonumber \\
&&\times \frac{q^{-s_3(a+b+f-h+k+z_1-z_2+2)}[2c-s_3]![b-c+f+k-z_2+s_3]!}{%
[f+g-e-z_3-s_2]![s_3]![c+h-a-z_1-s_3]!}  \label{trxb} \\
&=&(-1)^{b-c+f-g+h}q^{(k-c-1)(a+b+f-g)+(a-g)(b+f+1)+(e+f-g)(a-b-f+g+1)}  
\nonumber \\
&&\times q^{(a-c+h)(c-a+g+k+2)+Z_{deh}} \frac{\nabla [ach]\nabla [feg]
\nabla [kbd]}{\nabla [abe]\nabla [fcd]\nabla [kgh]}  \nonumber \\
&&\times \sum_{z_1,z_2,z_3}\frac{q^{-z_1(c-a+g+k+2)+z_2(c-b+f-k+1)
-z_3(a-b-f+g+1)}[2h-z_1]!}{[z_1]![b-d+k-z_2]![b+d+k-z_2+1]![z_3]!}  
\nonumber \\
&&\times \sum_{s_1,s_2,s_3}\frac{(-1)^{z_2+z_3+s_1+s_2+s_3}
[2e-z_3]![2b-s_1]!}{[s_1]![z_2-s_1]![b-a+e-z_3-s_1]![s_2]![e+f-g-z_3
-s_2]!}  \nonumber \\
&&\times q^{-s_1(a+b-e-z_2+z_3+1)-s_2(k-e-f-h+z_1+z_3)
-s_3(b+f-h+k-a+z_1-z_2+1)}  \nonumber \\
&&\times \frac{[g-h+k+z_1+s_2]![b+c+f+k-z_2+s_3+1]!}{[2g+s_2+1]![s_3]!
[a-c+h-z_1-s_3]![2c+s_3+1]!}.  \label{trxd} 
\end{eqnarray}
\end{mathletters}
The summations over $s_1,s_2,s_3$ give original expressions (\ref{trsa}), 
(\ref{trsb}), and (\ref{trsd}), when the summations of (\ref{trxa}), 
over $z_1,z_2,z_3$ give another expression for $q$-$9j$ coefficient, 
equivalent to (\ref{trsa}) after transpositions of two last rows and two 
last columns,
\begin{mathletters}
\begin{equation}
\left\{ \begin{array}{lll}
a & b & e \\ 
c & d & f \\ 
h & k & g
\end{array} \right\} _{\!\!q}=\left\{ \begin{array}{lll}
a & e & b \\ 
h & g & k \\ 
c & f & d
\end{array} \right\} _{\!\!q}.  \label{sim9a}
\end{equation}
Otherwise, the summations of (\ref{trxb}) over $z_1,z_2,z_3$ give 
expression, equivalent to (\ref{trse}), after changing the summation 
parameters and taking into account the same relation (\ref{sim9a}), as 
well as the summations of (\ref{trxd}) over $z_1,z_2,z_3$ give 
expression, equivalent to (\ref{trsc}), again after change of summation 
parameters and applying the relation 
\begin{equation}
\left\{ \begin{array}{lll}
a & b & e \\ 
c & d & f \\ 
h & k & g
\end{array} \right\} _{\!\!q}=\left\{ \begin{array}{lll}
f & d & c \\ 
e & b & a \\ 
g & k & h
\end{array} \right\} _{\!\!q}.  \label{sim9b}
\end{equation}
\end{mathletters}
Hence only three from these expressions for the $q$-$9j$ coefficients are
independent with respect to elementary rearrangements.

The quintuplet expansion by means of the Chu--Vandermonde summation 
formulas of expressions (\ref{trsxa}) or (\ref{trsxb}) leads to vanishing 
of the summation limit for $z_1$ and, therefore, it is not helpful for 
the rearrangement of $q$-$9j$ coefficients.

\subsection{Different expressions for the stretched $q$-$9j$ coefficients}

In the stretched cases, e.g., for $k=g+h$ in (\ref{trxa}) and 
(\ref{trxb}), or for $c=a+h$ in (\ref{trxd}) some couples of 
parameters $z_i$ and $s_j$ are fixed and summation over $z_l$ and $s_l$
(where $i,j,l$ is some permutation of $1,2,3$) is possible, using the 
Chu--Vandermonde formulas (see Appendix B). Hence we may derive 14 
versions of expressions (from which at least 13 are independent) for the 
stretched $q$-$9j$ coefficients as double sums over parameters $z_j$ and 
$s_i$ (where further the subscripts of the summation parameters will be 
omitted) as compositions of the both generic $_3F_2[\cdot \cdot 
\cdot ;q,x]$ series. For example, from (\ref{trxa}) and (\ref{trxb}) 
with $a=c+h$ and $z_1=s_3=0$, from (\ref{trxa}) with $k=h+g$ and 
$z_1=s_2=0$ or with $e=f+g$ and $z_3=s_2=0$, and from (\ref{trxd}) 
with $c=a+h$ and $z_1=s_3=0$ [using some symmetries (\ref{symev}) of the 
$q$-$9j$ coefficients and, in the last case, some change of summation 
parameter] we obtain, respectively, the following expressions:
\begin{mathletters}
\begin{eqnarray}
\left\{ \begin{array}{ccc}
a & b & e \\ 
c & d & f \\
h & k & e\!+\!f
\end{array} \right\} _{\!\!q} &=&\left( \frac{[2e]![2f]!}{[2g+1]!}%
\right) ^{1/2}\frac{\nabla [ghk]\nabla [bdk]\nabla [cah]}{\nabla [eab]
\nabla [fcd]}q^{(h-b-c-e)(b+c+f-k)+2bc+Z_{efhk}}  \nonumber \\
&&\times \sum_{s,z}\frac{(-1)^{s+z}[a-c+h+s]![k-b+d+z]!}{[s]![a+c-h-s]!
[b-c-e+h+s]![2h+s+1]!}  \nonumber \\
&&\times \frac{q^{s(b+c+f-k-z)+z(b+c+e-h)}[h-g+k+s+z]!}{[z]!
[b+d-k-z]![k-b+c-f+z]![2k+z+1]!}  \label{stra} \\
&=&(-1)^{b-a+f+h-k}\left( \frac{[2e]![2f]!}{[2g+1]!}\right) ^{1/2}
\frac{\nabla [ghk]\nabla [bdk]\nabla [cah]}{\nabla [eab]
\nabla [fcd]}  \nonumber \\
&&\times q^{(b+c-f)(k-g-h)+k(g+h+1)-2fh+Z_{bcf}+Z_{ghk}}  
\nonumber \\
&&\times \sum_{s,z}\frac{(-1)^{s}[2h-s]![e+h-b+c-s]![k-b+d+z]!}{[s]!
[c+h-a-s]![a+c+h-s+1]![g+h-k-s-z]!}  \nonumber \\
&&\times \frac{q^{s(b+c+f-k-z)+z(b+c+e+h+1)}}{[z]!
[b+d-k-z]![c-b-f+k+z]![2k+z+1]!}  \label{strb} \\
&=&(-1)^{d+f-c}\left( \frac{[2e]![2f]!}{[2g+1]!}\right) ^{1/2}
\frac{\nabla [cah]\nabla [ghk]}{\nabla [fcd]\nabla [eab]
\nabla [bdk]}  \nonumber \\
&&\times q^{(b-c+e+h+1)(g-h+k)-(e+k)(e+k+1)+Z_{bcg}}  \nonumber \\
&&\times \sum_{s,z}\frac{(-1)^{s+z}[a-c+h+s]![b+c+e-h-s]!}{[s]![a+c-h-s]!
[2h+s+1]![g-h+k-s-z]!}  \nonumber \\
&&\times \frac{q^{s(c-b+f+k-z)-z(b-c-e+h+1)}[b+d-k+z]![2k-z]!}{[z]!
[d+k-b-z]![b+c-f-k+z]!}  \label{strc} \\
&=&q^{2af-(a+b-e)(a+d+f-h)+Z_{bdgh}}\left( \frac{[2e]![2f]!}{[2g+1]!}
\right) ^{1/2}\frac{\nabla [ghk]\nabla [dbk]}{\nabla [hac]\nabla [fcd]
\nabla [eab]}  \nonumber \\
&&\times \sum_{z,s}\frac{(-1)^{a+c-h+s+z}[2b-z]![b+d-g+h-z]!
}{[z]![b+d-k-z]![b+d+k-z+1]![s]!}  \nonumber \\
&&\times \frac{q^{z(a+d+f-h-s)-s(g+h-b-d+1)}[2a-s]![h-a+c+s]!
}{[a+b-e-z-s]![a+c-h-s]![d-a-f+h+s]!}  \label{strd} \\
&=&q^{k(2f-2b+k-1)+(g-h+k)(a+b-f-h-k)+Z_{cdeh}}  \nonumber \\
&&\times \left( \frac{[2e]![2f]!}{[2g+1]!}\right) ^{1/2}
\frac{\nabla [ghk]\nabla [kbd]\nabla [ach]}{\nabla [fcd]\nabla [eab]}  
\nonumber \\
&&\times \sum_{s,z}\frac{(-1)^{c-d+e-g+z}q^{-s(a+b+e+z+1)-z(h-g+k)}}{%
[s]![g-h+k-s]![a+c-h-s]![2h+s+1]!}  \nonumber \\
&&\times \frac{[c+h-a+s]![b+d-k+z]![b+e-a+z]!}{[z]![d+k-b-z]!
[b-a-f+h-k+s+z]![2b+z+1]!}.  \label{stre}
\end{eqnarray}
\end{mathletters}
Expression (\ref{stra}) is invariant with respect to simultaneous 
permutations of the two first columns and rows of the stretched $q$-$9j$ 
coefficients in accordance with (\ref{symod}) and is related to the 
particular case of Eqs.\ (26)--(27) of Ref.\ \onlinecite{AlJ71}, that 
appeared in context of the stretched isoscalar factors of the Sp(4) or 
SO(5) group restricted to SU(2)$\times $SU(2), but (\ref{strb}) and 
(\ref{stre}) are more convenient, since separate sums are not 
alternating. 

The linear combinations of parameters $a+b-e$, $d-b+k$, and $c+d-f$ 
restrict the both summation parameters in expressions (\ref{strb}), 
(\ref{strc}), and (\ref{strd})--(\ref{stre}), respectively, for 
the $q$-$9j$ coefficients with the couples of adjacent consecutive 
stretched triangles [cf.\ Eq.\ (32.21) of Ref.\ \onlinecite{JB77}]. 
Otherwise, the summation of expression (\ref{stra}) for $g+h-k=0$, or 
$g-h+k=0$ is nontrivial, as well as Eq.\ (\ref{strb}) for $g-h+k=0$, 
Eq.\ (\ref{strc}) for $a+b-e=0$, or $g+h-k=0$, and Eqs.\ 
(\ref{strd})--(\ref{stre}) for $g+h-k=0$. 

The both separate sums only in Eqs.\ (\ref{stra}) and (\ref{strd})
correspond to the CG coefficients of $u_q(2)$.\cite{N89q,GKK90,STK91a} 
Hence the sum over $z$ in (\ref{strd}) [as well as the sum over $s$ in 
(\ref{strb})] may be included [using Eq.\ (5.17) of Ref.\ 
\onlinecite{STK91a}] into the Clebsch--Gordan coefficients of $u_q(2)$, 
reexpressed by means of Eq.\ (41a) of Ref.\ \onlinecite{GKK90} (with 
changed summation parameters) and the following expressions for the 
stretched $q$-$9j$ coefficients may be derived:
\begin{mathletters}
\begin{eqnarray}
\left\{ \begin{array}{ccc}
a & b & e \\ 
c & d & f \\
h & k & e\!+\!f
\end{array} \right\} _{\!\!q} &=&(-1)^{a+b+c+d-h-k}\frac{\left( [2e]!
[2f]![g+h-k]![g+h+k+1]!\right) ^{1/2}}{\left( [2g+1]![2k+1]\right) ^{1/2}
\nabla [hac] \nabla [fcd]\nabla [eab]}  \nonumber \\
&&\times q^{2af-(a+b-e)(a+d+f-h)+b(g-h)+(b+d-k)(b+d+k+1)/2+Z_{bdgh}} 
\nonumber \\
&&\times \sum_m(-1)^{a-e+m}q^{-m(g+h+1)}\frac{[a+e-m]![c-e+h+m]!}{[a-e+m]!
[c+e-h-m]!}  \nonumber \\
&&\times \left( \frac{[d+g-h-m]![b+m]!}{[d-g+h+m]![b-m]!}\right) ^{1/2}
\left[ \begin{array}{ccc}
d & b & k \\ 
g-h-m & m & g-h
\end{array} \right] _{\!q}  \label{stdCG} \\
&=&\frac{\left( [2e]![2f]![h+k-g]![g-h+k]!
[g+h-k]![g+h+k+1]!\right) ^{1/2}}{\left( [2g+1]!\right) ^{1/2}\nabla 
[hac] \nabla [fcd]\nabla [eab]\nabla [kbd]}  \nonumber \\
&&\times q^{(b+d-k)(a+b+f-h+k+1)-(a+b-e)(a+d+f-h)+2af+Z_{bdgh}}  
\nonumber \\
&&\times \sum_{s,z}\frac{(-1)^{a+b+c+d-h-k+s+z}q^{-s(g+h-k+1)-z(g-h+k+1)}%
}{[s]![a+c-h-s]![z]![b+d-k-z]!}  \nonumber \\
&&\times \frac{[2a-s]![c+h-a+s]![2d-z]![b-d+k+z]!}{[d-a-f+h+s-z]!
[a-d-e+k-s+z]!}  \label{stdK} \\
&=&\left( \frac{[2e]![2f]![a-b+e]![a+b-e]!}{[2g+1]![e-a+b]![a+b+e+1]!}%
\right) ^{1/2}\frac{\nabla [ghk]\nabla [bdk]}{\nabla [fcd]\nabla [ach]}  
\nonumber \\
&&\times q^{(c+d-f)(a+b-e+1)-(a+b+e+1)(b+d-k)+2ed+Z_{acgk}}  \nonumber \\
&&\times \sum_{s,z}\frac{[2c-s]![a-c+h+s]![k-b+d+z]!}{[s]![c+h-a-s]!
[z]![b+d-k-z]![2k+z+1]!}  \nonumber \\
&&\times \frac{(-1)^{c+k-b-f+s}q^{-s(a+b-e+1)+z(a+b+e+1)}}{[c-b-f+k+z-s]!
[a-c+g-k-z+s]!}  \label{stbK}
\end{eqnarray}
\end{mathletters}
Expression (\ref{stdK}) satisfies symmetry relation (\ref{symod}) for 
permutations of the two first columns or rows of the stretched $q$-$9j$ 
coefficient and is a $q$-generalization of standard formula (32.13) of 
Ref.\ \onlinecite{JB77} for the stretched $9j$ coefficients. The linear 
combinations of parameters $a+b-e$ or $c+d-f$ restrict the both 
summation parameters in expression (\ref{stdK}), when $g+h-k$ or 
$c+d-f$ restrict the summation limits in expression (\ref{stbK}). Hence 
expressions (\ref{stdK}) and (\ref{stbK}) for these cases of adjacent 
consecutive stretched triangles also turn into single terms [cf.\ Eq.\ 
(\ref{st2a})]. Note that the Chu--Vandermonde or Karlsson summation 
formulas are needed for the last summation of 11 doubly stretched cases 
of triple sum expressions (\ref{trsa})--(\ref{trse}), e.g., in 
(\ref{trsa}) for $k=g+h=b-d$, or in (\ref{trsd}) for $g=e+f=k-h$.

For the diverging adjacent stretched triangles, e.g., with $e=b-a=g-k$, 
summation parameters in (\ref{stbK}) are dependent and the doubly 
stretched $q$-$9j$ coefficients may be expressed as single sums with the 
alternating terms [cf.\ Eq.\ (32.18) of Ref.\ \onlinecite{JB77}], when 
using Eq.\ (\ref{trsa}) for $c=f-d=a-h$ an equivalent formula may be 
written directly. Again, for the merging adjacent stretched triangles, 
e.g., with $g=h+k=e+f$, the doubly stretched $q$-$9j$ coefficients may be 
expressed as single sums with the fixed sign of all terms [cf.\ Eq.\ 
(32.20) of Ref.\ \onlinecite{JB77}] by means of formula (\ref{stdK}) 
[as well as for $b=d+k=a+e$ by means of Eq.\ (\ref{trse}), in contrast 
with the remaining formulas of this paper]. These single sums in the both 
cases are related to the generic $_4F_3[\cdot \cdot \cdot ;q,x]$ series.

The $q$-$9j$ coefficients with the both stretched triangles appearing in 
the different layers, or rows of the $q$-$9j$-symbol are related to 
generic $_3F_2[\cdot \cdot \cdot ;q,x]$ series. In the case of two 
parallel stretched triangles [e.g., for $h=a+c$ in (\ref{stra}) or 
(\ref{stdK})] they are proportional [cf.\ Eq.\ (32.15a) of Ref.\ 
\onlinecite{JB77}] to the CG coefficients of $u_q(2)$ (see Refs.\
\onlinecite{N89q,GKK90,STK91a}), with appearing two different types of 
expressions. Otherwise, special $q$-$9j$ coefficients with two 
antiparallel stretched triangles may be expressed in four different forms 
[see (\ref{strb}) for $a=c+h$, or (\ref{strc}) for $b=d+k$ as 
generalizations of Eqs.\ (32.17a) and (32.17) of Ref.\ \onlinecite{JB77}, 
as well as Eq.\ (\ref{trsc}) for $h=a+c$ and $b=d+k$, or Eq.\ 
(\ref{trse}) for $b=d+k$ and $g=e+f$] and correspond to the CG 
coefficients of $u_q(1,1)$, with the expressions including either the 
alternating terms [with diverse distribution of summation parameter signs 
in two numerator $q$-factorial arguments, in analogy with CG coefficients 
of $u_q(2)$], or the fixed sign terms (with one or three numerator 
$q$-factorial arguments). Note, that expressions for the triply stretched 
$q$-$9j$ coefficients with the three mutually antiparallel stretched 
triangles [e.g., (\ref{stra}) or (\ref{stbK}) for $a=c+h$ and $d=b+k$] 
are not summable.

Expression related to (\ref{stdCG}) may be derived [in contrast with 
the intermediate version of (\ref{stbK})] also from expansion 
\cite{N91,STK90} [cf.\ Ref.\ \onlinecite{Sh67} in the SU(2) case] of the 
$q$-$9j$ coefficients in terms of the Clebsch--Gordan coefficients of 
$u_q(2)$ [cf.\ Eq.\ (3.12) of Ref.\ \onlinecite{N91}] which in the 
stretched case with $h=a+c$ (and with extreme CG coefficient for coupling 
$a\times c\rightarrow h$ in the r.h.s.\ equal to 1) may be written as 
follows:
\begin{eqnarray}
\left\{ \begin{array}{ccc}
a & b & e \\ 
c & d & f \\ 
a\!+\!c & k & g
\end{array} \right\} _{\!\!q^{-1}}&=& \frac{q^{-Z_{efhk}}}{([2e+1][2f+1]
[2h+1][2k+1])^{1/2}}\left\{ \left[ \begin{array}{ccc}
h & k & g \\ h & g\!-\!h & g
\end{array} \right] _{\!q} \right\} ^{-1}  \nonumber \\
&&\times \sum_{m}\left[ \begin{array}{ccc}
a & b & e \\ a & m\!-\!a & m
\end{array} \right] _{\!q}\left[ \begin{array}{ccc}
e & f & g \\ m & g\!-\!m & g
\end{array} \right] _{\!q}\left[ \begin{array}{ccc}
c & d & f \\ c & g\!-m\!-\!c & g\!-\!m
\end{array} \right] _{\!q}\nonumber \\
&&\times \left[ \begin{array}{ccc}
b & d & k \\ m\!-\!a & g\!-m\!-\!c & g\!-\!h
\end{array} \right] _{\!q}
\left( R^{cb}\right) _{c\;m\!-\!a}^{c\;m\!-\!a},  \label{vrCG}
\end{eqnarray}
where 
\[
\left( R^{cb}\right) _{c\;m\!-\!a}^{c\;m\!-\!a}=q^{2c(m-a)}
\]
is a diagonal extreme element of triangular braiding $R$-matrix and all 
the CG coefficients with exception of the last one may be expressed 
without sum. It should be noted that only in special stretched case 
(\ref{vrCG}) the summation over non-diagonal elements of $R$-matrix may 
be escaped.

Both expressions (\ref{stdK}) and (\ref{stbK}) correspond to 
$q$-generalizations of the Kamp\'e de F\'eriet\cite{K-F21} function 
$F_{1:1}^{1:2}$, which is defined as follows:
\begin{eqnarray}
&&^{\pm \!}F_{C:D}^{A:B}\left[ \begin{array}{c}
(a) \\ (c)
\end{array}:\begin{array}{c}
(b) \\ (d) 
\end{array};\begin{array}{c}
(b^{\prime }) \\ (d^{\prime }) 
\end{array};x,y;q\right]  \nonumber \\
&& \qquad =\sum_{s,t}^{\infty }\frac{\prod_{j=1}^A(a_j|q)_{s+t}}{%
\prod_{j=1}^C(c_j|q)_{s+t}}\frac{\prod_{j=1}^B(b_j|q)_s(b_j^{\prime %
}|q)_t}{\prod_{j=1}^D(d_j|q)_s(d_j^{\prime }|q)_t}\frac{x^{\pm s }
y^{\pm (1-2\delta )t}}{[s]![t]!}q^{\pm (A-C)st},  \label{sK}
\end{eqnarray}
with special parameters 
\begin{eqnarray*}
x&=&q^{p+1},\quad p=\tsum_{j=1}^Aa_j+\tsum_{j=1}^Bb_j-\tsum_{j=1}^Cc_j-
\tsum_{j=1}^Dd_j, \\
y&=&q^{p^{\prime }+1},\quad p^{\prime }=\tsum_{j=1}^Aa_j+\tsum_{j=1}^B%
b^{\prime }_j-\tsum_{j=1}^Cc_j-\tsum_{j=1}^Dd^{\prime }_j, \\
\delta &=&\delta _{AC}
\end{eqnarray*} 
for $A+B=C+D+1$ and $|A-C|\leq 1$. Of course, series (\ref{sK}) turn 
into usual Kamp\'e de F\'eriet function $F_{C:D}^{A:B}[\cdot \cdot \cdot %
;1,1]$ for $q=1$. Unfortunately, the standard definition 
\cite{V-JPS-R94,PV-J96}  
\begin{equation}
\Phi_{C:D}^{A:B}\left[ 
\begin{array}{c}
(\alpha ) \\ (\gamma )
\end{array}:\begin{array}{c}
(\beta ) \\ (\delta )
\end{array};\begin{array}{c}
(\beta ^{\prime }) \\ (\delta ^{\prime }) 
\end{array};x,y;q\right]  \label{dfqK}
\end{equation} 
(cf.\ Refs. \ \onlinecite{V-JPS-R94,PV-J96}) in terms of the asymmetric 
$q$-factorials\cite{GR90} 
\[
(a;q)_n=(1-a)(1-aq)\cdot \cdot \cdot (1-aq^{n-1}),\quad n=1,2,...,  
\]
derived after substitution
\begin{eqnarray}
&&q\rightarrow q^{1/2},\quad [n]!\rightarrow (q;q)_n(q^{-1/2}-q^{1/2})%
^{-n}q^{-n(n+1)/4},  \nonumber \\
&&(a|q)_n\rightarrow (q^a;q)_n(q^{-1/2}-q^{1/2})%
^{-n}q^{-n(2a+n-1)/4},  \label{sasqf}
\end{eqnarray}
may be not convenient in the double sums that appear in (\ref{stdK}) and 
(\ref{stbK}), since arguments $x$ and $y$ do not turn into $q$ both 
together, but turn either into $q$ and $q^{-p^{\prime }}$ or into 
$q^{-p}$ and $q$, respectively. 

Furthermore, the $q=1$ versions of Eqs.\ (\ref{stra})--(\ref{stre}) 
correspond to the Kamp\'e de F\'eriet \cite{K-F21} functions 
$F_{0:2}^{1:2}$ or [after reversing the order of summations, associated 
with spoiling some natural restrictions for both summation parameters in 
(\ref{strb})--(\ref{strd})] to $F_{1:1}^{0:3}$ (cf.\ Refs.\ 
\onlinecite{V-JPS-R94,PV-J96}). Otherwise, in the generic $q\neq 1$ case 
they can be expressed in terms of our (\ref{sK}) as 
$^{-}F_{0:2}^{1:2}[\cdot \cdot \cdot ;x,y,q]$ or $^{+}F_{1:1}^{0:3}[\cdot 
\cdot \cdot ;x,y;q]$, but only $^{+}F_{C:D}^{A:B}[\cdot \cdot \cdot 
;x,y,q]$ is equivalent to some $\Phi_{C:D}^{A:B}[\cdot \cdot \cdot 
;q,q;q]$ and for $^{-}F_{C:D}^{A:B}[\cdot \cdot \cdot ;x,y,q]$ with 
$A\neq C$ the factors $q^{mn}$, spoiling standard definition of 
$\Phi_{C:D}^{A:B}$ functions, cannot be eliminated in the new expansion 
[cf.\ Eq.\ (9) of Ref.\ \onlinecite{V-JPS-R94}], unless transition 
$q\rightarrow q^{-1}$ is performed preliminary.

Special rearrangement and summation formulas of the double $q$-factorial 
series and related Kamp\'e de F\'eriet \cite{K-F21} functions are given in
Appendix C.

\section{EXPRESSIONS FOR 12$\lowercase{j}$ COEFFICIENTS OF THE SECOND 
KIND OF SU(2) and \lowercase{$u_{q}$}(2)}

\subsection{Generic properties}

The $3nj$ coefficients of the second kind\cite{JLV60,LV57} ($n\geq 4$) 
whose graphs are planar [hence without braiding,\cite{N89bq} in contrast 
with the $3nj$ coefficients of the first kind ($n\geq 3$) whose graphs 
are possible only on the M\"{o}bius strip] usually are expanded 
\cite{JLV60,JB77,JB65} in terms of the factorized $n$ different $6j$ 
coefficients,
\begin{eqnarray}
&&\left[ \begin{array}{llllllll}
j_1\! &  & \!j_2\! &  & \!\cdot \cdot \cdot \! &  & \!j_n\! &  \\ 
& \!l_1\! &  & \!l_2\! &  & \!\cdot \cdot \cdot \! &  & \!l_n\! \\ 
k_1\! &  & \!k_2\! &  & \!\cdot \cdot \cdot \! &  & \!k_n\! & 
\end{array} \right]  =\sum_x(2x+1)(-1)^{R_n+nx}  \nonumber \\
&&\qquad \qquad \times \left\{ \begin{array}{lll}
j_1 & k_1 & x \\ 
k_2 & j_2 & l_1
\end{array} \right\} \left\{ \begin{array}{lll}
j_2 & k_2 & x \\ 
k_3 & j_3 & l_2
\end{array} \right\} \cdot \cdot \cdot \left\{ \begin{array}{lll}
j_{n-1} & k_{n-1} & x \\ 
k_n & j_n & l_{n-1}
\end{array} \right\} \left\{ \begin{array}{lll}
j_n & k_n & x \\ 
k_1 & j_1 & l_n
\end{array} \right\} ,  \label{d3nja} 
\end{eqnarray}
where
\[
R_n=\sum_{i=1}^n(j_i+k_i+l_i),
\]
and the triangular conditions are satisfied by the triplets of the 
nearest neighbors as $l_i,j_i,j_{i+1}$, or  $l_i,k_i,k_{i+1}$ 
($i=1,2,...,n-1$), or $l_n,j_n,j_1$, or  $l_i,k_n,k_1$, respectively. 

The $12j$ coefficients of the second kind, which may be extracted 
from the recoupling coefficients of the five irreps without braiding, 
\cite{N89bq} hence, with the cubic graph \cite{JLV60}
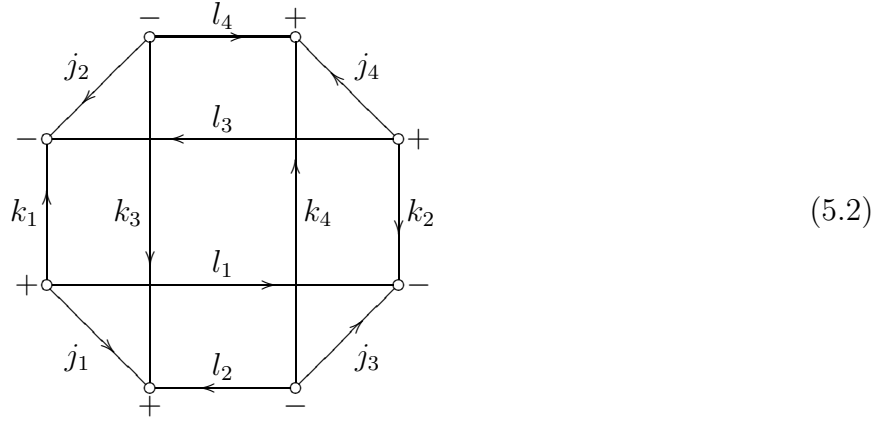
\begin{equation}
\xy
\xygraph{!{/r6pc/:}
[] 
!P8"A"{~>-\cir<2pt>{}}
!~:{@{-}|(0.65)@{>}}
"A1"*+!L{+}
"A2"*+!D{+}
"A3"*+!D{-}
"A4"*+!R{-}
"A5"*+!R{+}
"A6"*+!U{+}
"A7"*+!U{-}
"A8"*+!L{-}
"A1" (? :@{}_{\dpst j_4} "A2")
"A1" (? :@{}^{\dpst k_2} "A8")
"A3" (? :@{}^{\dpst l_4} "A2")
"A3" (? :@{}_{\dpst j_2} "A4")
"A5" (? :@{}^{\dpst k_1} "A4")
"A5" (? :@{}_{\dpst j_1} "A6")
"A7" (? :@{}_{\dpst j_3} "A8")
"A7" (? :@{}_{\dpst l_2} "A6")
"A4" ("A1" :_{\dpst l_3} ?)
"A8" ("A5" :^{\dpst l_1} ?)
"A6" ("A3" :_{\dpst k_3} ?)
"A2" ("A7" :_{\dpst k_4} ?)
}\endxy  \label{grp2r}
\end{equation}
were introduced by Elliott and Flowers \cite{EF55} and redefined by 
Vanagas and \v {C}iplys. \cite{VC58} These $12j$ coefficients and their 
$q$-generalizations satisfy 24 symmetries, \cite{JLV60,JB77,VC58} 
generated by the following substitutions:
\begin{mathletters} 
\begin{eqnarray}
\left[ \begin{array}{llll}
j_1 & j_2 & j_3 & j_4 \\ 
l_1 & l_2 & l_3 & l_4 \\ 
k_1 & k_2 & k_3 & k_4
\end{array} \right] _{\!q}&=&(-1)^{j_1-j_2-j_3+j_4+k_1-k_2-k_3+k_4}
\left[ \begin{array}{llllllll}
j_1\! &  & \!k_1\! &  & \!j_2\! &  & \!k_3\! &  \\ 
& \!l_1\! &  & \!l_3\! &  & \!l_4\! &  & \!l_2\! \\ 
j_3\! &  & \!k_2\! &  & \!j_4\! &  & \!k_4\! & 
\end{array} \right] _{\!q}  \label{dd0} \\
&=& \left[ \begin{array}{llll}
j_4 & j_3 & j_2 & j_1 \\ 
l_4 & l_3 & l_2 & l_1 \\ 
k_4 & k_3 & k_2 & k_1
\end{array} \right] _{\!q}=\left[ \begin{array}{llll}
l_1 & l_2 & l_3 & l_4 \\ 
k_1 & k_2 & k_3 & k_4 \\ 
j_1 & j_2 & j_3 & j_4
\end{array} \right] _{\!q}  \label{ddabc} \\
&=&\left[ \begin{array}{llll}
k_4 & k_2 & k_3 & k_1 \\ 
l_4 & l_2 & l_3 & l_1 \\ 
j_4 & j_2 & j_3 & j_1
\end{array} \right] _{\!q}=\left[ \begin{array}{llll}
j_1 & j_2 & j_3 & j_4 \\ 
l_2 & l_1 & l_4 & l_3 \\ 
k_3 & k_4 & k_1 & k_2
\end{array} \right] _{\!q}. \label{ddde}
\end{eqnarray}
\end{mathletters}
Eight triangular conditions may be visualized\cite{JB77} by means of the 
extended array
\[
\begin{array}{c}
\left[ \begin{array}{rrrr}
j_1 & j_2 & j_3 & j_4 \\ 
l_1 & l_2 & l_3 & l_4 \\ 
k_1 & k_2 & k_3 & k_4
\end{array} \right]  \\ 
\begin{array}{rrrr}
& j_2 \, & j_3 & \\ 
& l_2 \, & l_3 & 
\end{array} \end{array}
\]
and are satisfied by the triplets of parameters in the first and fourth
columns, as well as by the skew triplets descending from some parameter 
of the first or fourth column, e.g., by $l_1,k_2,j_3$, or by 
$j_4,l_3,k_2$.

Let's restrict ourselves to the following rearrangements of the $q$-$6j$ 
coefficients in expressions \cite{JLV60,JB77,JB65,N90} for the $q$-$12j$ 
coefficients of the second kind:
\begin{mathletters}
\begin{eqnarray}
\left[ \begin{array}{llll}
j_1 & j_2 & j_3 & j_4 \\ 
l_1 & l_2 & l_3 & l_4 \\ 
k_1 & k_2 & k_3 & k_4
\end{array} \right] _{\!q}&=&(-1)^{l_1-l_2-l_3+l_4}
\sum_{x}[2x+1]\left\{ \begin{array}{lll}
k_1 & j_1 & l_1 \\ 
j_3 & k_2 & x
\end{array} \right\} _{\!\!q}  \nonumber \\
&&\times \left\{ \begin{array}{lll}
k_3 & k_4 & x \\ 
j_3 & j_1 & l_2
\end{array} \right\} _{\!\!q}\left\{ \begin{array}{lll}
k_3 & j_2 & l_4 \\ 
j_4 & k_4 & x
\end{array} \right\} _{\!\!q}\left\{ \begin{array}{lll}
k_1 & k_2 & x \\ 
j_4 & j_2 & l_3
\end{array} \right\} _{\!\!q}  \label{dd6a} \\
&=&(-1)^{l_1-l_2-l_3+l_4}\sum_{x}[2x+1]\left\{ \begin{array}{lll}
k_1 & j_1 & l_1 \\ 
j_3 & k_2 & x
\end{array} \right\} _{\!\!q}  \nonumber \\
&&\times \left\{ \begin{array}{lll}
j_3 & x & j_1 \\ 
k_3 & l_2 & k_4
\end{array} \right\} _{\!\!q}\left\{ \begin{array}{lll}
k_3 & k_4 & x \\ 
j_4 & j_2 & l_4
\end{array} \right\} _{\!\!q}\left\{ \begin{array}{lll}
k_1 & x & k_2 \\ 
j_4 & l_3 & j_2
\end{array} \right\} _{\!\!q}  \label{dd6b} \\
&=&(-1)^{l_1-l_2-l_3+l_4}\sum_{x}[2x+1])\left\{ \begin{array}{lll}
k_1 & j_1 & l_1 \\ 
j_3 & k_2 & x
\end{array} \right\} _{\!\!q}  \nonumber \\
&&\times \left\{ \begin{array}{lll}
k_3 & k_4 & x \\ 
j_3 & j_1 & l_2
\end{array} \right\} _{\!\!q}\left\{ \begin{array}{lll}
j_4 & x & j_2 \\ 
k_3 & l_4 & k_4
\end{array} \right\} _{\!\!q}\left\{ \begin{array}{lll}
k_1 & x & k_2 \\ 
j_4 & l_3 & j_2
\end{array} \right\} _{\!\!q}  \label{dd6c} \\
&=&(-1)^{l_1-l_2-l_3+l_4}\sum_{x}[2x+1]\left\{ \begin{array}{lll}
j_3 & x & j_1 \\ 
k_1 & l_1 & k_2
\end{array} \right\} _{\!\!q}  \nonumber \\
&&\times \left\{ \begin{array}{lll}
k_3 & x & k_4 \\ 
j_3 & l_2 & j_1
\end{array} \right\} _{\!\!q}\left\{ \begin{array}{lll}
j_4 & x & j_2 \\ 
k_3 & l_4 & k_4
\end{array} \right\} _{\!\!q}\left\{ \begin{array}{lll}
l_3 & j_2 & k_1 \\ 
x & k_2 & j_4
\end{array} \right\} _{\!\!q}^{\!(\prime )},  \label{dd6d}
\end{eqnarray}
\end{mathletters}
with the asymmetric triangle coefficients depending on the summation
parameter $x$ distributed separately in the numerators or denominators of
each $q$-$6j$ coefficient in expansion (\ref{dd6a}), the mixed 
distribution of asymmetric triangle coefficients in the numerators and 
denominators of $q$-$6j$ coefficients in expansions (\ref{dd6b}) and 
(\ref{dd6c}), and resembling (\ref{dndcd}) distribution in (\ref{dd6d}). 
Using expression (\ref{f6jbb}) for $q$-$6j$ coefficients with summation 
parameter $x$ in the right lower position, Eq.\ (\ref{f6jb}) with 
inverted summation parameter for $q$-$6j$ coefficients with $x$ in the 
middle column, and Eq.\ (\ref{f6jb}) in the remaining cases, with 
exception of Eq.\ (\ref{f6jc}), used for the last $q$-$6j$ coefficient in 
(\ref{dd6d}), the asymmetric triangle coefficients depending on the 
summation parameter $x$ cancel. Then we may use the summation formula 
(\ref{p4s0}) for expansion (\ref{dd6a}) and formula \cite{Al97} 
(\ref{p4sa}) for expansions (\ref{dd6b})--(\ref{dd6d}).

\subsection{General expressions with fourfold sums}

This way we derived four different expressions for the $q$-$12j$
coefficients of the second kind,
\begin{mathletters}
\begin{eqnarray}
&&\left[ \begin{array}{llll}
j_1 & j_2 & j_3 & j_4 \\ 
l_1 & l_2 & l_3 & l_4 \\ 
k_1 & k_2 & k_3 & k_4
\end{array} \right] _{\!q}  \nonumber \\
&&\qquad =(-1)^{l_2-l_3+k_1-k_3-j_3+j_4}\frac{\nabla [k_3j_1l_2]
\nabla [j_3k_4l_2]\nabla [k_1j_2l_3]\nabla [j_4k_2l_3]}{\nabla [k_1j_1l_1]
\nabla [j_3k_2l_1]\nabla [k_3j_2l_4]\nabla [j_4k_4l_4]}  \nonumber \\
&&\qquad \quad \times \sum_{z_1,z_2,z_3,z_4}\frac{(-1)^{z_1+z_2+z_3+z_4}
[k_2+j_3-l_1+z_1]![k_1+j_1-l_1+z_1]!}{[z_1]![z_2]![z_3]![z_4]!
[l_1+k_2-j_3-z_1]![j_1+l_1-k_1-z_1]!}  \nonumber \\
&&\qquad \quad \times \frac{[2l_1-z_1]![j_1+l_2-k_3+z_2]![l_2-j_3+k_4+
z_2]!}{[j_1-l_2+k_3-z_2]![j_3+k_4-l_2-z_2]![2l_2+z_2+1]!}  \nonumber \\
&&\qquad \quad \times \frac{[j_4+k_4-l_4+z_3]![j_2+k_3-l_4+z_3]!
[2l_4-z_3]!}{[k_4+l_4-j_4-z_3]![j_2-k_3+l_4-z_3]![k_1+j_2-l_3-z_4]!
[2l_3+z_4+1]!}  \nonumber \\
&&\qquad \quad \times \frac{[j_2+l_3-k_1+z_4]![k_2+l_3-j_4+z_4]!}{%
[k_2-l_3+j_4-z_4]![k_1+k_3+j_3+j_4-l_1-l_4+z_1+z_3+1]!}  \nonumber \\
&&\qquad \quad \times \frac{[k_1+k_3+j_3+j_4-l_2-l_3-z_2-z_4]!}{%
[k_1+l_2-l_1-k_3+z_1+z_2]![l_3-l_1+j_3-j_4+z_1+z_4]!}  \nonumber \\
&&\qquad \quad \times \frac{[l_2+l_3-l_1-l_4+z_1+z_2+z_3+z_4]!}{[l_2
-l_4-j_3+j_4+z_2+z_3]![k_3-k_1+l_3-l_4+z_3+z_4]!}  \label{dasa} \\
&&\qquad =(-1)^{j_1-j_3-k_1+k_2-l_1-l_2-l_3+l_4}\frac{\nabla [j_3k_4l_2]
\nabla [k_3j_2l_4]\nabla [k_1j_2l_3]\nabla [j_4k_4l_4]}{\nabla [k_1j_1l_1]
\nabla [j_3k_2l_1]\nabla [k_3j_1l_2]\nabla [j_4k_2l_3]}  \nonumber \\
&&\qquad \quad \times \sum_{z_1,z_2,z_3,z_4}\frac{(-1)^{z_2+z_3+z_4}[k_2
+j_3-l_1+z_1]![k_1+j_1-l_1+z_1]!}{[z_1]![z_2]![z_3]![z_4]![l_1+k_2-j_3
-z_1]![j_1+l_1-k_1-z_1]!}  \nonumber \\
&&\qquad \quad \times \frac{[2l_1-z_1]![2l_2-z_2]![j_1-l_2+k_3+z_2]!}{%
[l_2+j_3-k_4-z_2]![j_1+l_2-k_3-z_2]![l_2+j_3+k_4-z_2+1]!}  \nonumber \\
&&\qquad \quad \times \frac{[j_2-k_3+l_4+z_3]![k_4+l_4-j_4+z_3]!}{[j_4+
k_4-l_4-z_3]![j_2+k_3-l_4-z_3]![2l_4+z_3+1]![k_1-j_2+l_3-z_4]!}  
\nonumber \\
&&\qquad \quad \times \frac{[2l_3-z_4]![k_2-l_3+j_4+z_4]!}{%
[k_2+l_3-j_4-z_4]![k_1+j_2+l_3-z_4+1]![k_1+k_3-l_1-l_2+z_1+z_2]!}  
\nonumber \\
&&\qquad \quad \times \frac{[j_3+j_4+l_2-l_4-z_2-z_3]![k_1+k_3+l_3-l_4
-z_3-z_4]!}{[k_1-k_3+j_3-j_4-l_1+l_4+z_1+z_3]![j_3+j_4-l_1-l_3+z_1+z_4]!}  
\nonumber \\
&&\qquad \quad \times \frac{[j_3-j_4+k_1-k_3+l_2+l_3-z_2-z_4]!}{
[l_1+l_2+l_3-l_4-z_1-z_2-z_3-z_4]!}  \label{dasb} \\
&&\qquad =(-1)^{k_1-k_2+l_1-l_2+l_3+l_4-j_2+j_4}\frac{\nabla [k_3j_1l_2]
\nabla [j_3k_4l_2]\nabla [j_4k_4l_4]\nabla [k_1j_2l_3]}{\nabla [k_1j_1l_1]
\nabla [j_3k_2l_1]\nabla [k_3j_2l_4]\nabla [j_4k_2l_3]}  \nonumber \\
&&\qquad \quad \times \sum_{z_1,z_2,z_3,z_4}\frac{(-1)^{z_2+z_3+z_4}
[j_3+k_2-l_1+z_1]![j_1+k_1-l_1+z_1]!}{[z_1]![z_2]![z_3]![z_4]![l_1+k_2
-j_3-z_1]![j_1-k_1+l_1-z_1]!}  \nonumber \\
&&\qquad \quad \times \frac{[2l_1-z_1]![j_1+l_2-k_3+z_2]![l_2-j_3+k_4+
z_2]!}{[j_1-l_2+k_3-z_2]![j_3+k_4-l_2-z_2]![2l_2+z_2+1]![j_4-k_4+l_4
-z_3]!}  \nonumber \\
&&\qquad \quad \times \frac{[2l_4-z_3]![j_2-l_4+k_3+z_3]!}{[j_2-k_3+l_4
-z_3]![j_4+k_4+l_4-z_3+1]![k_1-j_2+l_3-z_4]!}  \nonumber \\
&&\qquad \quad \times \frac{[2l_3-z_4]![k_2-l_3+j_4+z_4]!}{[k_2+l_3-j_4
-z_4]![k_1+j_2+l_3-z_4+1]![k_1-k_3-l_1+l_2+z_1+z_2]!}  \nonumber \\
&&\qquad \quad \times \frac{[j_3+j_4-l_2+l_4-z_2-z_3]![j_3-j_4+k_1
+k_3-l_2+l_3-z_2-z_4]!}{[j_3-j_4+k_1+k_3-l_1-l_4+z_1+z_3]!
[j_3+j_4-l_1-l_3+z_1+z_4]!}  \nonumber \\
&&\qquad \quad \times \frac{[k_1-k_3+l_3+l_4-z_3-z_4]!}{
[l_1-l_2+l_3+l_4-z_1-z_2-z_3-z_4]!}  \label{dasc} \\
&&\qquad =(-1)^{k_3+k_4-l_1+l_2+l_3-l_4-j_1-j_3}\frac{\nabla [j_3k_2l_1]
\nabla [k_3j_1l_2]\nabla [j_4k_4l_4]\nabla [k_1j_2l_3]}{\nabla [k_1j_1l_1]
\nabla [j_3k_4l_2]\nabla [k_3j_2l_4]\nabla [j_4k_2l_3]}  \nonumber \\
&&\qquad \quad \times \sum_{z_1,z_2,z_3,z_4}\frac{(-1)^{z_1+z_2+z_3}
[2l_1-z_1]![j_1+k_1-l_1+z_1]!}{[z_1]![l_1-k_2+j_3-z_1]![j_1-k_1+l_1-z_1]!
[l_1+k_2+j_3-z_1+1]!}  \nonumber \\
&&\qquad \quad \times \frac{[2l_2-z_2]![k_4+j_3-l_2+z_2]!}{[z_2]!
[k_3+l_2-j_1-z_2]![l_2-j_3+k_4-z_2]![j_1+l_2+k_3-z_2+1]!}  \nonumber \\
&&\qquad \quad \times \frac{[2l_4-z_3]![j_2+k_3-l_4+z_3]!}{[z_3]![z_4]!
[j_4+l_4-k_4-z_3]![j_2-k_3+l_4-z_3]![j_4+k_4+l_4-z_3+1]!}  \nonumber \\
&&\qquad \quad \times \frac{[2j_2-z_4]![j_2+j_4+k_1-k_2-z_4]![j_2+j_4+
k_1+k_2-z_4+1]!}{[k_1+j_2-l_3-z_4]![k_1+j_2+l_3-z_4+1]![j_2+k_3-l_4+z_3
-z_4]!}  \nonumber \\
&&\qquad \quad \times \frac{[l_1+l_2-k_1+k_3-z_1-z_2]![j_3+j_4-k_3-k_1
+l_1+l_4-z_1-z_3]!}{[j_2-j_3+j_4+k_1-l_1+z_1-z_4]![j_2+j_3+j_4-k_3-l_2
+z_2-z_4]!}  \nonumber \\
&&\qquad \quad \times \frac{[l_2+l_4-j_3+j_4-z_2-z_3]!}{%
[l_1+l_2+l_4-k_1-j_2-z_1-z_2-z_3+z_4]!}.  \label{dasd}
\end{eqnarray}
\end{mathletters}
The numerator--denominator distribution of factorials, depending on the 
summation parameters $z_1,z_2,z_3,z_4$, is different in each expression 
(\ref{dasa})--(\ref{dasd}). No single formula exhibites the full
symmetry (\ref{ddabc})--(\ref{ddde}) of the $q$-$12j$-symbol, but 
(\ref{dasa}) is invariant with respect to the transition from the main 
notation to the left array of (\ref{ddde}), as well as under 
transposition
\begin{equation}
\left[ \begin{array}{llll}
j_1 & j_2 & j_3 & j_4 \\ 
l_1 & l_2 & l_3 & l_4 \\ 
k_1 & k_2 & k_3 & k_4
\end{array} \right] _{\!q}=\left[ \begin{array}{rrrr}
k_2 & k_4 & k_1 & k_3 \\ 
l_1 & l_3 & l_2 & l_4 \\ 
j_3 & j_1 & j_4 & j_2
\end{array} \right] _{\!q}, \label{ddf}
\end{equation}
which, in turn, is a composition of symmetry relations 
(\ref{ddabc})--(\ref{ddde}). Expression (\ref{dasb}) is invariant 
with respect to the same symmetry (\ref{ddf}), but (\ref{dasc}) and 
(\ref{dasd}) do not satisfy any symmetry relations. Since all the sums 
in these expressions correspond to the balanced hypergeometric functions, 
the $q$-phases are also trivial.\cite{N89bq}

All the terms in the first sum of (\ref{dasb}) and (\ref{dasc}) are of 
the same sign, as well as in the last sum of (\ref{dasd}). Each separate 
sum corresponds in these expressions to the finite balanced basic 
hypergeometric series $_5F_4[q,1]$ (\ref{fbhs}), which also appeared in 
the elementary overlap coefficients of the definite biorthogonal coupled 
states\cite{Al97} of $u_q(3)$ and SU(3). The summation intervals are 
mainly restricted by 8 [in (\ref{dasa})--(\ref{dasc})], or 7 [in 
(\ref{dasd})] triangle linear combinations of parameters, respectively. 
In addition to correspondence of numerator and denominator factorials, 
determined by Eq.\ (\ref{p4s0}) or (\ref{p4sa}), definite correlation 
between the factorials under summation signs reveals itself in four 
quintuplets of factorials of each expression (\ref{dasa})--(\ref{dasd}), 
depending on the couples of adjacent summation parameters ($z_i$ and 
$z_{i+1}$, where $i=1,2,3$, or $z_1$ and $z_4$), although their expansion 
using the Chu--Vandermonde formulas is not helpful for further 
rearrangement of the generic expressions. 

\subsection{Stretched cases of the $q$-$12j$ coefficients of the second 
kind}

Let us consider the stretched cases of $q$-$12j$ coefficients. For 
definite stretched triangles some summation parameters in 
(\ref{dasa})--(\ref{dasd}) are either fixed (31 times), or expressions 
are partially summable (in the 11 cases) by means of Minton's 
summation formulas (\ref{sMint}) or (\ref{sMinta}) (see Ref.\ 
\onlinecite{GR90}).  One of three remaining sums turns into balanced 
basic hypergeometric series $_4F_3[q,1]$, the rearrangement\cite{GR90} of 
which enables us to transform a $_5F_4[q,1]$ type series into 
$_4F_3[q,1]$ type series, with only the last one remaining of the 
$_5F_4[q,1]$ type. Particularly, for $j_1+l_1=k_1$ with $z_1=0$, the sum 
over $z_2$ in expression (\ref{dasa}) corresponds to a $q$-$6j$ 
coefficient, which may be reexpressed in such a form (using Regge 
symmetry and change of the summation parameter) that the sum over $z_3$ 
also corresponds to a $q$-$6j$ coefficient. Hence, we obtain
\begin{eqnarray}
&&\left[ \begin{array}{cccc}
j_1 & j_2 & j_3 & j_4 \\ 
l_1 & l_2 & l_3 & l_4 \\ 
j_1\!+\!l_1 & k_2 & k_3 & k_4
\end{array}\right] _{\!q}  \nonumber \\
&&\qquad =\frac{(-1)^{j_1+l_2-k_3}\,\nabla [k_1j_2l_3]
\nabla [j_4k_2l_3]}{\nabla [j_1l_2k_3]\nabla [k_4j_3l_2]
\nabla [k_3j_2l_4]\nabla [j_4k_4l_4]\nabla [l_1k_2j_3]}
\left( \frac{[2l_1]![2j_1]!}{[2k_1+1]!}\right) ^{1/2}  \nonumber \\
&&\qquad \quad \times \sum_{z_1,z_3,z_4}\frac{(-1)^{z_1+z_3}
[l_2-j_3+k_4+z_1]![j_1+j_3+k_3-k_4-z_1]!}{[z_1]![z_3]![z_4]!
[l_2+j_3-k_4-z_1]![j_3-j_4-l_1+l_3+z_4-z_1]!}  \nonumber \\
&&\qquad \quad \times \frac{[2j_3-z_1]![j_4+k_4-l_4+z_3]!
[j_2+k_3-l_4+z_3]![2l_4-z_3]!}{[j_2-k_3+l_4-z_3]!
[k_3+j_1+j_3+j_4-l_4-z_1+z_3+1]!} \nonumber \\
&&\qquad \quad \times \frac{[j_2+l_3-k_1+z_4]![k_2+l_3-j_4+z_4]!}{%
[k_4+l_1-l_3+l_4-j_3+z_1-z_3-z_4]![k_3-k_1+l_3-l_4+z_3+z_4]!}  
\nonumber \\
&&\qquad \quad \times \frac{[k_1+k_3+k_4+j_4-l_3-z_4+1]!}{[k_1+j_2-
l_3-z_4]![k_2-l_3+j_4-z_4]![2l_3+z_4+1]!},  \label{dasr}
\end{eqnarray}
which is the composition of two balanced $_4F_3[q,1]$ series and the 
third balanced $_5F_4[q,1]$ series. 

After the summation over $z_3$ of the balanced $_3F_2[q,1]$ series is 
carried out [see Eqs.\ (\ref{sSbal}) and (\ref{sSbbl})] in this 
doubly stretched case of $q$-$12j$ coefficient with $k_1=j_1+l_1$ and 
$l_3=k_1+j_2$ [i.e., for adjacent consecutive stretched triangles in 
graph (\ref{grp2r})], we recognize some $q$-$6j$ coefficients, which may 
also be obtained using the symmetries (\ref{ddabc})--(\ref{ddde}) and the 
defining relations (\ref{dd6a})--(\ref{dd6d}) of the $q$-$12j$ 
coefficients. In this way, we derive following the relation:
\begin{eqnarray}
&&\left[ \begin{array}{cccc}
j_1 & j_2 & j_3 & j_4 \\ 
l_1 & l_2 & k_1\!+\!j_2 & l_4 \\ 
j_1\!+\!l_1 & k_2 & k_3 & k_4
\end{array} \right] _{\!q}  \nonumber \\
&&\qquad =\frac{(-1)^{j_1+l_2-k_3+j_4+k_4+l_4}\,\nabla [l_3j_4k_2]
\nabla [j_1\!+\!j_2,l_2,l_4]}{\nabla [j_1l_2k_3]\nabla [j_2k_3l_4]
\nabla [l_1k_2j_3]\nabla [j_1\!+\!j_2,j_3,j_4]}\left( \frac{[2l_1]!
[2j_1]![2j_2]!}{[2k_1+1][2l_3+1]!}\right) ^{1/2}  \nonumber \\
&&\qquad \quad \times \left\{ \begin{array}{ccc}
j_1\!+\!j_2 & l_4 & l_2 \\ 
k_4 & j_3 & j_4
\end{array} \right\} .  \label{dasrc}
\end{eqnarray}

In the doubly stretched case of the $q$-$12j$ coefficient for 
$j_1=k_1-l_1=l_2-k_3$ [i.e., when the adjacent stretched triangles in 
graph (\ref{grp2r}) are diverging], we obtain from Eq.\ (\ref{dasa}) or 
(\ref{dasc}), and from Eq.\ (\ref{dasb}) with fixed $z_1=z_2=0$, two 
different double sum expressions, each depending on 10 parameters and 
corresponding to the $q$-generalizations of the Kamp\'e de F\'eriet 
\cite{K-F21} function $F_{1:3}^{1:4}$, defined as (\ref{sK}). Each 
separate sum corresponds to the balanced basic hypergeometric 
$_5F_4[q,1]$ series. Again, we may identify the couples of quintuplets of 
factorials under summation signs in the numerator and denominator, each 
depending on the summation parameters $z_3$ and $z_4$.  

Otherwise, in the case of the merging adjacent stretched triangles 
(e.g., for $k_1=j_1+l_1=j_2+l_3$), the straightforwardly derived 
expressions include the triple sums; in particular all three sums in 
(\ref{dasrc}) correspond to the balanced basic hypergeometric 
$_4F_3[q,1]$ series. The $_4F_3[q,1]$ type sum over $z_4$ may be 
rearranged in analogy with expressions for the $q$-$6j$ coefficients 
\cite{JB77,AsST96} into another form (cf.\ Ref.\ \onlinecite{GR90}) in a 
such way that the sum over $z_3$ turns into summable balanced basic 
hypergeometric $_3F_2[q,1]$ series. Hence we obtain the doubly stretched 
$q$-$12j$ coefficient in terms of the double sum:
\begin{eqnarray}
&&\left[ \begin{array}{cccc}
j_1 & j_2 & j_3 & j_4 \\ 
l_1 & l_2 & k_1\!-\!j_2 & l_4 \\ 
j_1\!+\!l_1 & k_2 & k_3 & k_4
\end{array}\right] _{\!q}  \nonumber \\
&&\qquad =\frac{(-1)^{j_1+l_2-k_3}\left( [2l_1]![2j_1]![2j_2]![2l_3]!
\right) ^{1/2}\nabla [k_4l_4j_4]}{\nabla [l_3k_2j_4]\nabla [j_1l_2k_3]
\nabla [k_4j_3l_2]\nabla [j_2k_3l_4]\nabla [l_1k_2j_3]}  \nonumber \\
&&\qquad \quad \times \sum_{z,u}\frac{[l_2-j_3+k_4+z]![l_1+k_2-j_3+z]!
[j_1+j_3+k_3-k_4-z]!}{[z]![l_2+j_3-k_4-z]![k_2-l_1+j_3-z]!
[k_3+k_4-j_1-j_3+z]!}  \nonumber \\
&&\qquad \quad \times \frac{[2j_3-z]!(-1)^{z+u}[j_4-k_4+l_4+u]!
[j_4+l_3-k_2+u]!}{[u]![k_4-j_4+l_4-u]![k_3-k_4-j_2+j_4+u]![2j_4+u+1]!}  
\nonumber \\
&&\qquad \quad \times \frac{[k_3+k_4+j_2-j_4-u]![k_2+l_3+j_4+u+1]!}{%
[j_1+j_2+j_3-j_4-z-u]![l_1+l_3-j_3+j_4+z+u+1]!},  \label{dastm}
\end{eqnarray}
which again depends on 10 parameters and corresponds to the 
$q$-generalization of the Kamp\'e de F\'eriet \cite{K-F21} function 
$F_{1:3}^{1:4}$, with each separate sum corresponding to the balanced 
basic hypergeometric $_5F_4[q,1]$ series. Perhaps this expression is 
related to the above mentioned special case of (\ref{dasa}) with fixed 
$z_1=z_2=0$ and the adjacent diverging stretched triangles of $q$-$12j$ 
coefficient with respect to some composition of the usual and ``mirror 
reflection'' ($j\to -j-1$) symmetries.\cite{JB77,JB65} 

In the doubly stretched case of the $q$-$12j$ coefficient with 
$k_1=j_1+l_1$ and $j_4=k_4+l_4$ [i.e., for antipode stretched triangles 
of graph (\ref{grp2r})], we derive from Eq.\ (\ref{dasr}), with fixed 
$z_3=0$ and $z_4=l_1-l_3-j_3+j_4+z_1$, an expression with a single sum, 
which corresponds to the balanced basic hypergeometric $_6F_5[q,1]$ 
series and depends on 10 parameters: 
\begin{eqnarray}
&&\left[ \begin{array}{cccc}
j_1 & j_2 & j_3 & l_4\!+\!k_4 \\ 
l_1 & l_2 & l_3 & l_4 \\ 
j_1\!+\!l_1 & k_2 & k_3 & k_4
\end{array} \right] _{\!q}  \nonumber \\
&&\qquad =\frac{(-1)^{j_1+l_2-k_3}\,\nabla [k_1j_2l_3]\nabla
[j_4k_2l_3]}{\nabla [j_1l_2k_3]\nabla [k_4j_3l_2]\nabla [l_4k_3j_2]
\nabla [l_1k_2j_3]}\left( \frac{[2l_1]![2j_1]![2l_4]![2k_4]!}{[2k_1+1]!
[2j_4+1]!}\right) ^{1/2} \nonumber \\
&&\qquad \quad \times \sum_{z_1}\frac{(-1)^{z_1}[l_2-j_3+k_4+z_1]!
[j_4-j_3+j_2-j_1+z_1]!}{[z_1]![l_2+j_3-k_4-z_1]![j_1+j_2+j_3-j_4-z_1]!}  
\nonumber \\
&&\qquad \quad \times \frac{[l_1+k_2-j_3+z_1]![2j_3-z_1]!}{%
[k_2+j_3-l_1-z_1]![k_3+k_4-j_1-j_3+z_1]!}  \nonumber \\
&&\qquad \quad \times \frac{[j_1+j_3+k_3-k_4-z_1]!}{%
[l_1-l_3-j_3+j_4+z_1]![l_1+l_3-j_3+j_4+z_1+1]!}.  \label{dasra}
\end{eqnarray}
In the doubly stretched case with $k_1=j_1+l_1$ and $k_4=j_4+l_4$ (again 
for antipode stretched triangles) from (\ref{dasc}), after the summation 
over $z_2$ of the balanced $_3F_2[q,1]$ series (see Appendix B), we obtain 
\begin{eqnarray}
&&\left[ \begin{array}{cccc}
j_1 & j_2 & j_3 & j_4 \\ 
l_1 & l_2 & l_3 & l_4 \\ 
j_1\!+\!l_1 & k_2 & k_3 & j_4\!+\!l_4
\end{array} \right] _{\!q}  \nonumber \\
&&\qquad =\frac{(-1)^{l_1-l_2-k_2+k_4}\,\nabla [k_1j_2l_3]\nabla
[k_4j_3l_2]}{\nabla [j_1l_2k_3]\nabla [l_4k_3j_2]\nabla [j_4l_3k_2]
\nabla [l_1k_2j_3]}\left( \frac{[2l_1]![2j_1]![2l_4]![2j_4]!}{[2k_1+1]!
[2k_4+1]!}\right) ^{1/2}  \nonumber \\
&&\qquad \quad \times \sum_{z_4}\frac{(-1)^{z_4}[j_2+l_3-k_1+z_4]![k_1+
k_2-j_2+j_4-z_4]!}{[z_4]![k_1-j_2+l_3-z_4]![k_2+j_2-k_1-j_4+z_4]!}  
\nonumber \\
&&\qquad \quad \times \frac{[l_4-k_3+j_2+z_4]![j_3-j_4+j_2-j_1+z_4]!}{%
[j_1-j_2+j_3+j_4-z_4]![l_4-l_2-j_1+j_2+z_4]!}  \nonumber \\
&&\qquad \quad \times \frac{[j_2+k_3+l_4+z_4+1]!}{%
[2j_2+z_4+1]![l_2+l_4-j_1+j_2+z_4+1]!}.  \label{dasrd}
\end{eqnarray}
Again the summation corresponds to the balanced $_6F_5[q,1]$ series and 
depends on 10 parameters. Both expressions (\ref{dasra}) and 
(\ref{dasrd}) satisfy some Regge type symmetry relations. 

Consider also the doubly stretched cases of $q$-$12j$ coefficients, the 
case of the remote stretched triangles of graph (\ref{grp2r}) with 
touching angular momenta forming 4-cycles (quadrangles). There are four 
possible different distributions of the summarized angular momenta of the 
stretched triangles: both these momenta may be in the same 4-cycle, 
either (a) adjacent, (b) antiparallel, (c) the first one may be inside 
and the second one outside of the 4-cycle, or (d) both these momenta may 
be outside of the 4-cycle. For $q$-$12j$ coefficients of the ``adjacent'' 
(a) type, e.g., with $k_1=j_1+l_1$ and $l_3=k_2+j_4$, we derive from Eq.\ 
(\ref{dasr}) the relation
\begin{eqnarray}
&&\left[ \begin{array}{cccc}
j_1 & j_2 & j_3 & j_4 \\ 
l_1 & l_2 & k_2\!+\!j_4 & l_4 \\ 
j_1\!+\!l_1 & k_2 & k_3 & k_4
\end{array} \right] _{\!q}  \nonumber \\
&&\qquad =\frac{(-1)^{j_1+l_2-k_3}\,\nabla [l_3j_2k_1]}{\nabla
[j_1l_2k_3]\nabla [k_4j_3l_2]\nabla [k_3j_2l_4]\nabla [j_4k_4l_4]
\nabla [l_1k_2j_3]}\left( \frac{[2l_1]![2j_1]![2k_2]![2j_4]!}{[2k_1+1]!
[2l_3+1]!}\right) ^{1/2}  \nonumber \\
&&\qquad \quad \times \sum_{z_1,z_3}\frac{(-1)^{z_1+z_3}[l_2-j_3+k_4
+z_1]![j_1+j_3+k_3-k_4-z_1]!}{[z_1]![l_2+j_3-k_4-z_1]![j_3+k_2-l_1-z_1]!}  
\nonumber \\
&&\qquad \quad \times \frac{[2j_3-z_1]![2l_4-z_3]![j_4+k_4-l_4+z_3]!}{%
[z_3]![j_2-k_3+l_4-z_3]![k_4+l_1-l_3+l_4-j_3+z_1-z_3]!}  \nonumber \\
&&\qquad \quad \times \frac{[j_2+k_3-l_4+z_3]!}{[k_3-k_1+l_3-l_4+z_3]!
[k_3+j_1+j_3+j_4-l_4-z_1+z_3+1]!}.  \label{dasrb}
\end{eqnarray}
The double sum depends on 9 (from 10 free) parameters and corresponds to 
the $q$-generalization of the Kamp\'e de F\'eriet\cite{K-F21} function 
$F_{1:2}^{1:3}$, defined as (\ref{sK}) with $b_1+b_1^{\prime }=c_1$, and 
each separate sum corresponding to the balanced $_4F_3[q,1]$ or 
$_4\phi _3$ series. Different (i.e.\ not equivalent) expressions of the 
(a) type appear also for $k_3=j_1+l_2$ and $j_2=k_1+l_3$ from Eq.\ 
(\ref{dasb}) and for $j_3=l_1+k_2$ and $k_4=j_4+l_4$ from Eq.\ 
(\ref{dasc}). Furthermore, the doubly stretched $q$-$12j$ coefficients 
of the ``antiparallel'' (b) type, with $k_1=j_1+l_1$ and $k_3=j_2+l_4$, 
expressions (\ref{dasa}), (\ref{dasc}), and (\ref{dasd}) (with fixed 
parameters $z_1$ and $z_3$) also turn into (mutually different) double 
sums, again depending on 9 (from 10 free) parameters and related to the 
$F_{1:2}^{1:3}$ type functions. This is also the case for expression 
(\ref{dasb}) (with fixed $z_1=z_3=0$) for the doubly stretched $q$-$12j$ 
coefficients of the ``inside--outside'' (c) type, with $k_1=j_1+l_1$ 
and $l_4=j_2+k_3$ ([or, expression (\ref{dasc}) with $l_2=j_1+k_3$ and 
$j_2=k_1+l_3$). Finally, expression (\ref{dasa}) with fixed $z_2=z_4=0$ 
and $l_2=j_1+k_3$ and $l_3=k_1+j_2$ again turns into the double sums 
depending on 9 (from 10 free) parameters and related to the 
$F_{1:2}^{1:3}$ type function for the doubly stretched $q$-$12j$ 
coefficients with the both summarized angular momenta of the ``outside'' 
(d) type. These 8 independent expressions should be related to 
(\ref{dasrb}) by means of some compositions of the usual and ``mirror 
reflection'' ($j\to -j-1$) symmetries.\cite{JB77,JB65} Otherwise, 
many special versions of (\ref{dasa})--(\ref{dasd}) with fixed 
$z_i=z_{i+1}=0$ ($i=1,2,3$) or $z_1=z_4=0$ give expressions for the 
doubly stretched $q$-$12j$ coefficients with remote stretched triangles 
in terms of the double sums, related to compositions of the balanced 
$_4F_3[q,1]$ and $_5F_4[q,1]$ series.

Equation (\ref{dasb}) also turns into a single term for 
$l_1+l_2+l_3-l_4=0$ (when the all summation parameters $z_i$ are fixed), 
\begin{eqnarray}
&&\left[ \begin{array}{cccc}
j_1 & j_2 & j_3 & j_4 \\ 
l_1 & l_2 & l_3 & l_1\!+\!l_2\!+\!l_3 \\ 
k_1 & k_2 & k_3 & k_4
\end{array} \right] _{\!q}  \nonumber \\
&&\qquad =\frac{(-1)^{j_1-j_3-k_1+k_2}[2l_1]![2l_2]![2l_3]!
\nabla [l_4k_3j_2]\nabla [l_4j_4k_4]}{[2l_4+1]!\nabla [l_1j_1k_1]
\nabla [l_1k_2j_3]\nabla [l_2j_1k_3]\nabla [l_2j_3k_4]\nabla [l_3j_2k_1]
\nabla [l_3k_2j_4]}.  \label{dvosb}
\end{eqnarray}
For this special $q$-$12j$ coefficient [as well as in (\ref{dasc}) for 
$l_1-l_2+l_3+l_4=0$ and in (\ref{dasd}) for $l_1+l_2-l_3+l_4=0$], four 
linearly dependent angular momenta appear as disconnected in certain 
positions on a Hamilton line of graph (\ref{grp2r}). Actually, the single 
term expression of this virtually stretched case appears in accordance 
with symmetries (\ref{ddabc})--(\ref{ddde}) from expansion 
(\ref{dd6a}) with $j_1+j_2-j_3+j_4=0$ and fixed $x=j_3-j_1=j_2+j_4$.

\section{EXPRESSIONS FOR 12$\lowercase{j}$ COEFFICIENTS OF THE FIRST KIND}

\subsection{Generic properties}

Next we consider the rearrangement of expressions for the $q$-$12j$ 
coefficients of the first kind\cite{JLV60,JB77,JH54,O-S54} whose graphs 
are not planar: 
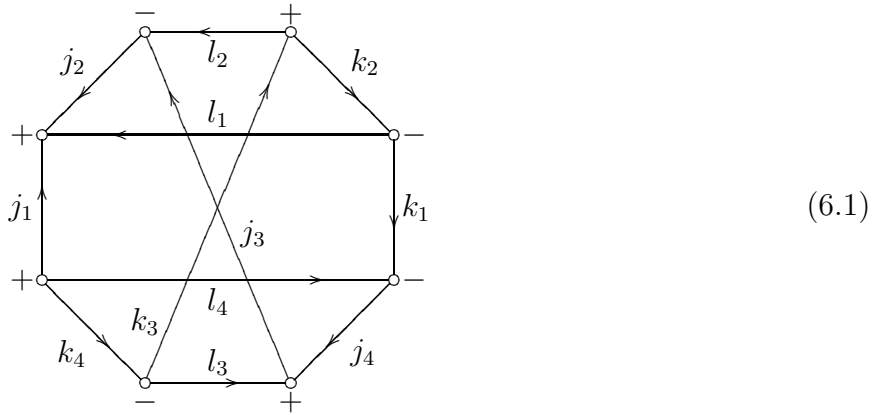
\begin{equation}
\xy
\xygraph{!{/r6pc/:}
[] 
!P8"A"{~>-\cir<2pt>{}}
"A1"*+!L{-}
"A2"*+!D{+}
"A3"*+!D{-}
"A4"*+!R{+}
"A5"*+!R{+}
"A6"*+!U{-}
"A7"*+!U{+}
"A8"*+!L{-}
"A2" (? :@{-}|(0.65)@{>}^{\dpst k_2} "A1")
"A1" (? :@{-}|(0.65)@{>}^{\dpst k_1} "A8")
"A2" (? :@{-}|(0.65)@{>}^{\dpst l_2} "A3")
"A3" (? :@{-}|(0.65)@{>}_{\dpst j_2} "A4")
"A5" (? :@{-}|(0.65)@{>}^{\dpst j_1} "A4")
"A5" (? :@{-}|(0.65)@{>}_{\dpst k_4} "A6")
"A8" (? :@{-}|(0.65)@{>}^{\dpst j_4} "A7")
"A6" (? :@{-}|(0.65)@{>}^{\dpst l_3} "A7")
"A4" ("A1" :@{-}|(0.8)@{>}_{\dpst l_1} ?)
"A8" ("A5" :@{-}|(0.8)@{>}_{\dpst l_4} ?)
"A3" ("A7" :@{-}|(0.83)@{>}_(0.4){\dpst j_3} ?)
"A2" ("A6" :@{-}|(0.85)@{>}^(0.15){\dpst k_3} ?)
}\endxy  \label{grp1r}
\end{equation}
(include some braiding\cite{N90}). These coefficients satisfy 16 
symmetries, generated by the following substitutions:
\begin{mathletters}
\begin{eqnarray}
\left\{ \begin{array}{llllllll}
j_1\! &  & \!j_2\! &  & \!j_3\! &  & \!j_4\! &  \\ 
& \!l_1\! &  & \!l_2\! &  & \!l_3\! &  & \!l_4\! \\ 
k_1\! &  & \!k_2\! &  & \!k_3 \! &  & \!k_4\! &
\end{array} \right\} _{\!\!q}&=&
\left\{ \begin{array}{llllllll}
j_2\! &  & \!j_3\! &  & \!j_4\! &  & \!k_1\! &  \\ 
& \!l_2\! &  & \!l_3\! &  & \!l_4\! &  & \!l_1\! \\ 
k_2\! &  & \!k_3\! &  & \!k_4 \! &  & \!j_1\! & 
\end{array} \right\} _{\!\!q}  \label{dpa} \\
&=&\left\{ \begin{array}{llllllll}
k_1\! &  & \!j_4\! &  & \!j_3\! &  & \!j_2\! &  \\ 
& \!l_4\! &  & \!l_3\! &  & \!l_2\! &  & \!l_1\! \\ 
j_1\! &  & \!k_4\! &  & \!k_3 \! &  & \!k_2\! &  
\end{array} \right\} _{\!\!q}.  \label{dpb}
\end{eqnarray}
\end{mathletters}
There expression\cite{JLV60,JB77,JB65,N90} in terms of the factorized 
four differently rearranged $q$-$6j$ coefficients is
\begin{mathletters}
\begin{eqnarray}
&&\left\{ \begin{array}{llllllll}
j_1\! &  & \!j_2\! &  & \!j_3\! &  & \!j_4\! &  \\ 
& \!l_1\! &  & \!l_2\! &  & \!l_3\! &  & \!l_4\! \\ 
k_1\! &  & \!k_2\! &  & \!k_3 \! &  & \!k_4\! &
\end{array} \right\} _{\!\!q}  \nonumber \\
&&\qquad =\sum_x(2x+1)(-1)^{R_4-x}q^{x(x+1)+Z_{j_1j_2j_3j_4}+
Z_{k_1k_2k_3k_4}}  \nonumber \\
&&\qquad \quad\times \left\{ \begin{array}{lll}
j_1 & j_2 & l_1 \\ 
k_2 & k_1 & x
\end{array} \right\} _{\!\!q}\left\{ \begin{array}{lll}
j_3 & k_3 & x \\ 
k_2 & j_2 & l_2
\end{array} \right\} _{\!\!q}\left\{ \begin{array}{lll}
j_3 & j_4 & l_3 \\ 
k_4 & k_3 & x
\end{array} \right\} _{\!\!q}\left\{ \begin{array}{lll}
k_1 & j_1 & x \\ 
k_4 & j_4 & l_4
\end{array} \right\} _{\!\!q}  \label{dvpja} \\
&&\qquad =\sum_{x}(2x+1)(-1)^{R_4-x} q^{x(x+1)+Z_{j_1j_2j_3j_4}+
Z_{k_1k_2k_3k_4}}  \nonumber \\
&&\qquad \quad\times \left\{ \begin{array}{lll}
j_1 & j_2 & l_1 \\ 
k_2 & k_1 & x
\end{array} \right\} _{\!\!q}\left\{ \begin{array}{lll}
k_2 & x & j_2 \\ 
j_3 & l_2 & k_3
\end{array} \right\} _{\!\!q}\left\{ \begin{array}{lll}
j_3 & k_3 & x \\ 
k_4 & j_4 & l_3
\end{array} \right\} _{\!\!q}\left\{ \begin{array}{lll}
k_1 & x & j_1 \\ 
k_4 & l_4 & j_4
\end{array} \right\} _{\!\!q},  \label{dvpjb}
\end{eqnarray}
\end{mathletters}
where the triangular conditions are to be satisfied by all the triplets 
of the nearest neighbors such as $l_i,j_i,j_{i+1}$, or $l_i,k_i,k_{i+1}$ 
($i=1,2,3$), or $l_4,j_1,k_4$, or  $l_4,k_1,j_4$, respectively. 

After using (\ref{f6jbb}) for the $q$-$6j$ coefficients with the 
summation parameter $x$ in the right lower position, Eq.\ (\ref{f6jb}) 
with inverted summation parameter for the $q$-$6j$ coefficients with $x$ 
in the middle column, and (\ref{f6jb}) directly in the remaining cases, 
the depending on the summation parameter $x$ asymmetric triangle 
coefficients [distributed separately in the numerators or denominators 
of each $q$-$6j$ coefficient in expansions (\ref{dvpja}) and (\ref{dvpjb})] 
cancel, with the exception of the factors
\[
\frac{\nabla [j_1k_1x]}{\nabla [k_1j_1x]}=
\frac{[j_1-k_1+x]!}{[k_1-j_1+x]!}.
\]
Then, the sums over $x$ correspond to the $q$-generalization of the very 
well-poised classical hypergeometric $_6F_5(-1)$ series (resembling the 
basic hypergeometric $_7\phi _6$ series) and may be rearranged into the 
$_3\phi _2$ or $_3F_2[q,x]$ type series using the following two formulas:
\begin{mathletters}
\begin{eqnarray}
&&\sum_{j}\frac{(-1)^{p_2+j+1}q^{j(j+1)}[2j+1][j-p_1-1]![j-p_2-1]!
[j-p_3-1]!}{[p_1+j+1]![p_2+j+1]![p_3+j+1]![p_4-j]![p_4+j+1]![p_5-j]!
[p_5+j+1]!}  \nonumber \\
&&\qquad =\frac{q^{-(p_4+1)(p_5+1)-p_2(p_4+p_5+1)}[-p_1-p_3-2]!}
{[p_1+p_4+1]![p_2+p_5+1]![p_3+p_4+1]!}  \nonumber \\
&&\qquad \quad \times \sum_{u}\frac{(-1)^uq^{u(p_2+p_5+1)}[p_4-p_3-1-u]!
[p_4-p_1-1-u]!}{[u]![-p_1-p_3-2-u]![p_2+p_4+1-u]![p_4+p_5+1-u]!}  
\label{p5sra}
\end{eqnarray}
with parameters 
\[
\begin{array}{c}
p_1=k_1-j_1-1,\;\;p_4=j_1+k_2-l_1+z_1,\;\;p_5=j_3+k_4-l_3+z_3, \\ 
p_2=l_2-k_2-j_3+z_2-1,\;\;p_3=l_4-k_1-k_4+z_4-1;
\end{array}
\]
\begin{eqnarray}
&&\sum_{j}\frac{q^{j(j+1)}[2j+1][j-p_1-1]![j-p_2-1]![j-p_3-1]!
[j-p_5-1]!}{[p_1+j+1]![p_2+j+1]![p_3+j+1]![p_4-j]![p_4+j+1]![p_5+j+1]!} 
\nonumber \\
&&\qquad =q^{-(p_4+1)(p_5+1)-p_2(p_4+p_5+1)}
\frac{[-p_1-p_3-2]![-p_2-p_5-2]!}{[p_1+p_4+1]![p_3+p_4+1]!}  \nonumber \\
&&\qquad \quad \times \sum_{u}\frac{(-1)^{u}q^{u(p_2+p_5+1)}
[p_4-p_3-1-u]![p_4-p_1-1-u]!}{[u]![p_4+p_5+1-u]![p_2+p_4+1-u]!
[-p_1-p_3-2-u]!}  \label{p5srb}
\end{eqnarray}
\end{mathletters}
with parameters 
\[
\begin{array}{c}
p_1=k_1-j_1-1,\;\;p_2=j_3-k_3-z_2-1,\;\;p_3=j_1-k_1-z_4-1, \\ 
p_4=j_1+k_2-l_1+z_1,\;\;p_5=l_3-j_3-k_4+z_3-1.
\end{array}
\]
Equation (\ref{p5sra}) corresponds to (5.5) [or (5.6), when $q=1$] of 
Ref.\ \onlinecite{Al97}, with the r.h.s.\ replaced using less symmetric 
expression\cite{JB77,STK91a} instead of the most symmetric (Van der 
Waerden) expression\cite{Ed57,JLV60,JB77,N89c,N89q,GKK90,Ru90} for the 
Clebsch--Gordan coefficients of SU(2) and $u_q(2)$ (cf.\ also Refs.\ 
\onlinecite{A-NS96,Al74}), when (\ref{p5srb}) is derived from 
(\ref{p5sra}) using the analytical continuation technique.

\subsection{General expressions with five sums}

Substituting the summation parameter $u$, which appeared after using 
(\ref{p5sra}) and (\ref{p5srb}) in (\ref{dvpja}) and (\ref{dvpjb}), by 
$u+z_1$, we obtain the following expressions for the $q$-$12j$ 
coefficients of the first kind: 
\begin{mathletters}
\begin{eqnarray}
&&\left\{ \begin{array}{llllllll}
j_1\! &  & \!j_2\! &  & \!j_3\! &  & \!j_4\! &  \\ 
& \!l_1\! &  & \!l_2\! &  & \!l_3\! &  & \!l_4\! \\ 
k_1\! &  & \!k_2\! &  & \!k_3 \! &  & \!k_4\! & 
\end{array}\right\} _{\!\!q}  \nonumber \\
&&\qquad =(-1)^{k_1+k_2-k_3-j_2+j_3+j_4-l_4}\frac{\nabla [j_3j_2l_2]
\nabla [k_2k_3l_2]\nabla [k_1j_4l_4]\nabla [k_4j_1l_4]}{\nabla 
[j_1j_2l_1]\nabla [k_2k_1l_1]\nabla [j_3j_4l_3]\nabla [k_4k_3l_3]}  
\nonumber \\
&&\qquad \quad \times q^{-(j_1+k_2-l_1+1)(k_4-l_3+l_2-k_2)-
(l_2-k_2-j_3-1)(j_3+k_4-l_3)+Z_{j_1j_2j_3j_4}+Z_{k_1k_2k_3k_4}}  
\nonumber \\
&&\qquad \quad \times \sum_{z_1,z_2,z_3,z_4,u}\frac{(-1)^{z_3+z_4+u}
[j_1+j_2-l_1+z_1]![2l_1-z_1]!}{[z_1]![z_2]![z_3]![z_4]![k_1-k_2+l_1-z_1]!
[l_1-j_1+j_2-z_1]!}  \nonumber \\
&&\qquad \quad \times \frac{q^{-z_2(j_1+j_3+k_2+k_4-l_1-l_3+z_3+1)}
[l_2+j_2-j_3+z_2]![l_2+k_3-k_2+z_2]!}{[j_2+j_3-l_2-z_2]!
[k_2+k_3-l_2-z_2]![2l_2+z_2+1]![2l_4+z_4+1]!}  \nonumber \\
&&\qquad \quad \times \frac{[k_3+k_4-l_3+z_3]![j_3+j_4-l_3+z_3]![2l_3-
z_3]!}{[k_3+l_3-k_4-z_3]![l_3-j_3+j_4-z_3]![k_1+j_4-l_4-z_4]!}  
\nonumber \\
&&\qquad \quad \times \frac{q^{z_3(l_1-l_2-j_1+j_3)}[j_4+l_4-k_1+z_4]!
[j_1+l_4-k_4+z_4]!}{[l_2-l_3-k_2+k_4+z_2+z_3]!
[j_1-k_1+k_2-k_4-l_1+l_4+z_1+z_4]!}  \nonumber \\
&&\qquad \quad \times \frac{q^{u(k_4-k_2-l_3+l_2+z_2+z_3)}[j_1+k_1+k_2
+k_4-l_1-l_4-z_4-u]!}{[u+z_1]![j_1+k_4-l_4-z_1-z_4-u]![j_1-j_3-l_1+
l_2+z_2-u]!}  \nonumber \\
&&\qquad \quad \times \frac{[2j_1-k_1+k_2-l_1-u]!}{[j_1+j_3+k_2+
k_4-l_1-l_3+z_3-u+1]!}  \label{dpsa} \\
&&\qquad =(-1)^{j_1+j_3+j_4+k_2-k_3-k_4+l_2-l_3+l_4}
\frac{\nabla [j_3j_4l_3]\nabla [k_4k_3l_3]\nabla [k_2k_3l_2]
\nabla [k_1j_4l_4]}{\nabla [j_1j_2l_1]\nabla [k_2k_1l_1]\nabla [j_3j_2l_2]
\nabla [k_4j_1l_4]}  \nonumber \\
&&\qquad \quad \times q^{(j_1+k_2-l_1+1)(l_2+k_2-l_3+k_4+1)
-(l_2+k_2-j_3+1)(j_3+k_4-l_3+1)+Z_{j_1j_2j_3j_4}+Z_{k_1k_2k_3k_4}}  
\nonumber \\
&&\qquad \quad \times \sum_{z_1,z_2,z_3,z_4,u}\frac{(-1)^{z_2+z_3+z_4+u}
[j_1+j_2-l_1+z_1]![2l_1-z_1]!}{[z_1]![z_2]![z_3]![k_1-k_2+l_1-z_1]!
[l_1-j_1+j_2-z_1]!}  \nonumber \\
&&\qquad \quad \times \frac{q^{-z_2(j_1-j_3+k_2-k_4-l_1+l_3+z_3)}
[j_2+j_3-l_2+z_2]![2l_2-z_2]!}{[l_2+k_2-k_3-z_2]![j_2-j_3+l_2-z_2]!
[l_2+k_2+k_3-z_2+1]!}  \nonumber \\
&&\qquad \quad \times \frac{q^{-z_3(j_1+j_3-l_1-l_2)}[k_3-k_4+l_3+z_3]!
[j_4-j_3+l_3+z_3]!}{[k_3+k_4-l_3-z_3]![j_3+j_4-l_3-z_3]!
[2l_3+z_3+1]!}  \nonumber \\
&&\qquad \quad \times \frac{[2l_4-z_4]![j_1+k_4-l_4+z_4]![k_2+k_4+l_2
-l_3-z_2-z_3]!}{[z_4]![k_1+l_4-j_4-z_4]![j_1-k_1+k_2+k_4-l_1-l_4+z_1
+z_4]!}  \nonumber \\
&&\qquad \quad \times \frac{q^{-u(k_2+k_4+l_2-l_3-z_2-z_3+1)}[2j_1-k_1
+k_2-l_1-u]!}{[k_1+l_4+j_4-z_4+1]![u+z_1]![j_1-k_4+l_4-z_1-z_4-u]!}  
\nonumber \\
&&\qquad \quad \times \frac{[j_1+k_1+k_2-k_4-l_1+l_4-z_4-u]!}{[j_1+j_3
-l_1-l_2+z_2-u]![j_1-j_3+k_2-k_4-l_1+l_3+z_3-u]!}.  \label{dpsb}
\end{eqnarray}
\end{mathletters}
Each of expressions (\ref{dpsa}) and (\ref{dpsb}) includes 5 summations, 
with four separate sums (over $z_1,z_2,z_3$, and $z_4$) corresponding to 
the finite (balanced in the first and last cases) basic hypergeometric 
series $_4\phi _3$ or $_4F_3[\cdot \cdot \cdot ;q,1]$ [cf.\ definition 
(\ref{fbhs})], and the fifth sum (over $u+z_1$) corresponding to the 
finite hypergeometric series $_3\phi _2$ or $_3F_2(1)$, related in the 
case of (\ref{dpsa}) to the Clebsch--Gordan coefficients of $u_q(2)$ or 
SU(2). However, it is impossible to rearrange all 5 sums together into 
standard basic hypergeometric series $_{p+1}\phi _p$. Some correlation 
between the factorials under the summation signs reveals itself in two 
quintuplets of factorials of each expression (\ref{dpsa}) and 
(\ref{dpsb}), depending on the couples of summation parameters $z_2,z_3$ 
and $z_1,z_4$. Definite correspondences may be observed in Eqs.\ 
(\ref{dpsa}) and (\ref{dpsb}) between the $q$-phase structure and three 
particular factorial arguments, depending on the couples of summation 
parameters $z_2,z_3$, and $u$, respectively [as well as in $q$-$9j$ 
coefficients (\ref{trsa})--(\ref{trse}) between the $q$-phases and three 
factorial arguments, depending on the couples of summation parameters 
$z_1,z_2$, and $z_3$]. The summation intervals in (\ref{dpsa}) and 
(\ref{dpsb}) are mainly restricted by 8 triangle linear combinations of 
parameters, respectively, but in the stretched cases only (\ref{dpsa}) 
for $k_4=l_4-j_1$ and (\ref{dpsb}) for $k_4=l_4+j_1$ (with $z_4=u+z_1=0$ 
in the both cases) turn into the triple sums. 

\subsection{Stretched cases of the $q$-$12j$ coefficients of the first 
kind}

When the total (maximal) angular momentum in a stretched triangle of the 
$q$-$12j$ coefficient of the first kind corresponds to a crossbar of the 
M\"{o}bius strip (\ref{grp1r}) [in the middle row of standard array 
(\ref{dpa}), e.g., for $l_4=k_4+j_1$], we obtain the following 
expression:
\begin{eqnarray}
&&\left\{ \begin{array}{cccccccc}
j_1\! &  & \!j_2\! &  & \!j_3\! &  & \!j_4\! &  \\ 
& \!l_1\! &  & \!l_2\! &  & \!l_3\! &  & \!k_4\!+\!j_1\! \\ 
k_1\! &  & \!k_2\! &  & \!k_3 \! &  & \!k_4\! & 
\end{array}\right\} _{\!\!q}  \nonumber \\
&&\qquad =(-1)^{k_1+k_2-k_3-j_2+j_3+j_4-l_4}\frac{\nabla [j_3j_2l_2]
\nabla [k_2k_3l_2]\nabla [l_4j_4k_1]}{\nabla [j_1j_2l_1]\nabla [k_2k_1l_1]
\nabla [j_3j_4l_3]\nabla [k_4k_3l_3]}\left( \frac{[2k_4]!
[2j_1]!}{[2l_4+1]!}\right) ^{1/2}  \nonumber \\
&&\qquad \quad \times q^{-(j_1+k_2-l_1+1)(k_4-k_2+l_2-l_3)
-(l_2-k_2-j_3-1)(j_3+k_4-l_3)+Z_{j_1j_2j_3j_4}+Z_{k_1k_2k_3k_4}}  
\nonumber \\
&&\qquad \quad \times \sum_{z_1,z_2,z_3}\frac{(-1)^{z_1+z_3}
[j_1+j_2-l_1+z_1]![k_1+k_2-l_1+z_1]![2l_1-z_1]!}{[z_1]![z_2]![z_3]!
[k_1-k_2+l_1-z_1]![l_1-j_1+j_2-z_1]!}  \nonumber \\
&&\qquad \quad \times \frac{[l_2+j_2-j_3+z_2]![l_2+k_3-k_2+z_2]!}{%
[j_2+j_3-l_2-z_2]![k_2+k_3-l_2-z_2]![2l_2+z_2+1]!}  \nonumber \\
&&\qquad \quad \times \frac{[k_3+k_4-l_3+z_3]![j_3+j_4-l_3+z_3]![2l_3-
z_3]!}{[k_3+l_3-k_4-z_3]![l_3-j_3+j_4-z_3]![j_1-j_3-l_1+l_2+z_1+z_2]!}  
\nonumber \\
&&\qquad \quad \times \frac{q^{z_3(l_1-l_2-j_1+j_3)-z_2(j_3+k_2-l_1-l_3
+l_4+z_3+1)-z_1(l_2-l_3-k_2+k_4+z_2+z_3)}}{[l_2-l_3-k_2+k_4+z_2+z_3]!
[j_3+k_2-l_1-l_3+l_4+z_1+z_3+1]!}.  \label{dpsta} 
\end{eqnarray}
This special case of (\ref{dpsa}) with 3 separate sums corresponds 
to the finite basic hypergeometric series $_4F_3[\cdot \cdot \cdot ;q,x]$. 
Expression (\ref{dpsta}) does not simplify noticeably for two adjacent 
merging stretched triangles (with $l_4=k_4+j_1=k_1+j_4$) in the same 
$q$-$6j$ coefficient of expansion (\ref{dvpja}), but one of the sums 
turns into a $_3F_2[\cdot \cdot \cdot ;q,x]$ series for two adjacent 
diverging stretched triangles (e.g., with $j_1=l_4-k_4=l_1-j_2$, or with 
$j_1=l_4-k_4=j_2-j_1$). The Chu--Vandermonde summation formula 
(\ref{sCVc}) may be used in Eq.\ (\ref{dpsta}) for $j_1=0,\ l_1=j_2,\ 
k_4=l_4$ and for $k_4=0,\ l_4=j_1,\ k_3=l_3$. In these cases, the couples 
of the $q$-$6j$ coefficients appear in accordance with Eq.\ (33.21) of 
Ref.\ \onlinecite{JB77}, as well as the consequences of expansions 
(\ref{dvpja})--(\ref{dvpjb}) for fixed $x$.

When the total angular momentum in a stretched triangle of the $q$-$12j$ 
coefficient of the first kind is located along the M\"{o}bius strip 
(\ref{grp1r}) (e.g., for $k_4=l_4+j_1$), we obtain from (\ref{dpsb}), 
after change of summation parameter $z_2\to l_2+k_2-k_3-z_2$, 
the following expression:
\begin{eqnarray}
&&\left\{ \begin{array}{cccccccc}
j_1\! &  & \!j_2\! &  & \!j_3\! &  & \!j_4\! &  \\ 
& \!l_1\! &  & \!l_2\! &  & \!l_3\! &  & \!l_4\! \\ 
k_1\! &  & \!k_2\! &  & \!k_3 \! &  & \!l_4\!+\!j_1\! & 
\end{array}\right\} _{\!\!q}  \nonumber \\
&&\qquad =\frac{(-1)^{j_3+j_4-l_3+l_2+k_2-k_3}\nabla [j_3j_4l_3]
\nabla [k_4k_3l_3]\nabla [k_2k_3l_2]}{\nabla [j_1j_2l_1]\nabla [k_2k_1l_1]
\nabla [j_3j_2l_2]\nabla [l_4k_1j_4]}\left( \frac{[2j_1]![2l_4]!}{%
[2k_4+1]!}\right) ^{1/2}   \nonumber \\
&&\qquad \quad \times q^{(j_1+k_2-l_1+1)(k_2+k_4+l_2-l_3+1)-(l_2+k_2-j_3
+1)(j_3-l_3+k_4+1)+Z_{j_1j_2j_3j_4}+Z_{k_1k_2k_3k_4}}  \nonumber \\
&&\qquad \quad \times \sum_{z_1,z_2,z_3}\frac{[j_1+j_2-l_1+z_1]!
[k_1+k_2-l_1+z_1]![2l_1-z_1]!}{[z_1]![z_2]![z_3]![k_1-k_2+l_1-z_1]!
[l_1-j_1+j_2-z_1]!}  \nonumber \\
&&\qquad \quad \times \frac{(-1)^{z_1+z_2+z_3}[j_2+j_3-l_2+z_2]!
[2l_2-z_2]!}{[l_2+k_2-k_3-z_2]![j_2-j_3+l_2-z_2]![l_2+k_2+k_3-z_2+1]!}  
\nonumber \\
&&\qquad \quad \times \frac{q^{z_1(k_2+k_4+l_2-l_3-z_2-z_3+1)}
[k_3-k_4+l_3+z_3]![j_4-j_3+l_3+z_3]!}{[k_3+k_4-l_3-z_3]![j_3+j_4-l_3-z_3]!
[2l_3+z_3+1]!}  \nonumber \\
&&\qquad \quad \times \frac{q^{-z_2(k_2-j_3-l_1+l_3-l_4+z_3)-z_3
(j_1+j_3-l_1-l_2)}[k_2+k_4+l_2-l_3-z_2-z_3]!}{[j_1+j_3-l_1-l_2+z_1+z_2]!
[k_2-j_3-l_1+l_3-l_4+z_1+z_3]!}.  \label{dpstb}
\end{eqnarray}
Expression (\ref{dpstb}) does not simplify for two couples of adjacent 
diverging stretched triangles [with $l_4=k_4-j_1=k_1-j_4$, or with 
$l_4=k_4-j_1=j_4-k_1$, respectively, again in the same $q$-$6j$ 
coefficient of expansion (\ref{dvpjb})], but for $l_4=0,\ k_4=j_1,\ 
j_4=k_1$ (after substitution $q\to q^{-1}$) it corresponds to the general 
expression (\ref{trsa}) for the $q$-$9j$ coefficients, in accordance with 
Eq.\ (33.20) of Ref.\ \onlinecite{JB77}. The triple sums in 
(\ref{dpsta}) and (\ref{dpstb}) resemble expressions for the $q$-$9j$ 
coefficients and definite correspondences may be observed between the 
$q$-phases and three factorial arguments, depending on the couples of 
summation parameters $z_1,z_2$, and $z_3$, respectively. Expansion of the 
present couples of the factorial quintuplets (depending on parameters 
$z_1,z_2$ and $z_2,z_3$, respectively) using the Chu--Vandermonde 
summation formulas enables one to perform the summation over $z_2$ and 
thus obtain expressions for the stretched $q$-$12j$ coefficients of the 
first kind as fourfold sums, related to compositions of 
$_3F_2[\cdot \cdot \cdot ;q,x]$ series.

Again, one of the sums in (\ref{dpstb}) turns into a 
$_3F_2[\cdot \cdot \cdot ;q,x]$ series for two adjacent diverging 
stretched triangles in two adjacent $q$-$6j$ coefficients of expansion 
(\ref{dvpja}) (e.g., with $j_1=k_4-l_4=l_1-j_2$, or with 
$j_1=k_4-l_4=j_2-l_1$), as well as for two adjacent merging stretched 
triangles (e.g., with $k_4=l_4+j_1=k_3+l_3$). The possible rearrangement 
of the $_3F_2[\cdot \cdot \cdot ;q,x]$ series is not helpful for reducing 
the remaining $_4F_3[\cdot \cdot \cdot ;q,x]$ series.

The doubly stretched $q$-$12j$ coefficients of the first kind with the 
adjacent consecutive stretched triangles may be expressed for 
$k_4=l_4+j_1$ and $j_1=l_1+j_2$ or $l_3=k_3+k_4$ by means of 
(\ref{dpstb}), as well as for $l_4=k_4+j_1$ and $j_1=l_1+j_2$ or 
$k_4=k_3+l_3$ by means of (\ref{dpsta}), as the double sums, related 
to the stretched $q$-$9j$ coefficients, respectively, of the type 
(\ref{trsa}) or (\ref{trsb}) (with some ``reflected'' parameters in the 
last cases). In this way we derive, for $k_4=l_4+j_1$ and $j_1=l_1+j_2$ 
[comparing (\ref{dpstb}) and (\ref{trsa})], the following relation: 
\begin{eqnarray}
&&\left\{ \begin{array}{cccccccc}
\!l_1\!+\!j_2\! &  & \!j_2\! &  & \!j_3\! &  & \!j_4\! &  \\ 
& \!l_1\! &  & \!l_2\! &  & \!l_3\! &  & \!l_4\! \\ 
k_1\! &  & \!k_2\! &  & \!k_3 \! &  & \!l_4\!+\!j_1\! & 
\end{array}
\right\} _{\!\!q}  \nonumber \\
&&\qquad =(-1)^{j_3+j_4-l_3+l_2+k_2-k_3}q^{Z_{j_1j_2j_3j_4}+
Z_{k_1k_2k_3k_4}-Z_{l_2,l_3,l_1\!+\!l_4}}\frac{\nabla [l_1\!+\!l_4,j_4,
k_2]}{\nabla [l_1k_1k_2]\nabla [l_4k_1j_4]}   \nonumber \\
&&\qquad \quad \times \left( \frac{[2l_1]![2l_4]!}{
[2l_1+2l_4]![2j_1+1]}\right) ^{1/2}\left\{ \begin{array}{ccc}
k_4 & k_3 & l_3 \\ 
j_2 & l_2 & j_3 \\ 
\!l_1\!+\!l_4\! & k_2 & j_4
\end{array}\right\} _{\!\!q^{-1}}.  \label{dvstb1}
\end{eqnarray}
Similarly, for $k_4=l_4+j_1$ and $l_3=k_3+k_4$ we obtain
\begin{eqnarray}
&&\left\{ \begin{array}{cccccccc}
\!j_1\! &  & \!j_2\! &  & \!j_3\! &  & \!j_4\! &  \\ 
& \!l_1\! &  & \!l_2\! &  & \!k_3\!+\!k_4\! &  & \!l_4\! \\ 
k_1\! &  & \!k_2\! &  & \!k_3 \! &  & \!l_4\!+\!j_1\! & 
\end{array}\right\} _{\!\!q}  \nonumber \\
&&\qquad =q^{Z_{j_1j_2j_3j_4}+Z_{k_1k_2k_3k_4}-Z_{l_2,j_1\!+\!k_3,l_1}}
\frac{(-1)^{j_3+j_4-l_4-l_2-k_2-j_2+l_1}\nabla [l_3j_3j_4]}{\nabla 
[l_4k_1j_4]\nabla [j_1\!+\!k_3,j_3,k_1]}   \nonumber \\
&&\qquad \quad \times \left( \frac{[2l_4]![2j_1+2k_3+1]!}{[2l_3+1]!
[2k_4+1]}\right) ^{1/2}\left\{ \begin{array}{ccc}
j_1 & k_3 & \!j_1\!+\!k_3\! \\ 
j_2 & l_2 & j_3 \\ 
l_1 & k_2 & k_1
\end{array}\right\} _{\!\!q^{-1}}.  \label{dvstb2}
\end{eqnarray}

We remark that  the doubly stretched $q$-$12j$ coefficients of the first 
kind for $l_4=k_4+j_1$ and $l_2=j_2+j_3$ or $l_2=k_2+k_3$ 
[expressed by means of (\ref{dpsta})], as obtained from (\ref{dpsta}),
as well as for $k_4=l_4+j_1$ and $j_3=l_2+j_2$ or $k_3=l_2+k_2$ as 
obtained from (\ref{dpstb}), turn into double sums equivalent to 
compositions of $_4F_3[\cdot \cdot \cdot ;q,x]$ and 
$_3F_2[\cdot \cdot \cdot ;q,x]$ series. However, although the 
$_3F_2[\cdot \cdot \cdot ;q,x]$ series may be rearranged into other 
forms, these double sums are not equivalent to any composition of two 
generic $_3F_2[\cdot \cdot \cdot ;q,x]$ series and, moreover, they are 
not related to any $q$-$9j$ coefficient. Minton's summation formula 
(\ref{sMint}) (see Ref.\ \onlinecite{GR90}) may be helpful in the 
analogical position of the stretched triangles in (\ref{dpsta}) for 
$l_4=k_4+j_1$ and $j_2=l_2+j_3$ or $k_3=k_2+l_2$, as well as in 
(\ref{dpstb}) for $k_4=l_4+j_1$ and $l_2=k_2+k_3$. Otherwise, the 
$_3F_2[\cdot \cdot \cdot ;q,x]$ series appearing in the triple sums, 
which remain in Eq.\ (\ref{dpsta}) for $j_3=j_2+l_2$ or $k_2=l_2+k_3$, 
and in Eq.\ (\ref{dpstb}) for $j_2=l_2+j_3$ or $k_2=l_2+k_3$, may be 
rearranged, but it is more reasonable in such a case to use the different 
expansions of the type (\ref{dvpja})--(\ref{dvpjb}) with inserted 
stretched $q$-$6j$ coefficients and adapted [e.g., expansion 
(\ref{dvpjb}) for $l_4=k_4+j_1$ and $j_3=j_2+l_2$ or expansion 
(\ref{dvpja}) for $k_4=l_4+j_1$ and $j_2=l_2+j_3$, with transposed the 
middle and the right columns of the third $q$-$6j$ coefficients in the 
both cases] for summation by means of (\ref{smDqb}) or (\ref{smDqc}).  

Furthermore, expression (\ref{dpsta}) turns into double sums equivalent 
to compositions of two $_4F_3[\cdot \cdot \cdot ;q,x]$ series for the 
remote stretched triangles of graph (\ref{grp1r}) $l_4=k_4+j_1$ and 
$j_3=l_3+j_4$ or $k_2=k_1+l_1$ (when two couples of the touching angular 
momenta form the 4-cycles, or quadrangles), as well as expression 
(\ref{dpstb}) for $k_4=l_4+j_1$ and $l_3=j_3+j_4$ or $k_2=k_1+l_1$. 
Again, Minton's summation formula (\ref{sMint}) may be used in 
(\ref{dpstb}) for $j_4=j_3+l_3$ and the $_3F_2[\cdot \cdot \cdot ;q,x]$ 
series may be rearranged in the triple sums, which remain in Eq.\ 
(\ref{dpsta}) for $l_1=k_1+k_2$ or $l_3=j_3+j_4$ and in Eq.\ 
(\ref{dpstb}) for $j_3=l_3+j_4$ or $k_1=k_2+l_1$. Expansion (\ref{dvpja}) 
with inserted stretched $q$-$6j$ coefficients may be used for summation 
by means of (\ref{smDqa}) for $l_4=k_4+j_1$ and $l_1=k_1+k_2$, as well 
as expansion (\ref{dvpjb}) for summation by means of (\ref{smDqb}) for 
$k_4=l_4+j_1$ and $k_1=k_2+l_1$, after in the both cases the middle and 
the right columns of the second $q$-$6j$ coefficients are transposed.

In general, the triply stretched cases of the $q$-$12j$ coefficients of 
the first kind may always be expressed in terms of the double sums, 
applying some ``mirror reflection'' ($j\to -j-1$) operations to the above 
mentioned expressions, or using the above mentioned rearrangements. 
In particular, the double sums (compositions of two 
$_3F_2[\cdot \cdot \cdot ;q,x]$ series) related to the most generic 
Kamp\'e de F\'eriet \cite{K-F21} functions of the type $F_{1:1}^{0:3}$ 
(with 9 independent parameters) appear instead of (\ref{dpstb}) in the 
triply stretched cases with $k_4=l_4+j_1=k_3+l_3$ and $j_3=j_2+l_2$, or 
with $j_1=k_4-l_4=j_2-j_1$ and $k_3=k_2+l_2$.

Let us return to the generic expressions (\ref{dpsa}) and (\ref{dpsb}). 
In both formulas, for $j_2=0,\ l_1=j_1,\ l_2=j_3$, the summation 
parameters $z_1,z_2$ and $u$ are fixed, the two remaining separate sums 
(over $z_3,z_4$) correspond to the balanced $_4F_3[q,1]$ series, and the 
couples of the $q$-$6j$ coefficients appear straightforwardly, in 
accordance with the expansion (\ref{dvpja})--(\ref{dvpjb}) for fixed $x$. 
Furthermore, in both (\ref{dpsa}) and (\ref{dpsb}), the parameters 
$z_1,z_2$ and $u$ are also fixed for $j_1+j_2+k_1-k_2=0$, i.e., for the 
adjacent consecutive stretched triangles $k_2=k_1+l_1=k_1+j_1+j_2$ with 
the intermediate angular momentum corresponding to a crossbar of the 
M\"obius strip (\ref{grp1r}). Two remaining separate sums (over 
$z_3,z_4$) correspond to the balanced (Saalsch\"utzian) $_4F_3[q,1]$ 
[related to $q$-$6j$ coefficient of the type (\ref{f6jb})] and summable 
[see Eq.\ (\ref{sSbal})] $_3F_2[q,1]$ series. In this way, we derived the 
following relation:
\begin{eqnarray}
&&\left\{ \begin{array}{cccccccc}
j_1\! &  & \!j_2\! &  & \!j_3\! &  & \!j_4\! &  \\ 
& \!j_1\!+\!j_2\! &  & \!l_2\! &  & \!l_3\! &  & \!l_4\! \\ 
k_1\! &  & \!k_1\!+\!l_1\! &  & \!k_3 \! &  & \!k_4\! & 
\end{array}\right\} _{\!\!q}  \nonumber \\
&&\qquad =(-1)^{k_1+j_3-l_3-l_4+2k_3}\frac{q^{2j_1k_1+Z_{j_2j_3j_4}+
Z_{k_2k_3k_4}}\nabla [k_2k_3l_2]\nabla [k_1\!+\!j_1,j_4,k_4]}{%
\nabla [j_1k_4l_4]\nabla [j_2l_2j_3]\nabla [k_1j_4l_4]
\nabla [k_1\!+\!j_1,j_3,k_3]}  \nonumber \\
&&\qquad \quad \times \left( \frac{[2j_2]![2j_1]![2k_1]!}{[2k_2+1]!
[2l_1+1]}\right) ^{1/2}\left\{ \begin{array}{ccc}
j_3 & j_4 & l_3 \\ 
k_4 & k_3 & k_1\!+\!j_1
\end{array} \right\} _{\!\!q}.  \label{dvsta} 
\end{eqnarray}

For 4 mutual positions of the couples of the stretched triangles 
in the graph of the $q$-$12j$ coefficient of the first kind, there are 
22 different orientations of the total (maximal) angular momenta. In 
seven cases, the expressions include triple sums, twice they are 
proportional to the stretched $q$-$9j$ coefficients, and once to the 
$q$-$6j$ coefficient. In the remaining 12 cases, double sums may be 
obtained. Otherwise, for 3 mutual positions of the stretched triangles 
in the $q$-$12j$ coefficient of the second kind, only 9 different 
orientations of the total angular momenta are possible, and in 6 cases 
double sums appear, and in 3 cases the expressions are single sums (only 
once proportional to the $q$-$6j$ coefficient).

Finally, the  summation parameters  $z_1,z_3$ and $u$ in (\ref{dpsb}) are 
fixed for $j_1-j_3+k_1+k_3=0$, as well as for $x=j_3-k_3=j_1+k_1$ in 
expansion (\ref{dvpja}), in which cases 4 linearly dependent angular 
momenta appear (disconnected) on the Hamilton line of graph 
(\ref{grp1r}). In this case, two remaining separate sums (over $z_2,z_4$) 
are summable Saalsch\"utzian series $_3F_2[q,1]$ [see Eqs.\ (\ref{sSbal}) 
and (\ref{sSbbl})], and the following expression may be derived: 
\begin{eqnarray}
&&\left\{ \begin{array}{cccccccc}
j_1\! &  & \!j_2\! &  & \!j_1\!+\!k_1\!+\!k_3\! &  & \!j_4\! &  \\ 
& \!l_1\! &  & \!l_2\! &  & \!l_3\! &  & \!l_4\! \\ 
k_1\! &  & \!k_2\! &  & \!k_3 \! &  & \!k_4\! & 
\end{array}\right\} _{\!\!q}=(-1)^{j_2+j_4-l_1-k_1-j_1-l_4}  \nonumber \\
&&\qquad \quad \times \frac{q^{2j_1k_1+Z_{j_2j_3j_4}+Z_{k_2k_3k_4}}[2j_1]!
[2k_1]![2k_3]!\nabla [j_3j_4l_3]\nabla [j_3j_2l_2]}{[2j_3+1]!
\nabla [j_1j_2l_1]\nabla [k_1k_2l_1]\nabla [j_1k_4l_4]\nabla [k_1j_4l_4]
\nabla [k_3k_2l_2]\nabla [k_3k_4l_3]}.  \label{dvpors}
\end{eqnarray}
As for (\ref{dvosb}), the virtually stretched triangles are 
seen only in the $q$-$6j$ coefficients, which appear in expansion 
of the $q$-$12j$ coefficient of the first kind with extreme parameters.

\section{CONCLUDING REMARKS}

Using Dougall's summation formula \cite{Sl66} of the very well-poised 
$_4F_3(-1)$ series and its generalization for the $q$-factorial series 
we derived six new (independent) triple sum expressions for $9j$ 
coefficients of SU(2) and seven independent triple sum expressions for 
$q$-$9j$ coefficients of the quantum algebra $u_q(2)$. Rearrangement 
technique of the multiple sum expressions give several classes of double 
sum expressions for the stretched $9j$ coefficients of SU(2) and 
$u_q(2)$, related to the Kamp\'e de F\'eriet functions. Hence the new 
multiple basic hypergeometric series are introduced (cf. Ref.\ 
\onlinecite{DMT00}). Otherwise, Dougall's summation formula \cite{Sl66} 
of the very well-poised $_5F_4(1)$ series and the transformation formula 
\cite{Al97} of the very well-poised $_6F_5(-1)$ series, allowed us to 
eliminate the factorial sums weighted with factor $(2j+1)$ and the very 
well-poised series in the traditional expressions of the $12j$ 
coefficients of SU(2) in terms of $6j$ coefficients, as well as in the 
expressions of the $q$-$12j$ coefficients of $u_q(2)$ weighted with 
factor $[2j+1]$. Although the obtained generic expressions for the 
$q$-$12j$ coefficients of the first kind include fivefold sums and the 
generic expressions for the $q$-$12j$ coefficients of the second kind 
contain fourfold sums, the stretched and doubly stretched $q$-$12j$ 
coefficients of both become considerably simpler and present some new 
versions of the triple, double, or single basic hypergeometric series, 
but more unusual and complicated than appear in the manuals of the 
angular momentum theory. The single term expressions for the $q$-$12j$ 
coefficients of the both kinds also embrace the virtually stretched 
cases with four extreme linearly dependent disconnected angular momenta 
parameters appearing on the Hamilton lines of graphs (\ref{grp2r}) or 
(\ref{grp1r}). 

The symmetry properties and variety of expressions for the $12j$ and 
$q$-$12j$ coefficients of the both kinds may inspire new possibilities of 
rearrangement of the multiple and usual classical and basic 
hypergeometric series. Expressions for $3nj$ coefficients of SU(2) (with 
$n>4$) in terms of $6j$ coefficients also may be rearranged using, e.g., 
Watson's transformation formula of the very well-poised $_7F_6(1)$ series 
[see Eq.\ (2.5.1) of Ref.\ \onlinecite{GR90} or Eq.\ (6.10) of Ref.\ 
\onlinecite{LB94}]. One may also use  Eq.\ (5.3) of Ref.\ 
\onlinecite{Al97} for rearrangement of the very well-poised $_8F_7(-1)$ 
series, or Eqs.\ (A5) of Ref.\ \onlinecite{Al92} for rearrangement of the 
very well-poised $_9F_8(1)$ series. The number of sums in two last cases 
increases with elimination of the sums dependent on $(2j+1)$.

\appendix
\section{Generalization of Dougall's $_4F_3(-1)$ and $_5F_4(1)$
summation formulas for $\lowercase{q}$-factorial series}

We present here 3 summation formulas of the twisted very well-poised 
$q$-factorial series, as generalizations of Dougall's summation formula 
(2.3.4.8) of Ref.\ \onlinecite{Sl66} of the very well-poised $_4F_3(-1)$ 
series. In particular, equation
\begin{eqnarray}
&&\sum_{j}\frac{(-1)^{p_1+j+1}q^{j(j+1)-p_1(p_1+1)}[2j+1][j-p_1-1]!}{%
[p_1+j+1]![p_2-j]![p_2+j+1]![p_3-j]![p_3+j+1]!}  \nonumber \\
&&\qquad \quad =\frac{q^{-(p_1+p_2+1)(p_1+p_3+1)}}{[p_1+p_2+1]!
[p_1+p_3+1]![p_2+p_3+1]!},  \label{smDqa}
\end{eqnarray}
(valid when the summation parameter $j$ is restricted naturally by the 
non-negative integer values of the denominator factorial arguments)
was derived by Ali\v {s}auskas.\cite{Al97} 

Two other summation formulas,
\begin{eqnarray}
&&\sum_{j}\frac{q^{j(j+1)-p_1(p_1+1)}[2j+1][j-p_1-1]![j-p_3-1]!}{%
[p_1+j+1]![p_2-j]![p_2+j+1]![p_3+j+1]!}  \nonumber \\
&&\qquad \quad =\frac{q^{-(p_1+p_2+1)(p_1+p_3+1)}[-p_1-p_3-2]!}{%
[p_1+p_2+1]![p_2+p_3+1]!}  \label{smDqb}
\end{eqnarray}
and 
\begin{eqnarray}
&& \sum_{j}\frac{(-1)^{p_2-j}q^{j(j+1)}[2j+1][j-p_1-1]![j-p_3-1]!
[-p_3-j-2]!}{[p_1+j+1]![p_2-j]![p_2+j+1]!}  \nonumber \\
&& \qquad \quad =\frac{q^{p_1(p_1+1)-(p_1+p_2+1)(p_1+p_3+1)}}{[p_1+
p_2+1]!}[-p_1-p_3-2]![-p_2-p_3-2]!,  \label{smDqc}
\end{eqnarray}
correspond to the analytical continuation of (\ref{smDqa}).

Besides we present here two summation formulas of the very well-poised 
$q$-factorial series,
\begin{mathletters}
\begin{eqnarray}
&&\sum_{j}\frac{[2j+1]\,[j-p_1-1]![j-p_2-1]!}{[p_1+j+1]![p_2+j+1]!
[p_3-j]![p_3+j+1]![p_4-j]![p_4+j+1]!}  \nonumber \\
&&\qquad \quad =\frac{[-p_1-p_2-2]![p_1+p_2+p_3+p_4+2]!}{[p_1+p_4+1]!
[p_2+p_4+1]![p_3+p_4+1]![p_2+p_3+1]![p_1+p_3+1]!},  \label{p4s0}
\end{eqnarray}
and\cite{Al97} 
\begin{eqnarray}
&&\sum_{j}\frac{(-1)^{p_4-j}\,[2j+1]\,[j-p_1-1]![j-p_2-1]![j-p_3-1]!%
}{[p_1+j+1]![p_2+j+1]![p_3+j+1]![p_4-j]![p_4+j+1]!}  \nonumber \\
&&\qquad \quad =\frac{[-p_1-p_2-2]![-p_2-p_3-2]![-p_1-p_3-2]!}{[p_1+p_4
+1]![p_2+p_4+1]![p_3+p_4+1]![-p_1-p_2-p_3-p_4-3]!}.  \label{p4sa}
\end{eqnarray}
\end{mathletters}
[Cf.\ Dougall's summation theorem of special very well-poised 
hypergeometric series $_5F_4(1)$ as (2.3.4.5) of Ref.\ \onlinecite{Sl66} 
and special very well-poised basic hypergeometric series $_6\phi _5$ as 
(2.4.2) of Ref.\ \onlinecite{GR90}. Note that very well-poised basic 
hypergeometric series $_{2k}\phi _{2k-1}$ cannot be expressed in terms of 
$_{2k}F_{2k-1}[\cdot \cdot \cdot ;q,x]$ series (\ref{fbhs})].

The very well-po                                                            ised $_6F_5(-1)$ and $_7F_6(1)$ series appear in context 
of the Clebsch--Gordan and $6j$ coefficients of SU(2) as presented in 
Ref.\ \onlinecite{Al74} (see also Ref.\ \onlinecite{JB77}), as well as 
their $q$-analogs in the CG (cf.\ Ref.\ \onlinecite{A-NS96}) and $6j$ 
coefficients (cf.\ Ref.\ \onlinecite{LB94}) of $u_q(2)$. The summation 
formulas (\ref{smDqa}) and (\ref{smDqb}) may be obtained after 
cancelling some factorials in the numerators and denominators of 
(\ref{p5sra}) and (\ref{p5srb}), respectively, for $p_3=-p_1-2$.

\section{Chu--Vandermonde, Saalsch\"utzian, and Minton's summation 
formulas}

We present here the Chu--Vandermonde--Gauss--Heine summation formulas, 
\cite{JB77,AsST96,GR90}
\begin{mathletters}
\begin{equation}
\sum_s\frac{q^{s(a+b+c)}}{[s]![b-s]![c-s]![a+s]!}=\frac{q^{bc}
[a+b+c]!}{[b]![c]![a+b]![a+c]!},  \label{sCVa}
\end{equation}
\begin{equation}
\sum_s\frac{(-1)^sq^{s(b+c-a-1)}[a-s]!}{[s]![b-s]![c-s]!}=q^{bc}
\frac{[a-b]![a-c]!}{[b]![c]![a-b-c]!}  \label{sCVb}
\end{equation}
for $a\geq b,c$,
\begin{equation}
\sum_s\frac{(-1)^sq^{s(b-a+c-1)}[a-s]!}{[s]![b-s]![c-s]!}=(-1)^c
q^{bc}\frac{[a-c]![b-a+c-1]!}{[c]![b]![b-a-1]!}  \label{sCVc}
\end{equation}
for $b>a\geq c$, and
\begin{equation}
\sum_sq^{s(a+b-c+2)}\frac{[a-s]![b+s]!}{[s]![c-s]!}=q^{(b+1)c}
\frac{[a-c]![b]![a+b+1]!}{[c]![a+b-c+1]!},  \label{sCVd}
\end{equation}
\end{mathletters}
valid for finite $q$-factorial series and needed for rearrangements of 
Section IV.

Under condition $c+d=a+b+e$, we may use also the summation formulas of 
the balanced (Saalsch\"utzian) $_3F_2[q,1]$ series [cf.\ Refs.\ 
\onlinecite{Sl66,GR90}],
\begin{mathletters}
\begin{equation}
\sum_s\frac{(-1)^s[c+s]![d-s]!}{[s]![a-s]![b-s]![e+s+1]!}=\frac{[c]!
[d-a]![d-b]![c+d+1]!}{[a]![b]![a+e+1]![b+e+1]![e-c]!}  \label{sSbal}
\end{equation}
for $e-c\geq 0$ and 
\begin{equation}
\sum_s\frac{(-1)^s[c-s]![d-s]!}{[s]![a-s]![b-s]![e-s+1]!}=\frac{[c-a]!
[c-b]![d-a]![d-b]!}{[a]![b]![e-c]![e-d]![e+1]!}  \label{sSbbl}
\end{equation}
\end{mathletters}
for $e-c\geq 0$ and $e-d\geq 0$.

It is sometimes useful to implement the $q$-version of Minton's summation 
formula or its inverse [cf.\ Eq.\ (1.9.6) of Ref.\ \onlinecite{GR90}]:
\begin{mathletters}
\begin{equation}
\sum_s\frac{(-1)^sq^{s(n-\sum_{i=1}^rm_i)}}{[s]![n-s]![S+s]}\tprod_{j=1}^m
(b_j+s|q)_{m_{j}}=\frac{q^{S(\sum_{i=1}^rm_i-n)}}{(S|q)_{n+1}}
\prod_{j=1}^m(b_j-S|q)_{m_{j}},  \label{sMint}
\end{equation}
if $S+s\neq 0$ for $s=0,1,...,n$, and $n\geq \sum_{i=1}^rm_i$; 
\begin{eqnarray}
&&\sum_s\frac{(-1)^sq^{s(a+b-c-m)}[c-s]!}{[s]![a-s]![b-s]![S-s+1]}  
\tprod_{j=1}^m[A_{j}-s]  \nonumber \\
&&\qquad =(-1)^{a+b-c-m}q^{(S+1)(a+b-c-m)}\frac{[S-a]![S-b]!}{[S-c]!
[S+1]!}\prod_{j=1}^m[S-A_{j}+1],  \label{sMinta}
\end{eqnarray}
\end{mathletters}
which is valid if $S-s+1\neq 0$ for $s=0,1,...,\min(a,b)$ and 
$a+b-c-m\geq 0$. Note that the analytical continuation of the summation 
formulas (\ref{sMint}) and (\ref{sMinta}) of the alternating series to 
related series with the fixed sign of all terms is impossible.

\section{Rearrangement formulas of some double sums and Kamp\'e 
de F\'eriet functions}

Comparing the double finite $q$-factorial series that appear in the most 
symmetric expression (\ref{stdK}) for the stretched $q$-$9j$ 
coefficients with its counterparts in Eqs.\ (\ref{stbK}) and 
(\ref{stra})--(\ref{stre}), respectively, and using single 
nonnegative integers $a,b,c,d,m,n,e-a,f-b$ (as the parameters restricting 
summation intervals in different situations) instead of the triangular 
linear combinations of angular momenta (after changing some summation 
parameters), the following rearrangement formulas may be written:
\begin{mathletters}
\begin{eqnarray}
&&\sum_{s,z}(-1)^{s+z}\frac{q^{s(b+c+e-m-n+1)-z(a+d+f-m-n+1)}[c+s]![e-s]!
[d+z]![f-z]!}{[s]![a-s]![z]![b-z]![n-s-z]![m-a-b+s+z]!}  \label{dsKa} \\
&&\quad =\frac{[e-a]![a+c-m]![d]![d+\!f\!+\!1]!}{[a]![m+n-a-b]!}
\sum_{s,z}\frac{[c+s]![e-s]![f-z]!}{[s]![e-a-s]![z]![b-z]!}  \nonumber \\
&&\qquad \times \frac{(-1)^{n+s}q^{n(b+c+1)-s(m+1)-z(a+b+c+d-m+1)}}{%
[d+f+1-z]![n-s-z]![a+c-m-n+s+z]!}  \label{dsKb} \\
&&\quad =q^{a(c+e-n+1)+bn}\frac{[d]![d+f+1]![e-a]![c+e+1]!}{[m+n-a-b]!}  
\nonumber \\
&&\qquad \times \sum_{s,z}\frac{(-1)^{a+s+z}q^{-s(e-n+z)-z(b+d-m)}[a+c-s]!
[f-z]![m+n-s-z]!}{[s]![a-s]![m-s]![c+e+1-s]![z]![b-z]![n-z]![d+f+1-z]!}  
\label{tst0k} \\
&&\quad =q^{n(c+1)-b(b+d-m)-(b+e-a)(c+e-n+1)}\frac{[d]![d+f+1]![e-a]!
[c+e+1]!}{[m+n-a-b]!}  \nonumber \\
&&\qquad \times \sum_{s,z}\frac{(-1)^{a+e+n+s}[c+e-a-s]![c+e-m-s]!
[f-b+z]!}{[s]![e-a-s]![c+e+1-s]![b+c+e-m-n-s-z]!}  \nonumber \\
&&\qquad \times \frac{q^{s(b+e-n-z)+z(b+c+d+e-m+1)}}{[z]![b-z]!
[n-b+z]![d+f-b+1+z]!}  \label{tstbk} \\
&&\quad =q^{a(c+1)-n(a+d-m)}\frac{[e-a]![c+e+1]![f-b]!}{[b]![m+n-a-b]!}  
\sum_{s,z}\frac{(-1)^{s+z}q^{s(e-z)-z(m-a+1)}}{[s]![a-s]![c+e-a+1+s]!}  
\nonumber \\
&&\qquad \times \frac{[c+s]![a+b+d-m-s]![f-n+z]![d+n-z]!}{[a+b+d-m-s-z]!
[z]![n-z]![f-b-n+z]!}  \label{tstck} \\
&&\quad =q^{a(c+e-m-n)-b(d+f-n+1)}\frac{[f-b]![d+f+1]!}{[m+n-a-b]!}
\sum_{z,s}(-1)^{a+b+s+z}\frac{[e-a+s]!}{[s]![a-s]!}  \nonumber \\
&&\qquad \times \frac{[a+c-s]![b+d-z]![m+n-a-z]!q^{z(a+f-n-s)-
s(c+e-m-n+1)}}{[n-a+s]![z]![b-z]![d+f+1-z]![m-s-z]!}  \label{tstdk} \\
&&\quad =q^{a(a+c+d+1)-n(a+d+1)+(a+b-f)(b-m)}
\frac{[c]![d]![f-b]![d+f+1]![c+e+1]!}{[b]![m+n-a-b]!}  \nonumber \\
&&\qquad \times \sum_{s,z}\frac{(-1)^{b-f+n+z}q^{-s(a+b+c+d-m+z+1)-
z(m+n-a-b)}}{[s]![a-s]![a+d+f-m-n-s]![c+e-a+s+1]!}  \nonumber \\
&&\qquad \times \frac{[e-a+s]![b+z]![b+d-m+z]!}{[z]![f-b-z]!
[b-a-f+n+s+z]![b+d+z+1]!}.  \label{tstek}
\end{eqnarray}
\end{mathletters}
We have only single terms in (\ref{dsKa}), (\ref{dsKb}), and 
(\ref{tstck}) for $n=0$, as well as in (\ref{dsKa}) and (\ref{tstdk}) 
for $m=0$. The bizarre restrictions $b+c+e-m-n\geq 0$ for (\ref{dsKb}) 
and (\ref{tstbk}) and $a+d+f-m-n\geq 0$ for (\ref{tstck}) and 
(\ref{tstek}) also correspond to some triangular conditions with 
remaining double series summable for their limit values. Otherwise,
restrictions $c+e-m\geq 0$ in (\ref{tstbk}), $a+b+d-m\geq 0$ in 
(\ref{tstck}), or $d+f-m\geq 0$ in (\ref{tstek}) correspond to some sums 
of triangular conditions.

Using above derived expressions for the stretched $q$-$9j$ coefficients, 
we may write in the notations (\ref{sK}) the following rearrangement 
formulas for the $q$-generalizations of special Kamp\'e de F\'eriet 
functions $F_{1:1}^{0:3}$, $F_{1:1}^{0:3}$, and $F_{1:1}^{1:2}$:
\begin{mathletters}
\begin{eqnarray}
&&\quad  ^{+\!}F_{1:1}^{0:3}\left[ \begin{array}{c}
- \\ b_1+b_1^{\prime }
\end{array}:\begin{array}{c}
b_1,\;b_2,\;-m \\ d 
\end{array};\begin{array}{c}
b_1^{\prime },\;b_2^{\prime },\;-n \\ d^{\prime } 
\end{array};x_a,y_a;q\right]  \label{hst0k} \\[5pt]
&&\qquad ={^{-\!}F_{0:2}^{1:2}}\left[ \begin{array}{c}
1\!-\!b_1\!-\!b_1^{\prime }\!-\!m\!-\!n \\ -
\end{array}\!:\begin{array}{c}
-m-d+1,\;-m \\ \!1\!-\!b_1\!-\!m,1\!-\!b_2\!-\!m 
\end{array}\!;\begin{array}{c}
-d^{\prime }-n+1,\;-n \\ 
\!1\!-\!b_1^{\prime }\!-\!n,1\!-\!b_2^{\prime }\!-\!n 
\end{array};x_b,y_b;q\right]  \nonumber \\
&&\qquad \quad \times (-1)^{m+n}q^{m(b_2-b_1^{\prime }-d-m+1)
+n(b_2^{\prime }-b_1-d^{\prime }-m-n+1)}  
\nonumber \\
&&\qquad \quad \times \frac{(b_1|q)_m(b_2|q)_m(b_1^{\prime }|q)_n
(b_2^{\prime }|q)_n}{(d|q)_m(d^{\prime }|q)_n
(b_1+b_1^{\prime }|q)_{(m+n)}}  \label{hstak} \\[5pt]
&&\qquad ={^{-\!}F_{0:2}^{1:2}}\left[ \begin{array}{c}
b_2\!-\!b_1\!-\!b_1^{\prime }\!-\!n\!+\!1 \\ -
\end{array}\!:\begin{array}{c}
b_2-d+1,\;b_2 \\ \!b_2\!+\!m\!+\!1,b_2\!-\!b_1\!+\!1
\end{array}\!;\begin{array}{c}
1-d^{\prime }-n,\;-n \\ 
1\!-\!b_1^{\prime }\!-\!n,1\!-\!b_2^{\prime }\!-\!n
\end{array};x_c,y_c;q\right]  \nonumber \\
&&\qquad \quad \times (-1)^{b_1^{\prime }+b_2+d+m-1}q^{n(b_2^{\prime }
-d^{\prime })+b_2(d+m-b_2-1)+b_1^{\prime }(2b_1+b_1^{\prime }+2n-1)}  
\nonumber \\
&&\qquad \quad \times \frac{(b_2^{\prime }|q)_n(d-b_2|q)_{(-d-m)}}{%
(d^{\prime }|q)_n(m+1|q)_{(-d-m)}}\QATOPD[ ] {b_1-b_2-1}{%
-b_1^{\prime }-n}_{q}\QATOPD[ ] {-b_1-b_1^{\prime }}{-b_1}_{q}^{-1}  
\label{hstbk} \\[5pt]
&&\qquad ={^{+\!}F_{1:1}^{0:3}}\left[ \begin{array}{c}
- \\ d^{\prime }\!-\!b_2^{\prime }\!+\!b_1\!+\!b_1^{\prime }\!+\!n
\end{array}\!:\begin{array}{c}
\!b_1\!-\!b_2^{\prime }\!+\!d^{\prime }\!+\!n,b_2,-m \\ d
\end{array}\!;\begin{array}{c}
\!d^{\prime }\!-\!b_2^{\prime },b_1^{\prime },d^{\prime }\!+\!n \\ 
d^{\prime } 
\end{array};x_d,y_d;q\right]  \nonumber \\
&&\qquad \quad \times (-1)^{b_1^{\prime }}q^{-b_1^{\prime }
(b_2^{\prime }-d^{\prime }-n)}\frac{(d^{\prime }-b_2^{\prime }
+b_1+b_1^{\prime }+n|q)_{(-b_1^{\prime })}}{(-b_1+1|q)_{(-b_1^{\prime })}}  
\label{hst0ck} \\[5pt]
&&\qquad ={^{+\!}F_{1:1}^{0:3}}\left[ \begin{array}{c}
- \\ 1\!-\!b_1\!-\!m\!-\!n
\end{array}\!:\begin{array}{c}
\!b_1^{\prime },1\!-\!d\!-m,-m \\ b_2-d-m+1
\end{array};\begin{array}{c}
\!d^{\prime }\!-\!b_2^{\prime },1\!-\!b_1\!-\!b_1^{\prime }\!-m\!-\!n,-n \\ 
1-b_2^{\prime }-n
\end{array};x_e,y_e;q\right]  \nonumber \\
&&\qquad \quad \times q^{m(b_2-b_1^{\prime }-1)+n(b_1^{\prime }
+b_2^{\prime }-d^{\prime })}\frac{(1-b_1-m-n|q)_{(-b_1^{\prime })}
(d-b_2|q)_m(b_2^{\prime }|q)_n}{(-b_1+1|q)_{(-b_1^{\prime })}(d|q)_m
(d^{\prime }|q)_n}  \label{hst0dk} \\[5pt]
&&\qquad ={^{-\!}F_{0:2}^{1:2}}\left[ \begin{array}{c}
b_1 \\ -
\end{array}\!:\begin{array}{c}
d-b_2,\;-m \\ 1-b_1^{\prime }-m,d 
\end{array};\begin{array}{c}
b_2^{\prime },\;-n \\ 
b_2^{\prime }-d^{\prime }-n+1,b_1+b_1^{\prime }+m 
\end{array};x_f,y_f;q\right]  \nonumber \\
&&\qquad \quad \times q^{m(b_1-1)+nb_2^{\prime }}
\frac{(b_1^{\prime }|q)_m(d^{\prime }-b_2^{\prime }|q)_n}{%
(b_1+b_1^{\prime }|q)_m(d^{\prime }|q)_n}  \label{hstdk} \\[5pt]
&&\qquad ={^{+\!}F_{1:1}^{1:2}}\left[ \begin{array}{c}
b_1^{\prime } \\ 1-b_1-m-n
\end{array}\!:\begin{array}{c}
1-d-m,\;-m \\ b_2-d-m+1 
\end{array}\!;\begin{array}{c}
d^{\prime }-b_2^{\prime },\,-n \\ d^{\prime } 
\end{array};x_g,y_g;q\right]  \nonumber \\
&&\qquad \quad \times q^{m(b_2-b_1^{\prime })+b_1^{\prime }n}
\frac{(d-b_2|q)_m(-b_1-m-n+1|q)_{(-b_1^{\prime })}}{(d|q)_m
(-b_1+1|q)_{(-b_1^{\prime })}}  \label{hsKa} \\[5pt]
&&\qquad ={^{-\!}F_{1:1}^{1:2}}\left[ \begin{array}{c}
b_1^{\prime } \\ b_1+b_1^{\prime }-d+1
\end{array}\!:\begin{array}{c}
b_2-d+1,1-d-m \\ b_2-d-m+1 
\end{array}\!;\begin{array}{c}
b_2^{\prime },-n \\ d^{\prime } 
\end{array};x_h,y_h;q\right]  \nonumber \\
&&\qquad \quad \times q^{b_1^{\prime }(d-1)+mb_2}\frac{(d-b_2|q)_m
(d-b_1|q)_{(-b_1^{\prime })}}{(d|q)_m(-b_1+1|q)_{(-b_1^{\prime })}}  
\label{hsKb} \\[5pt]
&&\qquad ={^{+\!}F_{1:1}^{1:2}}\left[ \begin{array}{c}
b_2-b_1-b_1^{\prime }-n+1 \\ b_2-b_1^{\prime }-d-n+2
\end{array}:\begin{array}{c}
b_2-d+1,m+1 \\ b_2+m+1 
\end{array};\begin{array}{c}
1-d^{\prime }-n,-n \\ 
1-b_2^{\prime }-n 
\end{array};x_h,y_h;q\right]  \nonumber \\
&&\qquad \quad \times (-1)^{b_1^{\prime }+b_2+d+m-1}
q^{b_2(2b_1^{\prime }-b_2+d+m-1)-b_1^{\prime }(b_1^{\prime }-1)+n
(b_2^{\prime }-d^{\prime })+2(b_2-b_1-b_1^{\prime }-n+1)
(b_2-b_1^{\prime }-d+1)}  \nonumber \\
&&\qquad \quad \times \frac{(b_2^{\prime }|q)_n(d-b_2|q)_{(-d-m)}
}{(d^{\prime }|q)_n(m+1|q)_{(-d-m)}}\QATOPD[ ] {b_1-d}
{b_1\!+\!b_1^{\prime }\!-\!b_2\!+\!n\!-\!1}_{q}\QATOPD[ ] {-b_1-
b_1^{\prime }}{-b_1}_{q}^{-1}  \label{hsvKb} \\[5pt]
&&\qquad ={^{-\!}F_{1:1}^{0:3}}\left[ \begin{array}{c}
- \\ b_1+b_1^{\prime }\end{array}:\begin{array}{c}
b_1,\;b_2^{\prime },\;-n \\ 
b_2^{\prime }-d^{\prime }-n+1
\end{array};\begin{array}{c}
b_1^{\prime },\;b_2,\;-m \\ b_2-d-m+1 
\end{array};x_j,y_j;q\right]  \nonumber \\
&&\qquad \quad \times q^{mb_2+nb_2^{\prime }}\frac{(d-b_2|q)_m
(b_2^{\prime }-d^{\prime }-n+1|q)_n}{(d|q)_m(-n-d^{\prime }+1|q)_n}  
\label{hst0kb} \\[5pt]
&&\quad ={^{+\!}F_{0:2}^{1:2}}\left[ \begin{array}{c}
b_2\!-\!b_1\!-\!b_1^{\prime }\!-\!n\!+\!1 \\ -
\end{array}\!:\begin{array}{c}
d^{\prime }-b_2^{\prime },\,-n \\ 
\!d\!-\!b_1\!-\!n,1\!-\!b_2^{\prime }\!-\!n 
\end{array}\!; 
\begin{array}{c}
b_2-d+1,\,m+1 \\ \!b_2\!+\!m\!+\!1,b_2\!-\!b_1^{\prime }\!-\!d\!+\!2 
\end{array};x_k,y_k;q\right] \nonumber \\
&&\qquad \quad \times (-1)^{b_1^{\prime }+b_2+d+m+1}
q^{2(b_2-b_1-b_1^{\prime }-n+1)(b_2-b_1^{\prime }-d+1)-b_2(b_2-d-m+1)
+n(b_2^{\prime }-d^{\prime })+b_1^{\prime }(2b_2-b_1^{\prime }+1)}  
\nonumber \\
&&\qquad \quad \times \frac{(b_2^{\prime }|q)_n(d-b_2|q)_{(-d-m)}}{%
(d^{\prime }|q)_n(m+1|q)_{(-d-m)}}
\QATOPD[ ] {b_1+n-d}{b_2-b_1^{\prime }-d+1}_{q}
\QATOPD[ ] {-b_1-b_1^{\prime }}{-b_1}_{q}^{-1},  \label{hsvck}
\end{eqnarray}
\end{mathletters}
where 
\begin{eqnarray*}
x_a&=&q^{b_2-b_1^{\prime }-d-m+1}\quad {\rm and}\quad 
y_a=q^{b_2^{\prime }-b_1-d^{\prime }-n+1}, \\ 
x_b&=&q^{b_2-b_1^{\prime }-d-m-n+1}\quad {\rm and}\quad 
y_b=q^{b_2^{\prime }-b_1-d^{\prime }-m-n+1}, \\
x_c&=&q^{b_2-b_1^{\prime }-d-m-n+1}\quad {\rm and}\quad y_c=q^{b_2-b_1+
b_2^{\prime }-d^{\prime }-n+1}, \\
x_d&=&q^{b_2-b_1^{\prime }-d-m+1}\quad {\rm and}\quad y_d=q^{-b_1+1}, \\
x_e&=&q^{b_1+b_1^{\prime }-b_2+n}\quad {\rm and}\quad y_e=q^{d^{\prime }-
b_1^{\prime }}, \\
x_f&=&q^{b_1+b_1^{\prime }-b_2}\quad {\rm and}\quad y_f=q^{d^{\prime }-
b_1^{\prime }-m}, \\
x_g&=&q^{b_1+b_1^{\prime }-b_2+n}\quad {\rm and}\quad y_g=q^{b_1+
b_1^{\prime }-b_2^{\prime }+m}, \\
x_h&=&q^{-b_1+1}\quad {\rm and}\quad y_h=q^{d-d^{\prime }-b_1+
b_2^{\prime }-n}, \\
x_j&=&q^{d^{\prime }-b_1^{\prime }}\quad {\rm and}\quad y_j=q^{d-b_1}, \\
{\rm and}\quad x_k&=&q^{d^{\prime }-d+b_2-b_1^{\prime }+1}\quad {\rm and} 
\quad y_k=q^{-n-b_1+1}, 
\end{eqnarray*}
respectively, where only parameters $m,n,-b_1,-b_1^{\prime }$ are 
apparently correlated with some triangular conditions. Special Kamp\'e de 
F\'eriet functions (\ref{hst0k}) and (\ref{hstak}) correspond, 
respectively, to the inverse and direct sums in (\ref{stra}), when 
function (\ref{hstbk}) corresponds to the direct sum in (\ref{strb}), 
function (\ref{hst0ck}) corresponds to the inverse sum in (\ref{strc}), 
and functions (\ref{hst0dk}) and (\ref{hstdk}) correspond, respectively, 
to the inverse and direct sums in (\ref{strd}). Further, functions 
(\ref{hsKa}), (\ref{hsKb}), and (\ref{hsvKb}) correspond, respectively, 
to the sums that appeared in (\ref{stdK}) and (\ref{stbK}), as well 
as in (\ref{dsKa}) and (\ref{dsKb}). The two last functions 
(\ref{hst0kb}) and (\ref{hsvck}) are derived from (\ref{hst0k}) and the 
direct sum in (\ref{strc}), respectively, after using the symmetry of 
$^{+}F_{1:1}^{1:2}$ function in (\ref{hsKa}) with fixed $b_1$ and 
$b_1^{\prime }$ under interchange of two sets,
\[
b_2,m,d,d^{\prime }\quad {\rm and}\quad b_2^{\prime },n,b_2^{\prime }
-d^{\prime }-n+1,b_2-d-m+1,
\]
together with transition to $^{-}F_{1:1}^{1:2}$ and $q\rightarrow q^{-1}$.

Finiteness of the Kamp\'e de F\'eriet series (\ref{hst0k}) is ensured 
either by the non-negative integer values of $m$ and $n$, or by the 
non-positive integer values of $b_1$ and $b_1^{\prime }$, or by some their 
couples ($m$ and $-b_1^{\prime }$, or $-b_1$ and $n$). The both 
summation parameters are also restricted by the non-negative integer 
values of $m$ and $n$ in series (\ref{hstak}), (\ref{hst0dk})--%
(\ref{hsKa}), and (\ref{hst0kb}), as well as by the non-negative values 
of $m$ and $-b_1^{\prime }$ in series (\ref{hst0k}) and (\ref{hst0ck}), or
by the non-negative values of $n$ and $d-b_2-1$ in series (\ref{hstbk}), 
(\ref{hsKb}), (\ref{hsvKb}), and (\ref{hsvck}). Furthermore, the 
parameter $b_1^{\prime }$ with the non-positive integer values restricts 
the double series in (\ref{hsKb}) and (\ref{hsvKb}), as well as separate 
series in (\ref{hst0k}), (\ref{hst0dk}), and (\ref{hst0kb}). Series 
(\ref{hstdk}) are finite for the non-positive integer values of single 
parameter $b_1$, as well as (\ref{hstbk}), (\ref{hsvKb}), and 
(\ref{hsvck}) for the non-positive integer values of 
$b_2\!-\!b_1\!-\!b_1^{\prime }\!-\!n\!+\!1$. Hence, special Kamp\'e de 
F\'eriet functions (\ref{hst0k})--(\ref{hsvck}) are summable for $b_1=0$, 
or $b_1^{\prime }=0$, or $b_2\!-\!b_1\!-\!b_1^{\prime }\!-\!n\!+\!1=0$ 
(i.e., for $c=b_1-b_1^{\prime }=b_2-n+1$) and, taking into account the 
symmetry of (\ref{hst0k}) with respect to the interchange of two sets, 
$b_1,b_2,-m;d$ and $b_1^{\prime },b_2^{\prime },-n;d^{\prime }$, for 
$b_2^{\prime }\!-\!b_1^{\prime }\!-\!b_1\!-\!m\!+\!1=0$ (i.e., for 
$c=b_1-b_1^{\prime }=b_2^{\prime }-m+1$). The subscripts of the 
$q$-Pochhammer symbols in the proportionality coefficients are accepted 
as non-negative integers, when they perform the restricting role or 
correspond to definite non-negative linear combinations of $9j$ 
parameters. Otherwise, for the negative integer subscripts $(-n)$ the 
following substitution may be used:
\[
\frac{(\alpha +n|q)_{(-n)}}{(\beta +n|q)_{(-n)}}\rightarrow 
\frac{(\beta|q)_n}{(\alpha |q)_n}.
\]
Hence in the $q=1$ case up to 5 or 6 parameters may be complex in the 
rearrangement formulas (\ref{hst0k})--(\ref{hst0kb}) of special Kamp\'e 
de F\'eriet series, with exception of (\ref{hstbk}), (\ref{hsvKb}), and 
(\ref{hsvck}). In these three cases, which ensure the summability of the 
remaining series for $b_1-b_1^{\prime }=b_2^{\prime }-m+1$, only 
$b_2^{\prime }$ and $d^{\prime }$ definitely may be taken the complex 
numbers. Extension problem to infinite series is open, since it is 
impossible to ensure the non-negative values of the all denominator 
arguments of $F_{1:1}^{0:3}$ series (\ref{hst0k})--(\ref{hsvck}) in the 
standard situation of the SU(2) stretched $9j$ coefficients, with 
exception of special Kamp\'e de F\'eriet series (with rather complicated 
parameters and proportionality coefficient), which could be written 
instead of (\ref{tstek}). 

Using the substitution (\ref{sasqf}) and different strategy for each 
mutual relation, we may transform the double finite series (\ref{hst0k}), 
(\ref{hst0ck}), (\ref{hst0dk}), (\ref{hsKa})--(\ref{hsvKb}), and 
(\ref{hsvck}) into standard functions $\Phi _{C:D}^{A:B}$, with partial 
cancelling of the $q$-phases of the proportionality coefficients. 
Preliminary in these situations only $q^{b_2^{\prime }}$ and 
$q^{d^{\prime }}$ can always be replaced by the complex numbers. 
Analogically, the double finite series (\ref{hstak}), (\ref{hstbk}), and 
(\ref{hstdk})--(\ref{hst0kb})  may be transformed into standard functions 
$\Phi_{C:D}^{A:B}$ after substituting $\mp \to \pm $ in the superscripts 
of $^{\mp }F_{C:D}^{A:B}$ series and $q^{-1}\to q$ in the corresponding 
$q$-phases.

Note, that the summable Kamp\'e de F\'eriet series $F_{1:1}^{0:3}$ (that 
appeared in Refs.\ \onlinecite{V-JPS-R94,PV-J96}) cannot be embedded 
into above presented versions of $F_{1:1}^{0:3}$ series, or be derived 
from the expressions of the stretched $9j$ coefficients given in Sec.~IV. 
Actually, expansion (\ref{trxa}) does not simplify under condition 
$a+b-e=0$, but the $^{\pm }F_{1:1}^{0:3}$ series \cite{V-JPS-R94,PV-J96} 
appear from expression (\ref{trsa}) in the doubly stretched case with 
$a+b-e=0$ and $g=k+h$, when $9j$ coefficients are proportional to the 
Clebsch--Gordan coefficients. In this particular case of expression 
(\ref{trsa}), we may also identify quintuplet of factorials under the 
summation sign in the numerator and denominator and reexpand it using 
the Chu--Vandermonde summation formulas given in Appendix B. As result of 
two alternative summations we obtain a $_3F_2[\cdot \cdot \cdot ;q,x]$ 
series, which is completely summable  for $k=b+d$.

Special cases of $^{\pm }F_{0:2}^{1:2}$ functions should be mentioned in 
context of the double sums (with 7 independent parameters) that appear in 
the extreme $u_q(3)$ canonical seed isofactors \cite{Al97,AlD98} and as
definite matrix elements of the $u_q(3)$ algebra (see Section 5 of Ref.\ 
\onlinecite{Al97}) and are related to some $q$-factorial series 
resembling the very well-poised $_9\phi _8$ basic hypergeometric series 
[which for $q=1$ are equivalent to the very well-poised $_8F_7(-1)$ 
classical hypergeometric series]. Further, the extreme denominator 
(normalization) functions of the $u_q(3)$ and SU(3) canonical tensor 
operators (with 5 independent parameters) may be expressed in terms of 
$^{\pm }F_{1:2}^{1:3}$ functions [cf.\ Eqs.\ (5.9c) and (5.9d) of Ref.\ 
\onlinecite{Al97}, or Eqs.\ (3.7) and (3.14) of Ref.\ \onlinecite{Al96}], 
or in terms of $^{\pm }F_{2:1}^{2:2}$ functions [see Eq.\ (2.8) and 
Section II of Ref.\ \onlinecite{Al99}, taking into account the definite 
controversies of the $q$-extension from the classical SU(3) case]. 
Besides, the summation possibilities for these special Kamp\'e de 
F\'eriet functions are elementary.

\section{Clebsch--Gordan coefficients of SU(2) and $u_q(2)$ and twisted 
very well-poised series}

The very well-poised $_6F_5(-1)$ and $_7F_6(1)$ series appear 
in context of the Clebsch--Gordan and $6j$ coefficients of SU(2) as 
presented in Ref.\ \onlinecite{Al74} (see also Ref.\ \onlinecite{JB77}), 
as well as their $q$-analogs in the CG (cf.\ Ref.\ \onlinecite{A-NS96}, 
where the dual Hahn $q$-polynomials are considered) and $6j$ coefficients 
(cf.\ Ref.\ \onlinecite{LB94}) of $u_q(2)$. 

We deduce here a new expression for the Clebsch--Gordan coefficients of 
SU(2) and $u_q(2)$ directly from the recoupling relation:
\begin{eqnarray}
&&\left[ \begin{array}{ccc}
(j_2+m_2)/2 & (j_2-m_2)/2 & j_2 \\ 
(j_2+m_2)/2 & (m_2-j_2)/2 & m_2
\end{array} \right] _{\!q} \left[ \begin{array}{ccc}
j_1 & j_2 & j \\ 
m_1 & m_2 & m
\end{array} \right] _{\!q}  \nonumber \\
&& \quad =\sum _x (-1)^{j_1+j_2+j}([2x+1][2j_2+1])^{1/2}
\left\{ \begin{array}{ccc}
(j_2+m_2)/2 & (j_2-m_2)/2 & j_2 \\ 
j & j_1 & x
\end{array} \right\} _{\!q}  \nonumber \\
&& \qquad \times \left[ \begin{array}{ccc}
j_1 & (j_2+m_2)/2 & x \\ 
m_1 & (j_2+m_2)/2 & m^{\prime }
\end{array} \right] _{\!q} \left[ \begin{array}{ccc}
x & (j_2-m_2)/2 & j \\ 
m^{\prime } & (m_2-j_2)/2 & m
\end{array} \right] _{\!q},  \label{nfiCG}
\end{eqnarray}
where $m^{\prime }=m_1+\frac 12(j_2+m_2)=m+\frac 12(j_2-m_2)$. Inserting 
the stretched $6j$ and extreme CG coefficients expressed without sums, 
we obtain the following expression,
\begin{eqnarray}
\left[ \begin{array}{ccc}
j_1 & j_2 & j \\ 
m_1 & m_2 & m
\end{array} \right] _{\!q}&=&\nabla [j_2j_1j]\left( \frac{[2j+1]
[j_2+m_2]![j_2-m_2]![j_1-m_1]![j-m]!}{[j_1+m_1]![j+m]!}\right) ^{1/2}  
\nonumber \\
&& \times q^{\{j_2(j_2+1)-j_1(j_1+1)-j(j+1)\}/2-m_1j_2-
(j_2+m_2)(j_2+m_2+2)/4}  \nonumber \\
&& \times \sum _x \frac{(-1)^{j_1+(j_2+m_2)/2-x}q^{x(x+1)}[2x+1]
[x+m^{\prime }]!}{\nabla ^2[\frac 12(j_2+m_2),j_1,x]\nabla ^2
[\frac 12(j_2-m_2),j,x][x-m^{\prime }]!},  \label{nexCG}
\end{eqnarray}
where the right-hand side is related to the left-hand side  of Eq.\ 
(\ref{p5sra}) with parameters
\[
\begin{array}{c}
p_1=\frac 12(j_2-m_2)-j-1,\;\;p_2=-m^{\prime }-1,\;\;
p_3=\frac 12(j_2+m_2)-j_1-1, \\ 
p_4=\frac 12(j_2-m_2)+j,\;\;p_5=\frac 12(j_2+m_2)+j_1.
\end{array}
\]
Expression (\ref{nexCG}) is invariant under 12 relations of the Regge 
symmetry, corresponding to the permutations in the sets $p_1,p_2,p_3$ or 
$p_4,p_5$. After expressing the CG coefficient of $u_q(2)$ by means of 
Eq.\ (5.17) of Ref.\ \onlinecite{STK91a} [which after some cyclic 
permutation, is for $q=1$, related to Eq.\ (13.1c) of Ref.\ 
\onlinecite{JB77}], and using the symmetry relation (4.13) of Ref.\ 
\onlinecite{STK91b} (which allows one to interchange the parameters 
$j_2,m_2$ and $j,-m$ in the CG coefficients), we derive our Eq.\ 
(\ref{p5sra}) straightforwardly. The remaining very well-poised series 
with different numerator and denominator distributions of $q$-factorial 
arguments [e.g., the non-alternating left-hand side  of (\ref{p5srb}) 
or the non-alternating right-hand side of  Eq.\ (5.3) of Ref.\ 
\onlinecite{Al97} with $p_3=-p_2-2$, and their other analytical 
continuations] are not related to the Clebsch--Gordan coefficients of 
$u_q(2)$, although sometimes they may be related to the CG coefficients 
of $u_q(1,1)$.


\begin{thebibliography}{99}
\bibitem{Al99}  S. Ali\v {s}auskas, J.\ Math.\ Phys.\ {\bf 40}, 5939 
(1999). Note, that the term ``$+ks^{\prime \prime }+$'' in the 
denominator of the last row of Eq.\ (1.8) should be written as 
``$+k\!+\!s^{\prime \prime }+$'' and the all terms in the first row of 
the r.h.s.\ of Eq.\ (1.12) should be in exponent of $q$.
\end{thebibliography}
\end{document}